\documentclass[10pt,twoside]{report}
\usepackage{amsfonts,amsmath,amssymb,amstext,amsthm}
\usepackage{mathpazo}
\usepackage[all]{xy}
\usepackage{makeidx}
\usepackage{nomencl}
\usepackage{graphicx,subfigure,color}
\usepackage{hyperref}
\usepackage[paper = letterpaper, left = 3.295cm, right = 3.295cm, headsep = 8mm, footskip = 7mm,
top = 30mm, bottom = 22.5mm, headheight = 2cm]{geometry}
\newtheorem{theorem}{Theorem}[chapter]
\newtheorem{lemma}[theorem]{Lemma}
\newtheorem{proposition}[theorem]{Proposition}
\newtheorem{corollary}[theorem]{Corollary}
\newtheorem{theorem*}{Theorem}[]
\theoremstyle{definition}
\newtheorem{remark}[theorem]{Remark}

\newtheorem{example}[theorem]{Example}
\newtheorem{definition}[theorem]{Definition}
\newtheorem*{notation}{Notation}
\renewcommand{\baselinestretch}{1.1}
\usepackage{fancyhdr}
\usepackage[clearempty]{titlesec}
\usepackage{float}
\restylefloat{figure}
\author{Kristina Frantzen}
\title{K3-surfaces with special symmetry}
\date{Oktober 2008}
\newcommand*{\titleKF}{\begingroup
 \vspace*{2cm}
 \centering
  {\huge K3-surfaces with special symmetry}\par
  \vspace*{3cm}
  {\LARGE {Dissertation} }\par
  \vspace{3cm}
 {\large zur Erlangung des\\\vspace{0.2cm} Doktorgrades der Naturwissenschaften\\\vspace{0.2cm} an der Fakult\"at f\"ur Mathematik \\\vspace{0.2cm} der Ruhr-Universit\"at Bochum}\par
 \vspace{3cm}
 {\large vorgelegt von \Large \\\vspace{1.2cm} \LARGE Kristina Frantzen \\ \vspace{1.2cm}\large im Oktober 2008}\par
 \vspace*{2cm}
\endgroup}

\begin{document}
%
%
\thispagestyle{empty}
\enlargethispage{2cm}
\titleKF
\newpage
\cleardoublepage
\pagestyle{fancy}
\fancyhf{}
\fancyhead[EL]{\thepage}
\fancyhead[ER]{\emph{Contents}}
\fancyhead[OR]{\thepage}
\fancyhead[OL]{\emph{Contents}}
%
%
\tableofcontents
\cleardoublepage
\thispagestyle{empty}
\newpage
\pagestyle{fancy}
\fancyhf{}
\fancyhead[EL]{\thepage}
\fancyhead[ER]{\emph{Introduction}}
\fancyhead[OR]{\thepage}
\fancyhead[OL]{\emph{Introduction}}
%
%
\chapter*{Introduction} 
\addcontentsline{toc}{chapter}{Introduction}
K3-surfaces are special two-dimensional holomorphic symplectic
manifolds. They come equip\-ped with a symplectic form $\omega$, which
is unique up to a scalar factor, and their symmetries are naturally
partitioned into symplectic and nonsymplectic transformations. An
important class of K3-surfaces consists of those possessing an
\emph{antisymplectic involution}, i.e., a holomorphic involution $\sigma$ such that
$\sigma^* \omega = - \omega$.

K3-surfaces with antisymplectic involution occur classically as branched double covers of the projective plane, or more generally of Del Pezzo surfaces. This construction is a prominent source of examples and plays a significant role in the classification of log Del Pezzo surfaces of index two (see the works of Alexeev and Nikulin e.g.\,in \cite{AlexNikulin} and the classification by Nakayama \cite{Nakayama}).
Moduli spaces of K3-surfaces with antisymplectic involution are
studied by Yoshikawa in \cite{yoshikawaInvent},
\cite{yoshikawaPreprint}, and lead to new developments in the area of automorphic forms.

In this monograph we study K3-surfaces with antisymplectic involution from the point of view of symmetry. On a K3-surface $X$ with antisymplectic involution it is natural the consider those holomorphic symmetries of $X$ compatible with the given structure $(X,\omega, \sigma)$. These are symplectic automorphisms of $X$ commuting with $\sigma$. 

Given a finite group $G$ one wishes to understand if it can act in the
above fashion on a K3-surface $X$ with antisymplectic involution
$\sigma$. If this is the case, i.e., if there exists a holomorphic
action of $G$ on $X$ such that $g^* \omega = \omega$ and $g \circ
\sigma = \sigma \circ g$ for all $g \in G$, then the structure of $G$
can yield strong constraints on the geometry of $X$. More precisely,
if the group $G$ has rich structure or large order, it is possible to
obtain a precise description of $X$. This can be considered the guiding
classification problem of this monograph. 

In Chapter \ref{chapterlarge} we derive a classification of K3-surfaces with antisymplectic involution centralized by a group of symplectic automorphisms of order greater than or equal to 96. We prove (cf. Theorem \ref{roughclassi}):
\begin{theorem*} \label{1}
Let $X$ be a K3-surface with a symplectic action of $G$ centralized by an
antisymplectic involution $\sigma$ such that
$\mathrm{Fix}(\sigma)\neq \emptyset$. If $|G|>96$, then $X/\sigma$ is a
Del Pezzo surface and $\mathrm{Fix}(\sigma)$ is a smooth connected curve $C$ with $ g(C)\geq
3$.
\end{theorem*}
By a theorem due to Mukai \cite{mukai} finite groups of symplectic
transformations on K3-surfaces are characterized by the existence of a
certain embedding into a particular Mathieu group and are 
subgroups of eleven specified finite groups of maximal symplectic symmetry.
This result naturally limits our considerations and has
led us to consider the above classification problem for a group $G$ from this list of eleven \emph{Mukai groups}.

Theorem \ref{1} above can be refined to obtain a complete classification of K3-surfaces with a symplectic action of a Mukai group centralized by an antisymplectic involution with fixed points (cf. Theorem \ref{thm mukai times invol}). 
\begin{theorem*}
Let $G$ be a Mukai group acting on a K3-surface $X$ by symplectic transformations. Let $\sigma $ be an antisymplectic involution on $X$ centralizing $G$ with $\mathrm{Fix}_X(\sigma) \neq \emptyset$. Then the pair $(X,G)$ can be found in Table \ref{Mukai times invol}.
\end{theorem*}
In addition to a number of examples presented by Mukai we find new examples of K3-surfaces with maximal symplectic symmetry as equivariant double covers of Del Pezzo surfaces.

It should be emphasized that the description of K3-surfaces
with given symmetry does however not necessary rely on the size of the group or its maximality and a classification can also be obtained for rather small subgroups of the Mukai groups.
In order to illustrate that the approach does rather depend on the structure of the group, we prove a classification of K3-surfaces with a symplectic action of the group $C_3 \ltimes C_7$ centralized by an antisymplectic involution in Chapter \ref{chapterC3C7}. The surfaces with this given symmetry are characterized as double covers of $\mathbb P_2$ branched along invariant sextics in a precisely described one-dimensional family $\mathcal M$ (Theorem \ref{mainthmc3c7}). 
\begin{theorem*}
The K3-surfaces with a symplectic action of $G = C_3 \ltimes C_7$ centralized by an antisymplectic involution $\sigma$ are parametrized by the space $\mathcal M$ of equivalence classes of sextic branch curves in $\mathbb P_2$.
\end{theorem*}
The group $C_3 \ltimes C_7$ is a subgroup of the simple group $L_2(7)$
of order 168 which is among the Mukai groups. The actions of $L_2(7)$ on
K3-surfaces have been studied by Oguiso and Zhang \cite{OZ168} in an a
priori more general setup. Namely, they consider finite groups containing $L_2(7)$ as a proper subgroup and obtain lattice theoretic classification results using the Torelli theorem. Since a finite group containing $L_2(7)$ as a proper subgroup posseses, in the cases considered, an antisymplectic involution centralizing $L_2(7)$, we can apply Theorem \ref{thm mukai times invol} and improve the existing result (cf. Theorem \ref{improve OZ}).

All classification results summarized above are proved by applying the following general strategy.

The quotient of a K3-surface by an antisymplectic involution $\sigma$ with fixed points centralized by a finite group $G$ is a rational $G$-surface $Y$. We apply an equivariant version of the minimal model program respecting finite symmetry groups to the surface $Y$. Chapter \ref{chapter mmp} is dedicated to a detailed derivation of this method, a brief outline of which can also be found in the book of Koll\'ar and Mori (\cite{kollarmori} Example 2.18, see also Section 2.3  in \cite{Mori}). 
In the setup of rational surfaces it leads to the well-known classification of $G$-minimal rational surfaces (\cite{maninminimal}, \cite{isk}). 

Equivariant Mori reduction and the theory of $G$-minimal models have applications in various different context and can also be generalized to higher dimensions.
Initiated by Bayle and Beauville in \cite{Bayle}, the methods have been employed in the classification of subgroups of the Cremona group $\mathrm{Bir}(\mathbb P_2)$ of the plane for example by Beauville and Blanc (\cite{Beauville}, \cite{BeauBlancPrime}, \cite{PhDBlanc}), \cite{Blanc1}, etc.), de Fernex \cite{fernex}, Dolgachev and Iskovskikh \cite{DolgIsk}, and Zhang \cite{ZhangRational}.

The equivariant minimal model $Y_\mathrm{min}$ of $Y$ is obtained from
$Y$ by a finite number of blow-downs of (-1)-curves. Since individual
(-1)-curves are not necessarily invariant, each reduction step blows
down a number of disjoint (-1)-curves. The surface $Y_\mathrm{min}$ is, in all cases considered, a Del Pezzo surface. 

Using detailed knowledge of the equivariant reduction map $Y \to Y_\mathrm{min}$, the shape of the invariant set $\mathrm{Fix}_X(\sigma)$, and the equivariant geometry of Del Pezzo surfaces, we classify $Y$, $Y_\mathrm{min}$ and $\mathrm{Fix}_X(\sigma)$ and can describe $X$ as an equivariant double cover of a possibly blown-up Del Pezzo surface. Besides the book of Manin, \cite{manin}, our analysis relies, to a certain extend, on Dolgachev's discussion of automorphism groups of Del Pezzo surfaces in \cite{dolgachev}, Chapter 10.

In addition to classification, this method yields a multitude of new
examples of K3-surfaces with given symmetry and a more geometric
understanding of existing examples. It should be remarked that a
number of these arise when the reduction $Y \to Y_{\mathrm{min}}$ is nontrivial.

In the last two chapters we present two different generalizations of our classification strategy for K3-surfaces with antisymplectic involution.

One of our starting points has been the 
study of K3-surfaces with $L_2(7)$-symmetry by Oguiso and Zhang mentioned above. Apart from a classification result for K3-surfaces with an action the group $L_2(7)\times C_4$, they also show that there does not exist a K3-surface with an action of a the group $L_2(7)\times C_3$. We give an independent proof of this result in Chapter \ref{chapter non exist}.
Assuming the existence of such a surface and following the strategy above,
we consider the quotient by the nonsymplectic action of $C_3$ and apply the equivariant minimal model program to its desingularization. Combining this with additional geometric consideration we reach a contradiction.

In the last chapter we consider K3-surfaces $X$ with an action of a
finite group $\tilde G$ which contains an antisymplectic involution
$\sigma$ but is not of the form $\tilde G_\mathrm{symp} \times \langle
\sigma \rangle$. Since the action of $\tilde G_\mathrm{symp}$ does not
descend to the quotient $X/\sigma$ we need to restrict our
considerations to the centralizer of $\sigma$ inside $\tilde G$. This
strategy is exemplified for a finite group $\tilde A_6$ characterized
by the short exact sequence $\{\mathrm{id}\} \to A_6 \to \tilde A_6
\to C_4 \to \{\mathrm{id}\}$. In analogy to the $L_2(7)$-case, the
action of $\tilde A_6$ on K3-surfaces has been studied by Keum,
Oguiso, and Zhang (\cite{KOZLeech}, \cite{KOZExten}), and a
characterization of $X$ using lattice theory and the Torelli theorem
has been derived. Since the existing realization of $X$ does however
not reveal its equivariant geometry, we reconsider the problem and,
though lacking the ultimate classification, find families of
K3-surfaces with $D_{16}$-symmetry, in which the $\tilde A_6$-surface
is to be found, as branched double covers. 
These families are of independent interest and should be studied further. In particular, 
it remains to find criteria to identify the $\tilde A_6$-surface inside these families. Possible approaches are outlined at the end of Chapter \ref{chapterA6}.

Since none of our results depends on the Torelli theorem, our approach to the classification problem allows generalization to fields of appropriate positive characteristic. This possible direction of further research was proposed to the author by Prof.\,Keiji Oguiso. Another potential further development would be the adaptation of the methods involved in the present work to related questions in higher dimensions.
\vfill
\begin{small}
\textit{Research supported by Studienstiftung des deutschen Volkes and Deutsche Forschungsgemeinschaft}.
\end{small}
%
%
%
%
\chapter{Finite group actions on K3-surfaces}
%
%
\pagestyle{fancy}
\fancyhf{}
\fancyhead[EL]{\thepage}
\fancyhead[ER]{\emph{\nouppercase{\leftmark}}}
\fancyhead[OR]{\thepage}
\fancyhead[OL]{\emph{\nouppercase{\rightmark}}}
This chapter is devoted to a brief introduction to finite groups actions on K3-surfaces and presents a number of basic, well-known results: 
We consider quotients of K3-surfaces by finite groups of symplectic or nonsymplectic automorphisms. It is shown that the quotient of a K3-surface by a finite group of symplectic automorphisms is a K3-surface, whereas the quotient by a finite group containing nonsymplectic transformations is either rational or an Enriques surface. Our attention concerning nonsymplectic automorphisms is then focussed on antisymplectic involutions and the description of their fixed point set.
The chapter concludes with Mukai's classification of finite groups of symplectic automorphisms on K3-surfaces and a discussion of basic examples. 
\section{Basic notation and definitions}
Let $X$ be a n-dimensional compact complex manifold. We denote by $\mathcal O_X$ the sheaf of holomorphic functions on $X$ \nomenclature{$\mathcal O_X$}{the sheaf of holomorphic functions on $X$} and by $\mathcal K_X$ its canonical line bundle\nomenclature{$\mathcal K_X$}{the canonical line bundle of $X$}. 
The \emph{ $i^{\text{th}}$ Betti number of $X$} \index{Betti number}\nomenclature{$b_i(X)$}{the $i^{\text{th}}$ Betti number of $X$} is the rank of the free part of $H_i(X)$ and denoted by $b_i(X)$. 

A \emph{surface}\index{surface} is a compact connected complex manifold of complex dimension two. A \emph{curve} \index{curve} on a surface $X$ is an irreducible 1-dimensional closed subspace of $X$. The (arithmetic) genus of a curve $C$ is denoted by $g(C)$ \nomenclature{$g(C)$}{the (arithmetic) genus of a curve $C$}.
\begin{definition}
A \emph{K3-surface} \index{K3-surface} is a surface $X$ with trivial canonical bundle $\mathcal K_X$ and $b_1(X)=0$.
\end{definition}
Note that a K3-surface is equivalently characterized if the condition $b_1(X)=0$ is replaced by $q(X) = \mathrm{dim}_\mathbb C H^1(X, \mathcal O_X) =0$ or $\pi_1(X) =\{\mathrm{id}\}$, i.e., $X$ is simply-connected. Examples of K3-surfaces arise as Kummer surfaces, quartic surfaces in $\mathbb P_3$ or double coverings of $\mathbb P_2$ branched along  smooth curves of degree six. 

Let $X$ be a K3-surface. Triviality of $\mathcal K _X$ is equivalent to the existence of a nowhere vanishing holomorphic 2-form $\omega$ on $X$. Any 2-form on $X$ can be expressed as a complex multiple of $\omega$. We will therefore mostly refer to $\omega$ (or $\omega _X$) as "the" holomorphic 2-form on $X$ \nomenclature{$\omega_X$}{the holomorphic 2-form on a K3-surface $X$}. We denote by $\mathrm{Aut}_\mathcal O (X)=\mathrm{Aut}(X)$ the group of holomorphic automorphisms of $X$ \nomenclature{$\mathrm{Aut}(X)$}{the group of holomorphic automorphisms of $X$} and consider a (finite) subgroup $G \hookrightarrow \mathrm{Aut}(X)$. If the context is clear, the abstract finite group $G$ is identified with its image in $\mathrm{Aut}(X)$. The group $G$ is referred to as a transformation group, symmetry group or automorphism group of $X$. Note that our considerations are independent of the question whether the group $\mathrm{Aut}(X)$ is finite or not. The order of $G$ is denoted by $|G|$\nomenclature{$|G|$}{the order of the group $G$}. 
\begin{definition}
The action of $G$ on $X$ is called \emph{symplectic} \index{symplectic action} if $\omega$ is $G$-invariant, i.e., $g^*\omega = \omega$ for all $g \in G$. 
\end{definition}
For a finite group $G < \mathrm{Aut}(X)$ we denote by $G_\mathrm{symp}$ the subgroup of symplectic transformations in $G$ \nomenclature{$G_\mathrm{symp}$}{the subgroup of symplectic transformations in $G$}. This group is the kernel of the homomorphism $\chi : G \to  \mathbb{C}^*$ defined by the action of $G$ on the space of holomorphic 2-forms $\Omega^2(X) \cong \mathbb C \omega$. It follows that $G$ fits into the short exact sequence
\begin{equation*}
\{\mathrm{id}\} \to G_\mathrm{symp} \to G \to C_n \to \{\mathrm{id}\}
\end{equation*}
for some cyclic group $C_n$\nomenclature{$C_n$}{the cyclic group of order $n$}. 
If both $G_\text{symp}$ and $C_n \cong G/G_\text{symp}$ are nontrivial, then $G$ is called a symmetry group of \emph{mixed type} \index{mixed type}. 
\section{Quotients of K3-surfaces}
Let $X$ be a surface and let $G < \mathrm{Aut}(X)$ be a finite subgroup of the group of holomorphic automorphisms of $X$. The orbit space $X/G$ carries the structure of a reduced, irreducible, normal complex space of dimension 2 where the sheaf of holomorphic functions is given by the sheaf $G$-invariant functions on $X$. In many cases, the quotient is a singular space. The map $X \to X/G$ is referred to as a quotient map or a covering (map). 

For reduced, irreducible complex spaces $X,Y$ of dimension 2 a proper holomorphic map $f: X \to Y$ is called \emph{bimeromorphic} \index{bimeromorphic map} if there exist proper analytic subsets $A \subset X$ and $B \subset Y$ such that $f: X \backslash A \to Y \backslash B$ is biholomorphic. A holomorphic, bimeromorphic map $f: X \to Y$ with $X$ smooth is a \emph{resolution of singularities of $Y$} \index{resolution of singularities}. 
\begin{definition}
A resolution of singularities $f: X \to Y$ is called \emph{minimal}\index{resolution of singularities! minimal resolution of singularities} if it does not contract any (-1)-curves, i.e., there is no curve $E \subset X$ with $E \cong \mathbb{P}_1$ and $E^2=-1$ such that $f(E) = \{\text{point}\}$. 
\end{definition}
Every normal surface $Y$ admits a minimal resolution of singularities $f: X \to Y$  which is uniquely determined by $Y$. In particular, this resolution is equivariant.
\subsection{Quotients by finite groups of symplectic transformations}

In the study and classification of finite groups of symplectic transformations on K3-surfaces, the following well-known result has proved to be very useful (see e.g. \cite{NikulinFinite})
\begin{theorem}\label{K3quotsymp}
Let $X$ be a K3-surface, $G$ be a finite group of automorphisms of $X$ and $f: Y \to X/G$ be the minimal resolution of $X/G$. Then $Y$ is a K3-surface if and only if $G$ acts by symplectic transformations. 
\end{theorem}
For the reader's convenience we give a detailed proof of this theorem.
We begin with the following lemma.
\begin{lemma}\label{betti}
Let $X$ be a simply-connected surface, $G$ be a finite group of
automorphisms and $f: Y \to X/G$ be an arbitrary
resolution of singularities of $X/G$. Then
$b_1(Y) =0$.
\end{lemma}
\begin{proof}
We denote by $\pi_1(Y)$ the fundamental group of $Y$ \nomenclature{$\pi_1(X)$}{the fundamental group of $X$} and by $[\gamma] \in \pi_1(Y)$ the homotopy equi\-va\-lence class of a closed continuous path $\gamma$.
The first Betti number is the rank
of the free part of 
$$
H_1(Y) = \pi_1(Y) / [\pi_1(Y), \pi_1(Y)].
$$
We show that
for each 
$ [\gamma] \in \pi_1(Y)$ there exists $N \in \mathbb N$ such that $[\gamma]^N =0$, i.e., $\gamma ^N$ is
homotopic to zero for some $N \in \mathbb N$ . It then
follows that $H_1(Y)$ is a torsion group and  $b_1(Y) =0$.

Let $C \subset X/G$ be the union of branch curves of the
covering $q: X \to X/G$, let $P \subset X/G$ be the set of isolated singularities of $X/G$, and $E \subset Y$ be the exceptional locus of $f$. Let $\gamma: [0,1] \to Y$ be a closed path in $Y$. By choosing a path homotopic to
$\gamma$ which does not intersect $E \cup f^{-1}(C)$ we may assume without loss of generality that $\gamma \cap (E \cup f^{-1}(C)) =
\emptyset$.

The path $\gamma$ is mapped to a
closed path in $(X/G)\backslash (C \cup P)$ which we denote also by $\gamma$. The quotient $q:X \to X/G$ is unbranched
outside $C \cup P$ and we can lift $\gamma$ to a path $\widetilde{\gamma}$ in $X$.
Let $\widetilde{\gamma}(0) = x \in X$, then $\widetilde{\gamma}(1) = g.x$ for some $g
\in G$. Since $G$ is a finite group, it follows that $\widetilde{\gamma^N}$ is
closed for some $N \in \mathbb N$.

As $X$ is simply-connected, we know that also $X \backslash q^{-1}(P)$ is
simply-connected. So $\widetilde{\gamma^N }$ is homotopic to zero in $X\backslash q^{-1}(P)$. We can map the corresponding homotopy to $(X/G)\backslash P$ and conclude that
$\gamma^N$ is homotopic to zero in $(X/G)\backslash P$. It follows that $\gamma
^N$ is homotopic to zero in $Y\backslash E$ and therefore in $Y$. 
\end{proof}
\begin{proof}[Proof of Theorem \ref{K3quotsymp}]
We let $E \subset Y$ denote the exceptional locus of the map $f:Y \to X/G$.
If $Y$ is a K3-surface, let $\omega_Y$ denote the nowhere vanishing holomorphic 2-form on $Y$. Let $(X/G)_\text{reg}$ denote the regular part of $X/G$.  Since $f|_{Y\backslash E}: Y\backslash E \to (X/G)_\text{reg}$ is biholomorphic, this defines a holomorphic 2-form $\omega_{(X/G)_\text{reg}}$ on  $(X/G)_\text{reg}$. Pulling this form back to $X$, we obtain a $G$-invariant holomorphic 2-form on $\pi^{-1}((X/G)_\text{reg}) = X \backslash \{p_1,\dots p_k\}$. This extends to a nonzero, i.e., not identically zero, $G$-invariant holomorphic 2-form on $X$. In particular, any holomorphic 2-form on $X$ is $G$-invariant and the action of $G$ is by symplectic transformations.

Conversely, if $G$ acts by symplectic transformations on $X$, then $\omega_X$ defines a nowhere vanishing holomorphic 2-form on $(X/G)_\text{reg}$ and on $ Y \backslash E$ . Our aim is to show that it extends to a nowhere vanishing holomorphic 2-form on $Y$. In combination with Lemma \ref{betti} this yields that $Y$ is a K3-surface.

 Locally at $p \in X$ the action of $G_p$ can be linearized. I.e., there exist a neighbourhood of $p$ in $X$ which is $G_p$-equivariantly isomorphic to a neighbourhood of $0 \in \mathbb C^2$ with a linear action of $G_p$. A neighbourhood of $\pi(p) \in X/G$ is isomorphic to a neighbourhood of the origin in $\mathbb C^2 / \Gamma$ for some finite subgroup $\Gamma < \mathrm{SL}(2,\mathbb C)$. In particular, the points with nontrivial isotropy are isolated. The singularities of $X/G$ are called simple singularities, Kleinian singularities, Du Val singularities or rational double points. Following   \cite{shafarevic} IV.4.3, we sketch an argument which yields the desired extension result.

Let $X \times_{(X/G)} Y = \{ (x,y) \in X\times Y \, | \, \pi(x) = f(y) \}$ and let $ N$ be its normalization. Consider the diagram 
$$
\begin{xymatrix}
{
X\ar[d]_{\pi}  & \ar[l]_{p_X}\ar[d]^{p_Y}N \\
X/G & \ar[l]^f Y.
}
\end{xymatrix}
$$

We let $\omega_X$ denote the nowhere vanishing holomorphic 2-form on $X$. Its pullback $p_X^*\omega_X$ defines a nowhere vanishing holomorphic 2-form on $N_\text{reg}$. Simultaneously, we consider the meromorphic 2-form $\omega _Y$ on $Y$ obtained by pulling back the 2-form on $X/G$ induced by the $G$-invariant 2-form $\omega_X$. By contruction, the pullback $p_Y^*\omega_Y $ coincides with the pullback $p_X^*\omega_X$ on $N_\text{reg}$. 

Consider the finite holomorphic map $p_Y|_{N_\text{reg}}: N_\text{reg} \to p_Y(N_\text{reg}) \subset Y$. Since $p_Y^*\omega_Y $ is holomorphic on $N_\text{reg}$, one checks (by a calculation in local coordinates) that $\omega_Y$ is holomorphic on $p_Y(N_\text{reg}) = Y \backslash \{y_1, \dots y_k\}$ and consequently extends to a holomorphic 2-form on $Y$. Since $p_X^*\omega_X =p_Y^*\omega_Y $ is nowhere vanishing on $N_\text{reg}$, it follows that $\omega_Y$ defines a global, nowhere vanishing holomorphic 2-form on $Y$. 
\end{proof}
\begin{remark}
Let $g$ be a symplectic automorphism of finite order on a K3-surface $X$.
Since K3-surfaces are simply-connected, the covering $X \to X/\langle g \rangle$ can never be unbranched. It follows that $g$ must have fixed points.
\end{remark}

Using Theorem \ref{K3quotsymp} we give an outline of Nikulin's classification of finite Abelian groups of symplectic transformations on a K3-surface \cite{NikulinFinite}. Let $C_p$ be a cyclic group of prime order acting on a K3-surface $X$ by symplectic transformations and $Y$ be the minimal desingularization of the quotient $X/C_p$. 

Notice that by adjunction the self-intersection number of a curve $D$ of genus $g(D)$ on a K3-surface is given by $D^2 = 2g(D)-2$. In particular, if $D$ is smooth, then $D^2 = -e(D)$.

 The exceptional locus of the map $Y \to X/G$ is a union of (-2)-curves and one can calculate their contribution to the topological Euler characteristic $e(Y)$ \nomenclature{$e(X)$}{the topological Euler characteristic of $X$} in relation to $e(X/C_p)$. Let $n_p$ denote the number of fixed point of $C_p$ on $X$. Then 
\begin{align*}
24 &= e(X) = p\cdot e(X/G) - n_p \\
24 &= e(Y) = e(X/G) + n_p \cdot p.
\end{align*}
Combining these formulas gives $n_p = 24/(p+1)$. For a general finite Abelian group $G$ acting symplectically on a K3-surface $X$, one needs to consider all possible isotropy groups $G_x$ for $x \in X$. By linearization, $G_x < \mathrm{SL}_2(\mathbb C)$. Since $G$ is Abelian, it follows that $G_x$ is cyclic and an analoguous formula relating the Euler characteristic of $X$, $X/G$, and $Y$ can be derived. A case by case study then yields Nikulin's classification. In particular, we emphasize the following remark.
\begin{remark}\label{order of symp aut}
 If $g \in \mathrm{Aut}(X)$ is a symplectic automorphism of finite order $n(g)$ on a K3-surface $X$, then  $n(g)$ is bounded by eight and the number of fixed points of $g$ is given by the following table:
\begin{table}[h]
\centering
\begin{tabular}{c|c|c|c|c|c|c|c}
$n(g)$ & 2 & 3 & 4 & 5 & 6 & 7 & 8 \\ \hline
$|\mathrm{Fix}_X(g)|$ & 8 & 6 & 4 & 4& 2& 3 & 2
\end{tabular}
\caption{Fixed points of symplectic automorphisms on K3-surfaces}\label{fix points symplectic}
\end{table}
\end{remark}
\subsection{Quotients  by finite groups of nonsymplectic transformations}
In this subsection we consider the quotient of a K3-surface $X$ by a finite group $G$ such that $G \neq G_\text{symp}$, i.e., there exists $g  \in G$ such that $g^* \omega \neq \omega$.  We prove
\begin{samepage}
\begin{theorem}\label{K3quotnonsymp}
Let $X$ be a K3-surface and let $G < \mathrm{Aut}(X)$ be a finite group such that $g^* \omega \ne \omega$ for some $g \in G $. Then either 
\begin{itemize}
\item[-]
$X/G$ is rational, i.e., bimeromorphically equivalent to $\mathbb P_2$, or
\item[-]
the minimal desingularisation of $X/G$ is a minimal Enriques surface and $$G/G_\text{symp}  \cong  C_2.$$ In this case,  $\pi: X \to X/G$ is unbranched if and only if $G_\text{symp} = \{ \mathrm{id}\}$.
\end{itemize}
\end{theorem}
\end{samepage}
Before giving the proof, we establish the necessary notation and state two useful lemmata.
We denote by $\pi: X \to X/G$ the quotient map. This map can be ramified at isolated points
and along curves. Let $P=
\{p_1, \dots, p_n\}$ denote the set of singularities of $X/G$. For simplicity, the denote the correspondig subset $\pi^{-1}(P)$ of $X$ also by $P$. Outside $P$, the map $\pi$ is ramified along curves $C_i$ of ramification order $c_i+1$. We write $C = \sum c_i C_i$.

Let $r: Y \to X/G$ denote a minimal resolution of singularities of
$X/G$. The exceptional locus of $r$ in $Y$ is denoted by $D$.
As $Y$ is not necessarily a minimal surface, we denote by $p: Y \to Y_ \text{min}$
the sucessive blow-down of (-1)-curves. The union of exceptional curves of $p$
is denoted by $E$. 
\[
\begin{xymatrix}{
C \subset \mathbf{X}\supset P \ar[d]^\pi  \\ 
\pi(C) \subset \mathbf{X/G} \supset P & D \subset \mathbf{Y} \supset E \ar[l]_>>>>>>{r} \ar[d]^p\\
 &  \mathbf{Y_\text{min}}
}
\end{xymatrix}
\]
The following two lemmata (cf. e.g. \cite{BPV} I.16 and Thm. I.9.1) will be useful in order to relate the canonical bundles of the spaces $X$, $(X/G)_\text{reg}$, $Y$ and $Y_ \text{min}$. For a divisor $D$ on a manifold $X$ we denote by $\mathcal O_X(D)$ the line bundle associated to $D$ \nomenclature{$\mathcal O_X(D)$}{the line bundle associated to the divisor $D$}.
\begin{lemma}\label{adj1}
Let $X,Y$ be surfaces and let $\varphi:X \to Y$ be a surjective
finite proper holomorphic map ramified along a curve $C$ in $X$ of ramification order
$k$. Then
\[
\mathcal{K}_X = \varphi^*( \mathcal K_Y ) \otimes \mathcal O_X (C)^{\otimes(k-1)}.
\]
More generally, if $\pi$ is ramified along a ramification divisor $R = \sum_i r_i R_i$, where $R_i$ is an irreducible curve and $r_i + 1$ is the ramification order of $\pi$ along $R_i$, then
\[
\mathcal{K}_X = \pi^*( \mathcal K_Y ) \otimes \mathcal O_X (R).
\]
\end{lemma}
\begin{lemma}\label{adj2}
Let $X$ be a surface and let $b: X
\to Y$ be the blow-down of a (-1)-curve $E \subset X$. Then 
\[
\mathcal{K}_X = b^* ( \mathcal K_Y )  \otimes \mathcal O_X (E).
\]
\end{lemma}
We present a proof of Theorem \ref{K3quotnonsymp} using the Enriques Kodaira classification of surfaces.
\begin{proof}[Proof of Theorem \ref{K3quotnonsymp}]
The Kodaira dimension of the K3-surface $X$ is $\mathrm{kod}(X)=0$. 
The Kodaira dimension of $X/G$, which is defined as the Kodaira dimension of some resolution of $X/G$, is less than or equal to the Kodaira dimension of $X$. (c.f. Theorem 6.10 in \cite{ueno}),
\[
0=\mathrm{kod}(X) \geq \mathrm{kod}(X/G) = \mathrm{kod}(Y) =
\mathrm{kod}(Y_\text{min}) \in \{0, -\infty\}.
\]
By Lemma \ref{betti}, the first Betti number of $Y$ and $Y_\text{min}$ is zero. 
If $\mathrm{kod}(Y) = - \infty$, then $Y$ is a smooth rational surface. If $\mathrm{kod}(Y) =\mathrm{kod}(Y_\text{min})= 0$, then, since $Y$ is not a K3-surface by Theorem \ref{K3quotsymp}, it follows that $Y_\text{min}$ is an Enriques surface. 

If $Y_\text{min}$ is an Enriques surface, then $\mathcal K_{Y_\text{min}}^{\otimes 2}$ is trivial. Let $ s \in
  \Gamma(Y_\text{min}, \mathcal{K}_{Y_\text{min}}^{\otimes 2})$ be a nowhere vanishing section.
Consecutive application of Lemma \ref{adj2} yields the following formula
\[
\mathcal{K}_Y^{\otimes 2} = (p^* \mathcal{K}_{Y_\mathrm{min}})^{\otimes 2} \otimes \mathcal O_Y(E)^{\otimes 2}= p^* (\mathcal{K}_{Y_\mathrm{min}}^{\otimes 2}) \otimes \mathcal O_Y(E)^{\otimes 2}.
\]
Let $e \in \Gamma(Y, \mathcal O_Y(E)^{\otimes 2})$ and write
$\tilde{s}= p^*(s) \cdot e$. This global section of $\mathcal{K}_Y^{\otimes 2}$ vanishes along $E$ and is
nowhere vanishing outside $E$.
By restricting $\tilde{s}$ to $Y\backslash D$ we obtain a section of
$\mathcal{K}_{Y \backslash D}^{\otimes 2}$. Since $\pi$ is biholomorphic
outside $D$, we can map the restricted section to $(X/G)\backslash P =
(X/G)_\text{reg}$ and obtain a section $\hat{s}$ of
$\mathcal{K}_{(X/G)_\text{reg}}^{\otimes 2}$. Note that $\hat{s}$ is not the
zero-section. If $E \neq \emptyset$, i.e., $Y$ is not minimal, let $E_1 \subset E$ be a (-1)-curve. The minimality of the resolution $r: Y \to X/G$ implies $E_1 \nsubseteq D$. It follows that $\hat{s}$ vanishes along the image of $E_1$ in $(X/G)_\text{reg}$

We may now apply Lemma \ref{adj1} to the map $\pi|_{X\backslash P}$ to see
\begin{align*}
\mathcal{K}_{X \backslash P}^{\otimes 2} &= (\pi^*
\mathcal{K}_{(X/G)_\text{reg}})^{\otimes 2} \otimes \mathcal O_{X\backslash P}(C)^{\otimes 2}\\
 &=\pi^* (\mathcal{K}_{(X/G)_\text{reg}}^{\otimes 2}) \otimes \mathcal O_{X\backslash P}(C)^{\otimes 2}.
\end{align*}
Let $c \in \Gamma(X\backslash P,
\mathcal O_{X\backslash P}(C))^{\otimes 2}$. Then $t := \pi^*
\hat{s} \cdot c \in \Gamma(X\backslash P, \mathcal{K}_{X\backslash P}^{\otimes2})$ 
is not the zero-section but
vanishes along $C$ and along the preimage of the zeroes of $\hat{s}$.

Now $t$ extends to a holomorphic section $\tilde{t} \in \Gamma(X,
\mathcal{K}_X^{\otimes 2})$. Since $X$ is K3, it follows that both $\mathcal{K}_X$ and  $\mathcal{K}_X^{\otimes 2}$ are trivial and $\tilde t$ must be nowhere vanishing. Consequently,
both $E$ and $C$ must be empty.
It follows that the map $\pi$ is at worst branched at points $P$ (not along curves) and the
minimal resolution  $Y$ of $X/G$ is a minimal surface. 
\[
\begin{xymatrix}{
P \subset \mathbf{X} \ar[d]^\pi \\
P  \subset \mathbf{X/G} &  \mathbf{Y} \supset D \ar[l]_>>>>>>{r}  
}
\end{xymatrix}
\]
The section $\tilde{t}$ on $X$ is $G$-invariant by construction. Let $\omega$
be a nonzero section of the trivial bundle $\mathcal{K}_X$ such that $\tilde t = \omega ^2$. The action of $G$ on $X$ is
nonsymplectic, therefore $\omega$ is not invariant but $\tilde{t}$ is. Hence
$G$ acts on $\omega$ by multiplication with $\{1,-1\}$ and $G/G_\text{symp} \cong C_2$.

If $\pi: X \to X/G$ is unbranched, it follows that $\mathrm{Fix}_X(g) = \emptyset$ for all $g \in G \backslash \{\mathrm{id}\}$. Since symplectic automorphisms of finite order necessarily have fixed points, this implies $G_\text{symp} = \{\mathrm{id}\}$.

Conversely,
if $G$ is isomorphic to $C_2$,
it remains to show that the set $P=\{p_1,\dots,p_n\}$ is empty. 
Our argument uses the Euler characteristic $e$ of $X$, $X/G$, and $Y$. By chosing a triangulation of $X/G$
such that all points $p_i$ lie on vertices we calculate $24= e(X) = 2e(X/G)-n$.
Blowing up the $C_2$-quotient singularities of $X/G$ we obtain $12=e(Y) = e(X/G) +n$.
This implies $e(X/G) =12$ and $n=0$ and completes the proof of the theorem.
\end{proof}
\section{Antisymplectic involutions on K3-surfaces}
As a special case of the theorem above we consider the quotient of a K3-surface $X$ by an involution $\sigma \in \mathrm{Aut}(X)$ which acts on the 2-form $\omega$ by multiplication by $-1$ and is therefore called \emph{antisymplectic involution}\index{antisymplectic involution}.
\begin{proposition}\label{K3quotnonsympinvo}
Let $\pi:X \to X/\sigma$ be the quotient of a K3-surface by an antisymplectic involution $\sigma$. If $\mathrm{Fix}_X(\sigma) \neq  \emptyset$, then $\mathrm{Fix}_X(\sigma)$ is a disjoint union of smooth curves and $X/\sigma$ is a smooth rational surface. Furthermore, $\mathrm{Fix}_X(\sigma) = \emptyset$ if and only if $X/\sigma$ is an Enriques surfaces.
\end{proposition}
\begin {proof}
If $\mathrm{Fix}_X(\sigma) \neq  \emptyset$, then Theorem \ref{K3quotnonsymp} and linearization of the $\sigma$-action at its fixed points yields the proposition.
If $\mathrm{Fix}_X(\sigma) = \emptyset$, then $X \to X/\sigma$ is unbranched and $\mathrm{kod}(X) = \mathrm{kod}(X/G)$. It follows that $X/G$ is an Enriques surface.
\end {proof}
In order to sketch Nikulin's description of the fixed point set of an anti\-symplec\-tic involution we summarize some information about the Picard lattice of a K3-surface.
\subsection{Picard lattices of K3-surfaces}
Let $X$ be a complex manifold. The \emph{Picard group of $X$} \index{Picard group} \nomenclature{$\mathrm{Pic}(X)$}{the Picard group of $X$} is the group of isomorphism classes of line bundles on $X$ and denoted by $\mathrm{Pic}(X)$. It is isomorphic to $H^1(X, \mathcal O_X^*)$. Let $\mathbb Z_X$ denote the constant sheaf on $X$ corresponding to $\mathbb Z$, then the exponential sequence $0 \to \mathbb Z_X \to \mathcal O_X \to \mathcal O_X^* \to 0$ induces a map
\[
\delta: H^1(X, \mathcal O_X^*) \to H^2(X, \mathbb Z).
\]
Its kernel is the identity component $\mathrm{Pic}^0(X)$ of the Picard group. The quotient $\mathrm{Pic}(X) / \mathrm{Pic}^0(X)$ is isomorphic to a subgroup of $ H^2(X, \mathbb Z)$ and referred to as the \emph{N\'eron-Severi group $NS(X)$ of $X$} \index{N\'eron-Severi group}\nomenclature{$NS(X)$}{the N\'eron-Severi group of $X$}. On the space $ H^2(X, \mathbb Z)$ there is the natural intersection or cupproduct pairing. The rank of the N\'eron-Severi group of $X$ is denoted by $\rho(X)$ and referred to as the \emph{Picard number of $X$}\nomenclature{$\rho(X)$}{the Picard number of $X$}

If $X$ is a K3-surface, then $H^1(X, \mathcal O_X)=\{0\}$ and $\mathrm{Pic}(X)$ is isomorphic to $NS(X)$. In particular, the Picard group carries the structure of a lattice, the \emph{Picard lattice} \index{Picard lattice} of $X$, and is regarded as a sublattice of $ H^2(X, \mathbb Z)$, which is known to have signature $(3,19)$ (cf. VIII.3 in \cite{BPV}). 

If $X$ is an algebraic K3-surface, i.e., the transcendence degree of the field of meromorphic functions on $X$ equals 2, then $\mathrm{Pic}(X)$ is a nondegenerate lattice of signature $(1, \rho -1)$ (cf. \S3.2 in \cite{NikulinFinite}).
\subsection{The fixed point set of an antisymplectic involution}
We can now present Nikulin's classification of the fixed point set of an antisymplectic involution on a K3-surface \cite{NikulinFix}. 
\begin{theorem}\label{FixSigma}
The fixed point set of an antisymplectic involution $\sigma$ on a K3-surface $X$ is one of the following three types:
\[
\text{1.)\ }\ \mathrm{Fix}(\sigma) = D_g \cup \bigcup_{i=1}^n R_i,
\quad \quad
\text{2.)\ }\ \mathrm{Fix}(\sigma) = D_1 \cup D'_1,
\quad ¸\quad
\text{3.)\ }\ \mathrm{Fix}(\sigma)= \emptyset,
\]
 where $D_g$ denotes a smooth curve of genus $g \geq 0$ and $\bigcup_{i=1}^n R_i$
 is a possibly empty union of smooth disjoint rational curves.  In case 2.), $D_1$ and $D'_1$ denote disjoint elliptic curves.
\end{theorem}
\begin{proof}
Assume there exists a curve $D_g$ of genus $g \geq 2$ in  $\mathrm{Fix}(\sigma)$. By adjunction, this curve has positive self-intersection. We claim that each curve $D$ in $\mathrm{Fix}(\sigma)$ disjoint from $D_g$ is rational.

First note that the existence of an antisymplectic automorphism on $X$ implies that $X$ is algebraic (cf. Thm. 3.1 in \cite{NikulinFinite}) and therefore $\mathrm{Pic}(X)$ is a nondegenerate lattice of signature $(1, \rho -1)$.

If $D$ is elliptic, then $D^2 =0$, $D_g^2>0$ and $D \cdot D_g=0$ is contrary to the fact that $\mathrm{Pic}(X)$ has signature $(1,\rho-1)$. If $D$ is of genus $\geq 2$, then $D^2 >0$  and we obtain the same contradiction.

Now assume that there exists an elliptic curve $D_1$ in $\mathrm{Fix}(\sigma)$. By the considerations above, there may not be curves of genus $\geq 2$ in $\mathrm{Fix}(\sigma)$. If there are no further elliptic curves in $\mathrm{Fix}(\sigma)$, we are in case 1) of the classification. If there is another elliptic curve $D'_1$ in $\mathrm{Fix}(\sigma)$, this has to be linearly equivalent to $D_1$, as otherwise the intersection form of $\mathrm{Pic}(X)$ would degenerate on the span of $D_1$ and $D'_1$.
The linear system of $D_1$ defines an elliptic fibration $X \to \mathbb{P}_1$. The induced action of $\sigma$ on the base may not be trivial since this would force $\sigma$ to act trivially in a neighbourhood of $D_1$ in $X$. It follows that the induced action of $\sigma$ on $\mathbb P_1$ has precisely two fixed points and that $\mathrm{Fix}(\sigma)$ contains no other curves than $D_1$ and $D'_1$.
This completes the proof of the theorem.
\end{proof}
\section{Finite groups of symplectic automorphisms}
In preparation for stating Mukai's classification of finite groups of symplectic automorphisms on K3-surfaces we present his list \cite{mukai} of symplectic actions of finite groups $G$ on K3-surfaces $X$. It is an important source of examples, many of these will occur in our later discussion.
For the sake of brevity, at this point we do not introduce the notation of groups used in this table. 
\renewcommand{\baselinestretch}{1.5}
\begin{table}[h]
\centering
\begin{tabular}{l|l|l|l}
  & $G$ & $|G|$ & \textbf{K3-surface} $X$ \\ \hline 
1 & $L_2(7)$ & 168 & $\{x_1^3x_2+x_2^3x_3+x_3^3x_1+x_4^4 =0\} \subset \mathbb P_3$\\ \hline
2 & $A_6$ & 360 & $\{\sum_{i=1}^6 x_i = \sum_{i=1}^6 x_1^2 = \sum_{i=1}^6 x_i^3=0\} \subset \mathbb P_5$\\  \hline
3 & $S_5$ & 120 & $\{\sum_{i=1}^5 x_i = \sum_{i=1}^6 x_1^2 = \sum_{i=1}^5 x_i^3=0\} \subset \mathbb P_5$\\ \hline
4 & $M_{20}$ & 960 & $\{ x_1^4+ x_2 ^4 + x_3^4 +x_4^4 +12 x_1x_2x_3x_4 = 0\} \subset \mathbb P_3$\\ \hline
5 & $F_{384}$ & 384 & $\{ x_1^4+ x_2 ^4 + x_3^4 +x_4^4 = 0\} \subset \mathbb P_3$\\ \hline
6 & $A_{4,4}$ & 288 & $\{x_1^2+x_2^2 +x_3^2 = \sqrt{3}x_4^2\} \cap $\\
& & & 			$ \{x_1^2+ \omega x_2^2 +\omega^2 x_3^2 = \sqrt{3}x_5^2\} \cap $\\
& & & 			$ \{x_1^2+\omega^2 x_2^2 +\omega x_3^2 = \sqrt{3}x_6^2\}  \subset \mathbb P_5$ \\ \hline
7 & $T_{192}$ & 192 & $\{ x_1^4+ x_2 ^4 + x_3^4 +x_4^4 - 2 \sqrt{-3}(x_1^2x_2^2 + x_3^2x_4^2 = 0\} \subset \mathbb P_3$\\\hline
8 & $H_{192}$ & 192 & $ \{x_1^2+x_3^2+x_5^2 = x_2^2 + x_4^2 + x_6^2\} \cap $ \\
& & & 			$ \{x_1^2+x_4^2 = x_2^2+x_5^2=x_3^2+x_6^2\} \subset \mathbb P_5$\\ \hline
9 & $N_{72}$ & 72 & $\{ x_1^3+ x_2 ^3 + x_3^3 +x_4^3= x_1x_2 + x_3x_4+ x_5^2 = 0 \} \subset \mathbb P_4$ \\ \hline
10 & $M_9$ & 72 & Double cover of $\mathbb P_2$ branched along \\
& & & 		$\{x_1^6+y_2^6 +x_3^6 -10(x_1^3x_2^3 + x_2^3x_3^3 +x_3^3x_1^3) =0\}$\\ \hline
11 & $ T_{48}$ & 48 & Double cover of $\mathbb P_2$ branched along \\
& & & 		$\{x_1x_2(x_1^4-x_2^4)+ x_3^6 =0\}$
\end{tabular}
\caption{Finite groups of symplectic automorphisms on K3-surfaces}\label{TableMukai}
\end{table}
\renewcommand{\baselinestretch}{1.1}

The following theorem (Theorem 0.6 in \cite{mukai}) characterizes finite groups of symplectic automorphisms on K3-surfaces.
\begin{theorem}\label{mukaithm}
A finite group $G$ has an effective sympletic actions on a K3-surface if and only if it is isomorphic to a subgroup of one of the eleven groups in Table \ref{TableMukai}. 
\end{theorem}
The "only if"-implication of this theorem follows from the list of eleven examples summarized in Table \ref{TableMukai}. This list of examples is, however, far from being exhaustive. It is therefore desirable to find further examples of K3-surfaces where the groups from this list occur and describe or classify these surfaces with maximal symplectic symmetry.. 
\begin{definition}\label{maxsymp}
By Proposition 8.8 in \cite{mukai} there are no subgroup relations among the eleven groups in Mukai's list. Therefore, the groups are \emph{maximal finite groups of symplectic transformations}\index{maximal group of symplectic trans\-for\-ma\-tions}.  We refer to the groups in this list also as \emph{Mukai groups}. \index{Mukai group}
\end{definition}
\subsection{Examples of K3-surfaces with symplectic symmetry}
We conclude this chapter by presenting two typical examples of K3-surface with symplectic symmetry.
\begin{example}\label{L2(7)example}
The group $L_2(7) = \mathrm{PSL}(2, \mathbb F_7)=\mathrm{GL}_3(\mathbb F_2)$ is a simple group of order 168. It is generated by the three projective transformations $\alpha, \beta, \gamma$ of $\mathbb P_1( \mathbb F_7)$ given by
\[
\alpha(x) = x+1; \quad \beta(x) =2x; \quad \gamma(x) = -x^{-1}.
\]
In terms of these generators, we define a three-dimensional representation of $L_2(7)$ by
\[
\alpha \mapsto
\begin{pmatrix}
\xi & 0 & 0\\
0 & \xi^2 & 0\\
0 & 0 & \xi^4
\end{pmatrix};\quad
\beta \mapsto 
\begin{pmatrix}
0 & 0 & 1\\
1 & 0 & 0\\
0 & 1 & 0
\end{pmatrix};\,
\gamma \mapsto \frac{-1}{\sqrt{-7}}
\begin{pmatrix}
a&b&c\\
b&c&a\\
c&a&b
\end{pmatrix}
\]
where 
$
\xi=e^{\frac{2\pi i }{7}},\,
a=\xi^2-\xi^5,\,
b=\xi-\xi^6,\,
c=\xi^4-\xi^3,
$
and $\sqrt{-7}= \xi+\xi^2+\xi^4-\xi^3-\xi^5-\xi^6$. Klein's quartic curve
\[
C_\text{Klein}= \{x_1x_2^3 + x_2x_3^3 + x_3x_1^3=0\} \subset \mathbb P_2
\]
is invariant with respect to induced action of $L_2(7)$ on $\mathbb P_2$. Mukai's example of a K3-surface with symplectic $L_2(7)$-symmetry is the smooth quartic hypersurface in $\mathbb P_3$ defined by  
\[
X_\text{KM} = \{x_1x_2^3 + x_2x_3^3 + x_3x_1^3 + x_4^4=0\} \subset \mathbb P_3,
\]
where the action of $L_2(7)$ is defined to be trivial on the coordinate $x_4$ and defined as above on $x_1,x_2,x_3$. Since $L_2(7)$ is a simple group, it follows 
that the action is effective and symplectic. The surface $ X_\text{KM}$ is called the \emph{Klein-Mukai surface}\index{Klein-Mukai surface}. By construction, it is a cyclic degree four cover of $\mathbb P_2$ branched along Klein's quartic curve. In fact, there is an action of the group $L_2(7) \times C_4$ on $X_\text{KM}$, where the action of $C_4$ is by nonsymplectic transformations. The Klein-Mukai surface will play an important role in Sections \ref{KMsurface} and \ref{168}. 
\end{example}
\subsubsection{Cyclic coverings}
Since many examples of K3-surfaces are constructed as double covers we discuss the construction of branched cyclic covers with emphasis on group actions induced on the covering space.

Let $Y$ be a surface such that Picard group of $Y$ has no torsion, i.e., there does not exist a nontrivial line bundle $E$ on $Y$ such that $E^{\otimes n}$ is trivial for some $n \in \mathbb N$.

Let $B$ be an effective and reduced divisor on $Y$ and suppose there exists a line bundle $L$ on $Y$ such that $\mathcal O_Y(B) = L^{\otimes n}$ and a section $s \in \Gamma( Y , L^{\otimes n})$ whose zero-divisor is $B$. Let $p: L \to L^{\otimes n}$ denote the bundle homomorphism mapping each element $(y,z) \in L$ for $y \in Y$ to $(y,z^n ) \in L^{\otimes n}$. The preimage $X = p^{-1}(\mathrm{Im}(s))$ of the image of $s$ is an analytic subspace of $L$. The bundle projection $L \to Y$ restricted to $X$ defines surjective holomorphic map $X \to Y$ of degree $n$. 
\[
\begin{xymatrix}{
 X \subset L \ar[r]^p \ar[d]& L^{\otimes n} \supset \mathrm{Im}(s) \ar[d]\\
Y \ar[r]_{\mathrm{id}}& Y \ar@/_1pc/[u]_s
} 
\end{xymatrix}
\]
Since $\mathrm{Pic}(Y)$ is torsion free, the line bundle $L$ is uniquely determined by $B$. It follows than $X$ is uniquely determined and we refer to $X$ as \emph{the} cyclic degree $n$ covering of $Y$ branched along $B$. We note that $X$ is normal and irreducible. It is smooth if the divisor $B$ is smooth. (cf. I.17 in \cite{BPV})

Let $G$ be a finite group in $\mathrm{Aut}(Y)$ and assume that the divisor $B$ is invariant, i.e., $gB =B$ for all $g \in G$. Then the pull-back bundle $g^* L^{\otimes n}$ is isomorphic to $L^{\otimes n}$. We consider the group $\mathrm{BAut}(L^{\otimes n})$ of bundle maps of $ L^{\otimes n}$ and the homomorphism $\mathrm{BAut}(L^{\otimes n}) \to \mathrm{Aut}(Y)$ mapping each bundle map to the corresponding automorphism of the base. Its kernel is isomorphic to $\mathbb C^*$. The observation $g^* L^{\otimes n} \cong L^{\otimes n}$ implies that the group $G$ is contained in the image of $\mathrm{BAut}(L^{\otimes n})$ in $\mathrm{Aut}(Y)$. 

By assumption, the zero set of the section $s$ is $G$-invariant. 
The bundle map induced by $g^*$ maps the section $s$ to a multiple $\chi(g) s$ of $s$ for some character $\chi: G \to \mathbb C^*$. 
It follows that the bundle map $\tilde g$ induced by $\chi(g)^{-1} g^*$ stabilizes the section.
The group $\tilde G = \{ \tilde g \, | \, g \in G \} \subset  \mathrm{BAut}(L^{\otimes n})$ is isomorphic to $G$ and stabilizes $\mathrm{Im}(s) \subset L^{\otimes n} $. 

In order to define a corresponding action on $X$, first observe that $g^* L \cong L$ for all $g \in G$. This follows from the observation that $ g^* L \otimes L^{-1}$ is a torsion bundle and the assumption that $\mathrm{Pic}(Y)$ has no torsion. As above, we deduce that the group $G$ is contained in the image of $\mathrm{BAut}(L)$ in $\mathrm{Aut}(Y)$. Let $\overline G$ be the preimage of $G$ in $\mathrm{BAut}(L)$. Then $\overline G$ is a central $\mathbb C^*$-extension of $G$,
\[
 \{\mathrm{id}\} \to \mathbb C^* \to \overline G  \to G \to \{\mathrm{id}\}.
\]
The map $p: L \to L^{\otimes n}$ induces a homomorphism $p_*:\mathrm{BAut}(L) \to \mathrm{BAut}(L^{\otimes n})$. Its kernel is isomorphic to $C_n < \mathbb C^*$ and we consider the preimage $H = p_*^{-1}(\tilde G)$ in $\mathrm{BAut}(L)$. The group $H < \overline G$ is a central $C_n$-extension of $\tilde G \cong G$,
\[
 \{\mathrm{id}\} \to C_n \to H \to G \to \{\mathrm{id}\}.
\]
By construction, the subset $X \subset L$ is invariant with respect to $H$. This discussion proves the following proposition.
\begin{proposition}
Let $Y$ by a surface such that $\mathrm{Pic}(Y)$ is torsion free and $G < \mathrm{Aut}(Y)$ be a finite group. If $B \subset Y$ is an effective, reduced, $G$-invariant divisor defined by a section $s \in \Gamma( Y, L^{\otimes n})$ for some line bundle $L$, then the cyclic degree $n$ covering $X$ of $Y$ branched along $B$ carries the induced action of a central $C_n$-extension $H$ of $G$ such that the covering map $\pi: X \to Y$ is equivariant.
\end{proposition}
\begin{example}[Double covers]
For any finite subgroup $G < \mathrm{PSL}(3,\mathbb C)$ and any $G$-invariant smooth curve $C \subset \mathbb P_2$ of degree six, the double cover $X$ of $\mathbb P_2$ branched along $C$ is a K3-surface with an induced action of a degree two central extension of the group $G$. 
Many interesting examples (no. 10 and 11 in Mukai's table) can be contructed this way. For example, the Hessian of Klein's curve  $\mathrm{Hess}(C_\text{Klein})$ is an $L_2(7)$-invariant sextic curve and the double cover of $\mathbb{P}_2$ branched along $\mathrm{Hess}(C_\text{Klein})$ is a K3-surface with a symplectic action of $L_2(7)$ centralized by the antisymplectic covering involution (cf. Section \ref{168}). 
\end{example}
%
%
%
%
%
\chapter{Equivariant Mori reduction}\label{chapter mmp}
This chapter deals with a detailed discussion of Example 2.18 in \cite{kollarmori} (see also Section 2.3 in \cite{Mori}) and introduces a minimal model program for surfaces
respecting finite groups of symmetries. Given a projective algebraic
surface $X$ with $G$-action, in analogy to the usual minimal model program, one obtains from $X$ a $G$-minimal model $X_{G\text{-min}}$ by a finite number of $G$-equivariant blow-downs, each contracting a finite number of disjoint (-1)-curves. The surface $X_{G\text{-min}}$ is either a conic bundle over a smooth curve, a Del Pezzo surface or has nef canonical bundle. The case $G \cong C_2$ is also discussed in \cite{Bayle}, the case $G \cong C_p$ for $p$ prime in \cite{fernex}. As indicated in the introduction, applications can be found throughout the literature.
\section{The cone of curves and the cone theorem}

Throughout this chapter we let $X$ be a smooth projective algebraic surface and let $\mathrm{Pic}(X)$ denote the group of isomorphism classes of line bundles on $X$.  
\begin{definition}
A \emph{divisor}\index{divisor} on $X$ is a formal linear combination of irreducible curves $D = \sum a_i C_i$ with $a_i \in \mathbb Z$.
A \emph{1-cycle} \index{1-cycle} on $X$ is a formal linear combination of irreducible curves $C = \sum b_i C_i$ with $b_i \in \mathbb R$. A 1-cycle is \emph{effective} if $b_i \geq 0$ for all $i$. 

We define a pairing $\mathrm{Pic}(X) \times \{\text{divisors}\} \to \mathbb Z$ by $(L,D) \mapsto L \cdot D = \deg(L|_D)$. Extending by linearity, this defines a pairing $\mathrm{Pic}(X) \times \{\text{1-cycles}\} \to \mathbb R$.
\nomenclature{$L\cdot C$}{the intersection number of a line bundle $L$ and a 1-cycle $C$}\index{intersection number}
We use this notation for the intersection number also for pairs of divisors $C$ and $D$ and write $C\cdot D = \deg(\mathcal O_X(D)|_C)$.
Two 1-cycles $C,C'$ are called \emph{numerically equivalent} \index{numerically equivalent} if $L\cdot C = L \cdot C'$ for all $L \in \mathrm{Pic}(X)$.  We write $C \equiv C'$. The numerical equivalence class of a 1-cycle $C$ is denoted by $[C]$.
The space of all 1-cycles with real coefficients modulo numerical equivalence is a real vector space denoted by $N_1(X)$. 
Note that $N_1(X)$ is finite-dimensional.
\begin{remark}
Let $L$ be a line bundle on $X$ and let $L^{-1}$ denote its dual bundle. Then $L^{-1} \cdot C = -L \cdot C$ for all $[C]\in N_1(X)$. We therefore write $L^{-1} = -L$ in the following.
\end{remark}
\begin{definition}
A line bundle $L$ is called \emph{nef}\index{nef} if $L \cdot C \geq 0$ for all irreducible curves $C$. 
\end{definition}
We set
\[
NE(X) = \{ \sum a_i[C_i] \ | \ C_i \subset X \text{ irreducible curve},\, 0 \leq a_i \in \mathbb R\} \subset N_1(X).
\]
The closure $\overline{NE}(X)$ of $NE(X)$ in $N_1(X)$ is called \emph{Kleiman-Mori cone} or \emph{cone of curves} \index{cone of curves} on $X$. \nomenclature{$\overline{NE}(X)$}{the cone of curves on $X$}
 
For a line bundle $L$, we write $\overline{NE}(X)_{L\geq 0} = \{ [C]\in N_1(X) \ |\ L \cdot C\geq 0 \} \cap \overline{NE}(X)$. Analogously, we define $\overline{NE}(X)_{L\leq 0}$, $\overline{NE}(X)_{L > 0}$, and $\overline{NE}(X)_{L < 0}$.
\end{definition}
Using this notation we phrase Kleiman's ampleness criterion (cf. Theorem 1.18 in \cite{kollarmori})
\begin{theorem}
A line bundle $L$ on $X$ is ample if and only if $\overline{NE}(X)_{L>0} = \overline{NE}(X)\backslash \{0\}$. 
\end{theorem}
\begin{definition}
Let $V$ be a finite-dimensional real vector space . A subset $N \subset V$ is called \emph{cone} if $0\in N$ and $N$ is closed under multiplication by positive real numbers. 
A subcone $M \subset N$ is called \emph{extremal} if $u,v \in N $ satisfy $u,v  \in M$ whenever $u+v \in M$. An extremal subcone is also referred to as an \emph{extremal face}. A 1-dimensional extremal face is called \emph{extremal ray}\index{extremal ray}. 
For subsets $A, B \subset V$ we define $A+B := \{ a+b \,|\, a\in A, b\in B\}$.
\end{definition}
The cone of curves $\overline{NE}(X)$ is a convex cone in $N_1(X)$ and the following cone theorem, which is stated here only for surfaces, describes its geometry (cf. Theorem 1.24 in \cite{kollarmori}).
\begin{theorem}\label{conethm}
Let $X$ be a smooth projective surface and let $\mathcal{K}_X$ denote the canonical line bundle on $X$. 
There are countably many rational curves $C_i \in X$ such that $0 < -\mathcal{K}_X \cdot C_i \leq \mathrm{dim}(X) +1 $ and 
\[
\overline{NE}(X) = \overline{NE}(X) _{\mathcal{K}_X \geq 0} + \sum_i \mathbb R_{\geq 0} [C_i]. 
\]
For any $\varepsilon>0$ and any ample line bundle $L$
\[
\overline{NE}(X) = \overline{NE}(X) _{(\mathcal{K}_X+\varepsilon L) \geq 0} + \sum_{\text{finite}} \mathbb R_{\geq 0} [C_i]. 
\]
\end{theorem}
\section{Surfaces with group action and the cone of invariant curves}
Let $X$ be a smooth projective surface and let $G< \mathrm{Aut}_{\mathcal{O}}(X)$ be a group of holomorphic transformations of $X$. We consider the induced action on the space of 1-cycles on $X$. For $g \in G$ and an irreducible curve $C_i$ we denote by $g C_i$ the image of $C_i$ under $g$. For a 1-cycle $C = \sum a_i C_i$ we define $gC = \sum a_i (g C_i)$. This defines a $G$-action on the space of 1-cycles. 
\begin{lemma}
Let $C_1, C_2$ be 1-cycles and $C_1 \equiv C_2$. Then $gC_1 \equiv gC_2$ for any $g \in G$. 
\end{lemma}
\begin{proof}
The 1-cycle $gC_1$ is numerically equivalent to $g C_2$ if and only if
$L \cdot (gC_1) = L \cdot (gC_2)$ for all $L \in \mathrm{Pic}(X)$. For $g \in G$ and
$L \in \mathrm{Pic}(X)$ let $g^*L$ denote the pullback of $L$ by
$g$. The claim above is equivalent to  $((g^{-1})^*L) \cdot (gC_1) = ((g^{-1})^*L) \cdot (gC_2)$ for all $L \in \mathrm{Pic}(X)$. Now 
$$((g^{-1})^*L) \cdot (gC_1) = \mathrm{deg}((g^{-1})^* L|_{gC_1}) = \mathrm{deg}(L|_{C_1}) = L \cdot C_1 = L \cdot C_2 = (g^{-1})^*L(gC_2)$$ 
for all $L \in\mathrm{Pic}(X)$. 
\end{proof}
This lemma allows us to define a $G$-action on $N_1(X)$ by setting $g[C] := [gC]$ and extending by linearity. We write $N_1(X)^G = \{ [C] \in N_1(X) \ | \ [C]=[gC] \text{ for all } g \in G\}$, the set of invariant 1-cycles modulo numerical equivalence. This space is a linear subspace of $N_1(X)$.

\begin{figure}[H]
\centering
     \subfigure[The cone of curves and its \textcolor{red}{extremal rays}]
     	{\includegraphics[width=0.45\textwidth]{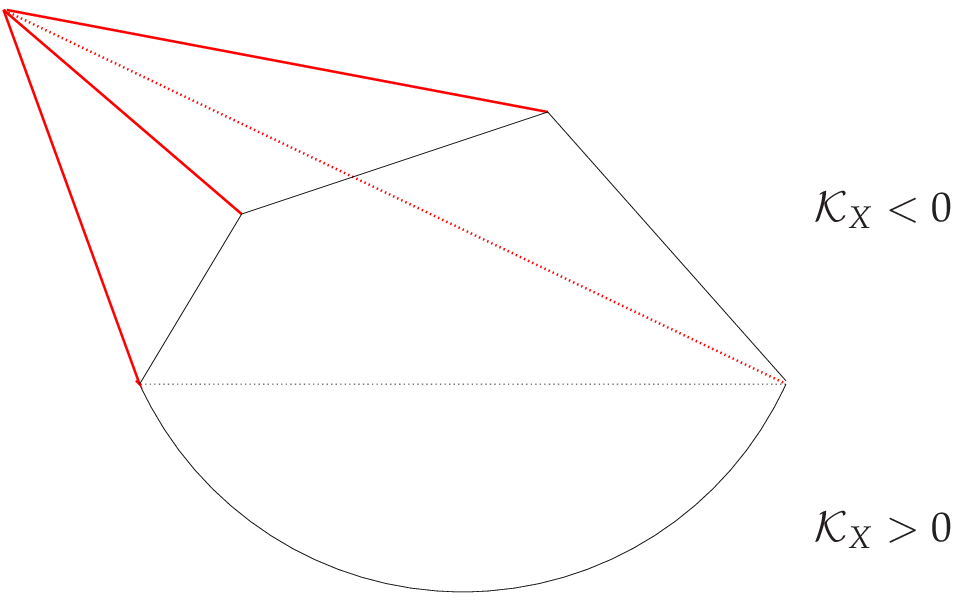}}\hspace{1cm}
     \subfigure[The cone of curves and the \textcolor{green}{invariant subspace $N_1(X)^G$}. 
		Their intersection $\overline{NE}(X)^G$ has a \textcolor{blue}{new extremal ray}.]
     	{\includegraphics[width=0.45\textwidth]{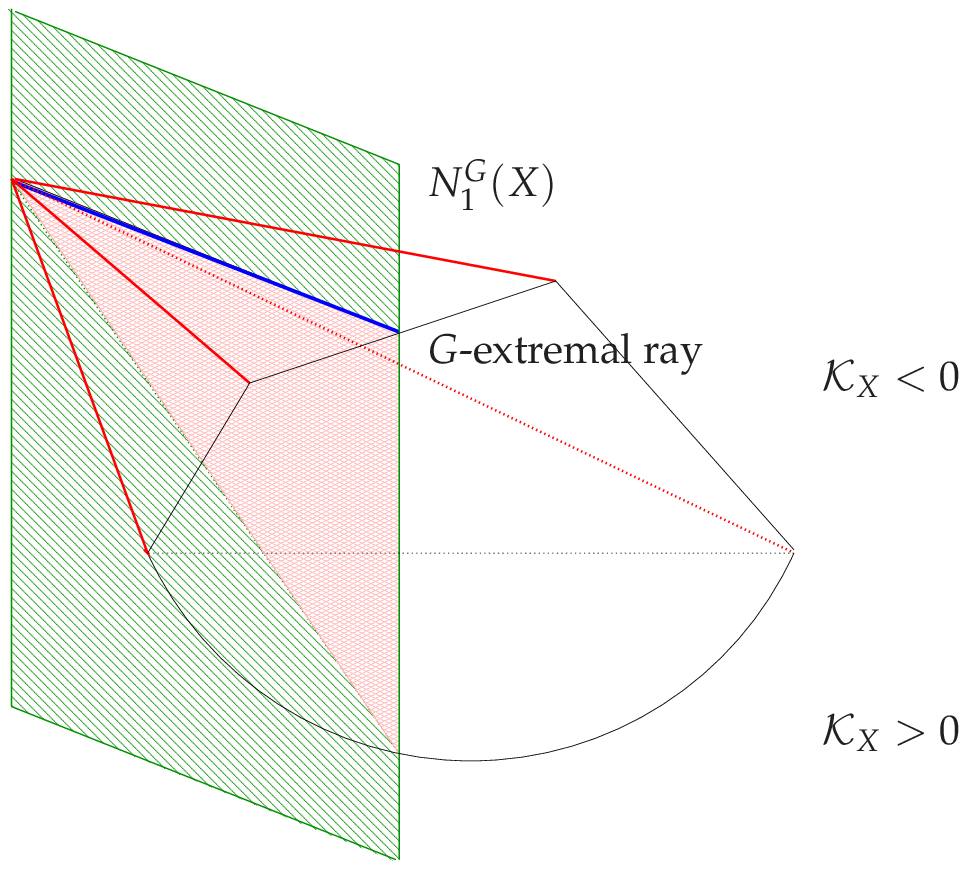}}
      \caption{The extremal rays of $\overline{NE}(X)^G$ are not extremal in $\overline{NE}(X)$}\label{moribild}
\end{figure}

Since the cone $NE(X)$ is a $G$-invariant set it follows that its closure $\overline{NE}(X)$ is $G$-invariant. The subset of invariant elements in $\overline{NE}(X)$ is denoted by $\overline{NE}(X)^G$. \nomenclature{$\overline{NE}(X)^G$}{the intersection of $\overline{NE}(X)$ with the space of invariant numerical equivalence classes of 1-cycles}
\begin{remark}
$
\overline{NE}(X)^G = \overline{NE(X) \cap N_1(X)^G}=\overline{NE}(X) \cap N_1(X)^G .
$
\end{remark}
The subcone $\overline{NE}(X)^G$ of $\overline{NE}(X)$ is seen to inherit the geometric properties of $\overline{NE}(X)$ established by the cone theorem. 
Note however that the extremal rays of $\overline{NE}(X)^G$ are in general neither extremal in $\overline{NE}(X)$ (cf. Figure \ref{moribild}) nor generated by classes of curves but by classes of 1-cycles.
\begin{definition}
The extremal rays of $\overline{NE}(X)^G$ are called \emph{$G$-extremal rays} \index{extremal ray!$G$-extremal rays}.
\end{definition}
\begin{lemma}\label{Gextremalray} Let $G$ be a finite group and let $R$ be a $G$-extremal ray with $\mathcal{K}_X \cdot R<0$. Then there exists a rational curve $C_0$ such that $R$ is generated by the class of the 1-cycle $C = \sum_{g \in G} gC_0$.
\end{lemma}
\begin{proof}
Consider an $G$-extremal ray $R = \mathbb{R }_{\geq 0}[E]$ where $[E] \in \overline{NE}(X)^G \subset \overline{NE}(X)$. By the cone theorem (Theorem \ref{conethm}) it can be written as $[E] = [\sum_i a_i C_i] +[F]$,
where $ \mathcal{K}_X \cdot F \geq 0$, $a_i \geq 0$ and $C_i$ are rational curves. Let $|G|$ denote the order of $G$ and let $[GF] = G[F] = \sum_{g\in G} g[F]$. Since $g[E]=[E]$ for all $g\in G$ we can write
\[
|G| [E] = \sum_{g\in G} g[E] = \sum_{g\in G}([\sum_i a_igC_i] + g[F]) = \sum_i a_i G[C_i] + G[F].
\]
The element $[\sum a_i (GC_i)] + [GF]$ of the extremal ray $\mathbb{R }_{\geq 0}[E]$ is decomposed as the sum of two elements in $\overline{NE}(X)^G$. Since $R$ is extremal in $\overline{NE}(X)^G$ both must lie in $R=\mathbb{R }_{\geq 0}[E]$ .

Consider $[GF] \in R$. Since $g^*\mathcal{K}_X \equiv \mathcal{K}_X$ for all $g \in G$, we obtain
$$\mathcal{K}_X \cdot (GF) = \sum_{g\in G} \mathcal{K}_X \cdot (gF) = \sum_{g \in G}(g^*\mathcal{K}_X) \cdot F= |G|\mathcal{K}_X \cdot F \geq 0.$$ 
Since $\mathcal{K}_X \cdot R < 0$ by assumption this implies $[F]=0$ and $\mathbb{R }_{\geq 0}[E] = \mathbb{R }_{\geq 0}[\sum a_i (GCi)]$. Again using the fact that $R$ is extremal in $\overline{NE}(X)^G$, we conclude that each summand of $[\sum a_i (GC_i)]$ must be contained in $R=\mathbb{R }_{\geq 0}[E]$ and the extremal ray $\mathbb{R }_{\geq 0}[E]$ is therefore generated by $[GC_i]$ for some $C_i$ chosen such that $[GC_i] \neq 0$. This completes the proof of the lemma. 
\end{proof}
\section{The contraction theorem and minimal models of surfaces}\label{mmp}
In this section, we state the contraction theorem for smooth
projective surfaces. The proof of this theorem can be found e.g. in
\cite{kollarmori} and needs to be modified slightly in order to give an equivariant contraction theorem in the next section. 
\begin{definition}
Let $X$ be a smooth projective surface and let $F \subset \overline{NE}(X)$ be an extremal face. A morphism $\mathrm{cont}_F: X \to Z$ is called the \emph{contraction of $F$} \index{contraction morphism} if \nomenclature{$\mathrm{cont}_F$}{the contraction of an extremal face $F$}
\begin{itemize}
\item $(\mathrm{cont}_F)_*\mathcal{O}_X = \mathcal{O}_Z$ and 
\item $\mathrm{cont}_F(C) = \{\text{point}\}$ for an irreducible curve $C\subset X$ if and only if $[C] \in F$. 
\end{itemize}
\end{definition}
The following result is known as the contraction theorem (cf. Theorem 1.28 in \cite{kollarmori}).
\begin{theorem}\label{contractionthm} \index{contraction theorem}
Let $X$ be a smooth projective surface and $R \subset \overline{NE}(X)$ an extremal ray such that $\mathcal{K}_X \cdot R<0$. Then the contraction morphism $\mathrm{cont}_R: X \to Z$ exists and is one of the following types:
\begin{enumerate}
\item{$Z$ is a smooth surface and $X$ is obtained from $Z$ by blowing up a point. }
\item{$Z$ is a smooth curve and $\mathrm{cont}_R:X \to Z $ is a minimal ruled surface over $Z$.}
\item{$Z$ is a point and $-\mathcal{K}_X$ is ample. }
\end{enumerate}
\end{theorem}
The contraction theorem leads to the minimal model program for surfaces: Starting from $X$, if $\mathcal{K}_X$ is not nef, i.e, there exists an irreducible curve $C$ such that $\mathcal{K}_X C < 0$, then $\overline{NE}(X)_{\mathcal{K}_X <0}$ is nonempty and there exists an extremal ray $R$ which can be contracted. The contraction morphisms either gives a new surface $Z$ (in case 1) or provides a structure theorem for $X$ which is then either a minimal ruled surface over a smooth curve (case 2) or isomorphic to $\mathbb P^2$ (case 3). Note that the contraction theorem as stated above only implies $-\mathcal{K}_X$ ample in case 3. It can be shown that $X$ is in fact $\mathbb P^2$. This is omitted here since the statement does not transfer to the equivariant setup. In case 1, we can repeat the procedure if $K_Z$ is not nef. Since the Picard number drops with each blow down, this process terminates after a finite number of steps. The surface obtained from $X$ at the end of this program is called a \emph{minimal model} \index{minimal model} of $X$. 
\begin{remark}
 Let $E$ be a (-1)-curve on $X$. If $C$ is any irreducible curve on $X$, then $E \cdot C < 0$ if and only if $C =E$. It follows that $\overline{NE}(X) = \mathrm{span}(\mathbb R_{\geq 0}[E], \overline{NE}(X)_{E \geq 0})$. Now $E^2 = -1$ implies $ E \not\in \overline{NE}(X)_{E \geq 0}$ and $E$ is seen to generate an extremal ray in $\overline{NE}(X)$. By adjunction, $\mathcal K_X \cdot E < 0$. The contraction of the extremal ray $R = \mathbb R_{\geq 0}[E]$ is precisely the contraction of the (-1)-curve $E$. Conversely, each extremal contraction of type 1 above is the contraction of a (-1)-curve generating the extremal ray $R$. 
\end{remark}
\section{Equivariant contraction theorem and $G$-minimal models}
We state and prove an equivariant contraction theorem for smooth projective surfaces with finite groups of symmetries. Most steps in the proof are carried out in analogy to the proof of the standard contraction theorem.
\begin{definition}\label{equiContraction}
Let $G$ be a finite group, let $X$ be a smooth projective surface with $G$-action and let $R \subset \overline{NE}(X)^G$ be $G$-extremal ray. A morphism $\mathrm{cont}_R^G: X \to Z$ is called the \emph{$G$-equivariant contraction of $R$} \index{contraction morphism!equivariant contraction morphism} if
\begin{samepage}
\begin{itemize}
\item $\mathrm{cont}_R^G$ is equivariant with respect to $G$
\item $(\mathrm{cont}_R^G)_*\mathcal{O}_X = \mathcal{O}_Z$ and 
\item $\mathrm{cont}_R(C) = \{\text{point}\}$ for an irreducible curve $C\subset X$ if and only if $[GC] \in R$. 
\end{itemize}
\end{samepage}
\end{definition}
\begin{theorem}\index{contraction theorem!equivariant contraction theorem}
Let $G$ be a finite group, let $X$ be a smooth projective surface with $G$-action and let $R$ be a $G$-extremal ray with $ \mathcal{K}_X \cdot R <0$.  Then $R$ can be spanned by the class of $C= \sum_{g \in G} gC_0$ for a rational curve $C_0$, the equivariant contraction morphism $\mathrm{cont}_R^G: X \to Z$ exists and is one of the following three types:
\begin{enumerate}
\item $C^2 <0$ and $gC_0$ are smooth disjoint (-1)-curves. The map $\mathrm{cont}_R^G: X \to Z$ is the equivariant blow down of the disjoint union $\bigcup_{g \in G} gC_0$.
\item $C^2 =0$ and any connected component of $C$ is either irreducible or the union of two (-1)-curves intersecting transversally at a single point. The map $\mathrm{cont}_R^G: X \to Z$ defines an equivariant conic bundle over a smooth curve .
\item $C^2 >0$ , $N_1(X)^G = \mathbb{R }$ and $\mathcal{K}_X^{-1}$ is ample, i.e., $X$ is a Del Pezzo surface. The map $\mathrm{cont}_R^G: X \to Z$ is constant, $Z$ is a point. 
\end{enumerate}
\end{theorem}
\begin{proof}
Let $R$ be a $G$-extremal ray with $\mathcal{K}_X \cdot R <0$.
It follows from Lemma \ref{Gextremalray} that the ray $R$ can be
 spanned by a 1-cycle of the form $C = GC_0$ for a rational curve
 $C_0$. Let $n = |GC_0|$ and write $C = \sum_{i=1}^n C_i$
 where the $C_i$ correspond to $gC_0$ for some $g \in G$.
We distinguish three cases according to the sign of the self-intersection of $C$.

\medskip
\textbf{The case $C^2 <0$}

 We write $0 > C^2 = \sum_i C_i^2 + \sum_{i\neq j}C_i \cdot C_j$. Since $C_i$ are effective curves we know $C_i \cdot C_j \geq 0$ for all $i \neq j$. Since all curves $C_i$ have the same negative self-intersection and by assumption, $\mathcal{K}_X \cdot C = \sum_i \mathcal{K}_X \cdot C_i = n(\mathcal{K}_X \cdot C_i) <0$ the adjunction formula reads
$
2g(C_i) -2 = -2 = \mathcal{K}_X  \cdot C_i + C_i^2
$. Consequently, $\mathcal{K}_X \cdot C_i = -1$ and $C_i^2 = -1$.
 It remains to show that all $C_i$ are disjoint. We assume the contrary and without loss of generality 
$C_1 \cap C_2 \neq \emptyset$. Now $gC_1 \cap g C_2 \neq \emptyset $ for all $g \in G$ and $\sum_{i \neq j} C_i \cdot C_j \geq n$. This is however contrary to $0>C^2 = \sum_i C_i^2 + \sum_{i\neq j}C_i \cdot C_j = -n +  \sum_{i\neq j}C_i \cdot C_j$.

We let $\mathrm{cont}_R^G: X \to Z$ be the blow-down of $\bigcup_{g \in G} gC_0$ which is equivariant with respect to the induced action on $Z$ and fulfills $(\mathrm{cont}_R^G)_*\mathcal{O}_X = \mathcal{O}_Z$. If $D$ is an irreducible curve such that $\mathrm{cont}_R^G(D)= \{\text{point}\}$, then $D= gC_0$ for some $g \in G$. In particular, $GD= GC_0 =C$ and $[GD] \in R$. Conversely, if $[GD] \in R$ for some irreducible curve $D$, then $[GD]= \lambda [C]$ for some $\lambda \in \mathbb R _{\geq 0}$. Now $(GD)\cdot C = \lambda C^2 <0$. It follows that $D$ is an irreducible component of $C$. 

\medskip
\textbf{The case $C^2 >0$}

This case is treated in precisely the same way as the corresponding case in the standard contraction theorem.
Our aim is to show that $[C]$ is in the interior of $\overline{NE}(X)^G$. This is a consequence of the following lemma.
\begin{lemma}
Let $X$ be a projective surface and let $L$ be an ample line bundle on $X$. Then the set $Q = \{[E] \in N_1(X) \ |\ E^2 >0\}$ has two connected components
$Q^+=\{ [E] \in Q \ |L \cdot E >0\}$ and $Q^-=\{ [E] \in Q \ |L \cdot E <0\}$. Moreover, $Q^+ \subset \overline{NE}(X)$. 
\end{lemma}
This result follows from the Hodge Index Theorem (cf. Theorem IV.2.14 in \cite{BPV}) and the fact, that $E^2 >0$ implies that either $E$ or $-E$ is effective.
For a proof of this lemma, we refer the reader to Corollary 1.21 in \cite{kollarmori}. 

We consider an effective cycle $C = \sum C_i$ with $C^2 >0$. By the above lemma, $[C]$ is contained in $Q^+$ which is an open subset of $N_1(X)$ contained in $\overline{NE}(X)$. It follows that $[C]$ lies in the interior of $\overline{NE}(X)$. The $G$-extremal ray $R = \mathbb{R }_{\geq 0} [C]$ can only lie in the interior if $\overline{NE}(X)^G= R$. By assumption $\mathcal{K}_X \cdot R <0$, so that $\mathcal{K}_X$ is negative on $\overline{NE}(X)^G \backslash \{0\}$ and therefore on $\overline{NE}(X) \backslash \{0\}$. The anticanonical bundle $\mathcal{K}_X^{-1}$ is ample by Kleiman ampleness criterion and $X$ is a Del Pezzo surface. 

We can define a constant map $\mathrm{cont}_R^G$ mapping $X$ to a point $Z$ which is the equivariant contraction of $R = \overline{NE}(X)$ in the sense of Definition \ref{equiContraction}. 

\medskip
\textbf{The case $C^2 =0$}

Our aim is to show that for some $m >0$ the linear system $|mC|$ defines a conic bundle structure on $X$. The argument is seperated into a number of lemmata.
For the convenience of the reader, we include also the proofs of well-known preparatory lemmata which do not involve group actions. 
Recall that $\mathcal{O}(D)$ denotes the line bundle associated to the divisor $D$ on $X$.
\begin{lemma}
$H^2(X, \mathcal{O}(mC)) =0$ for $m \gg 0$.
\end{lemma}
\begin{proof}
By Serre's duality (cf. Theorem I.5.3 in \cite{BPV}) 
\[
h^2(X, \mathcal{O}(mC))= h^0(\mathcal{O}(-mC)\otimes \mathcal K_X).
\]
Since $C$ is an effective divisor on $X$, it follows that $h^0(\mathcal{O}(-mC)\otimes \mathcal K_X)=0$ for $m \gg 0$.
\end{proof}
\begin{lemma}
 For $m \gg 0$ the dimension $h^0(X,\mathcal{O}(mC))$ of $H^0(X,\mathcal{O}(mC))$ is at least two.
\end{lemma}
\begin{proof}
Let $m$ be such that $h^2(X, \mathcal{O}(mC)) =0$. For a line bundle $L$ on $X$ we denote by $\chi(L) = \sum_i(-1)^i h^i(X,L)$ the Euler characteristic of $L$. Using the theorem of Riemann-Roch (cf. Theorem V.1.6 in \cite{hartshorne}), 
\begin{align*}
h^0(X, \mathcal{O}(mC)) &\geq h^0(X, \mathcal{O}(mC))- h^1(X, \mathcal{O}(mC))\\
               & = h^0(X, \mathcal{O}(mC)) - h^1(X, \mathcal{O}(mC)) + h^2(X, \mathcal{O}(mC))\\
               & = \chi(\mathcal{O}(mC))\\
               & = \chi (\mathcal{O}) + \frac{1}{2}(\mathcal{O}(mC)\otimes\mathcal{K}_X^{-1})\cdot(mC)\\
               & \overset{C^2=0}{=} \chi (\mathcal{O}) - \frac{m}{2}\mathcal{K}_X \cdot C.
\end{align*}
Now $\mathcal{K}_XC<0$ implies the desired behaviour of $h^0(X, \mathcal{O}(mC))$. 
\end{proof}
For a divisor $D$ on $X$ we denote by $|D|$ \nomenclature{$|D|$}{the linear system of a divisor $D$}\index{linear system}the complete \emph{linear system of $D$}, i.e., the set of all effective divisors linearly equivalent to $D$. A point $ p \in X$ is called a \emph{base point}\index{base point} of $|D|$ if $p \in \mathrm{support}(C)$ for all $C \in |D|$.
\begin{lemma}
 There exists $m' >0$ such that the linear system $|m'C|$ is base point free.
\end{lemma}
\begin{proof}
We first exclude a positive dimensional set of base points. Let $m$ be chosen such that $h^0(X, \mathcal{O}(mC)) \geq 2$. We denote by $B$ the \emph{fixed part} of the linear system $|mC|$, i.e., 
the biggest divisor $B$ such that each $D \in |mC|$ can be decomposed as $D = B + E_D$ for some effective divisor $E_D$.
The support of $B$ is the union of all positive dimensional components of the set of base points of $|mC|$. We assume that $B$ is nonempty.
The choice of $m$ guarantees that $|mC|$ is not fixed, i.e., there exists $D \in |mC|$ with $D \neq B$.
Since $\mathrm{supp}(B) \subset \{s=0\}$ for all $ s \in \Gamma(X, \mathcal O(mC))$, each irreducible component of $\mathrm{supp}(B)$ is an irreducible component of $C$ and $G$-invariance of $C$ implies $G$-invariance of the fixed part of $|mC|$. It follows that $B=m_0 C$ for some $m_0 < m$.
Decomposing $|mC|$ into the fixed part $B = m_0C$  and the remaining \emph{free part} $|(m-m_0)C|$ shows that some multiple $|m'C|$ for $m' >0$  has no fixed components.
The linear system $|m'C|$ has no isolated base points since these would
correspond to isolated points of intersection of divisors linearly equivalent to $m'C$. Such intersections are excluded by $C^2=0$.
 \end{proof}
We consider the base point free linear system $|m'C|$ and the induced morphism 
\begin{align*}
\varphi =\varphi_{|m'C|}: &X \to \varphi(X) \subset \mathbb P(\Gamma(X,\mathcal{O}(m'C))^*)\\
&x \mapsto \{ s \in \Gamma(X,\mathcal{O}(m'C)) \, | \,  s(x)=0 \}
\end{align*}
Since $C$ is $G$-invariant, it follows that $\varphi$ is an equivariant map with respect to action of $G$ on $\mathbb P(\Gamma(X,\mathcal{O}(m'C))^*)$ induced by pullback of sections.

Let us study the fibers of $\varphi$. Let $z$ be a linear hyperplane in $\Gamma(X,\mathcal{O}(m'C))$. By definition, $\varphi^{-1}(z)= \bigcap _{s\in z}(s)_0$ where $(s)_0$ denotes the zero set of the section $s$. Since $(s)_0$ is linearly equivalent to $m'C$ and $C^2=0$, the intersection $\bigcap _{s\in z}(s)_0$ does not consist of isolated points but all $(s)_0$ with $s \in z$ have a common component. In particular, each fiber is one-dimensional.

Let $f: X \to Z$ be the Stein factorization of $\varphi: X \to \varphi(X)$. The space $Z$ is normal and 1-dimensional, i.e., $Z$ is a smooth curve.  Note that there is a $G$-action on the smooth curve $Z$ such that $f$ is equivariant. 
\begin{lemma}
The map $f: X \to Z$ defines an equivariant conic bundle\index{conic bundle}, i.e., an equivariant fibration with general fiber isomorphic to $\mathbb{P}_1$.
\end{lemma}
\begin{proof}
Let $F$ be a smooth fiber of $f$. By construction, $F$ is a component of $(s)_0$ for some $s \in \Gamma(X,\mathcal{O}(m'C))$. We can find an effective 1-cycle $D$ such that $(s)_0 = F+D$. Averaging over the group $G$ we obtain
$
\sum_{g\in G}gF + \sum_{g\in G}gD = \sum_{g\in G}g(s)_0
$. 
Recalling $(s)_0 \sim m'C $ and $[C] \in \overline{NE}(X)^G$ we deduce
\[
[\sum_{g\in G}gF + \sum_{g\in G}gD] = [\sum_{g\in G}g(s)_0] = m'[\sum_{g\in G}gC]= m |G| [C]
\]
showing that $[\sum_{g\in G}gF + \sum_{g\in G}gD]$ in contained in the $G$-extremal ray generated by $[C]$. Now by the definition of extremality $[\sum_{g\in G}gF] = \lambda [C] \in \mathbb{R }^{>0}[C]$ and therefore $\mathcal{K}_X \cdot (\sum_{g\in G}gF) <0$. This implies $\mathcal{K}_X F<0$.

In order to determine the self-intersection of $F$, we first observe $(\sum_{g\in G}gF)^2= \lambda^2 C^2 =0$. Since $F$ is a fiber of a $G$-equivariant fibration, we know that $\sum_{g\in G}gF = kF + kF_1 + \dots + kF_l$ where $F, F_1, \dots F_l$ are distinct fibers of $f$ and $k \in \mathbb N ^{>0}$. Now
$0=(\sum_{g\in G}gF)^2 = (l+1)k^2F^2$
shows $F^2=0$. The adjunction formula then implies $g(F)=0$ and $F$ is isomorphic to $\mathbb P_1$.
\end{proof}
The map $\mathrm{cont}_R^G:=f$ is equivariant and fulfills $f_* \mathcal{O}_X = \mathcal{O}_Z$ by Stein's factorization theorem. Let $D$ be an irreducible curve in $X$ such that $f$ maps $D$ to a point, i.e., $D$ is contained in a fiber of $f$. Going through the same arguments as above one checks that $[GD] \in R$. Conversely, if $D$ is an irreducible curve in $X$ such that $[GD] \in R$ it follows that $(GD) \cdot C=0$. If $D$ is not contracted by $f$, then $f(D) = Z$ and $D$ meets every fiber of $f$. In particular, $D \cdot C >0$, a contradiction. It follows that $D$ must be contracted by $f$.

This completes the proof of the equivariant contraction theorem.
\end{proof}
The singular fibers of the conic bundle in case 2 of the theorem above are characterized by the following lemma. 
\begin{lemma}\label{singular fibers of conic bundle}
Let $R = \mathbb R ^{>0} [C]$ be a $\mathcal K_X$-negative $G$-extremal ray with $C^2=0$. Let $\mathrm{cont}_R^G:=f: X \to Z$ be the equivariant contraction of $R$ defining a conic bundle structure on $X$. Then every singular fiber of $f$ is a union of two (-1)-curves intersecting transversally. 
\end{lemma}
\begin{proof}
Let $F$ be a singular fiber of $f$. The same argument as in the previous lemma yields that $\mathcal{K}_X \cdot F<0$ and $F^2 =0$. Since $F$ is connected, adjunction implies that the arithmetic genus of $F$ is zero and $\mathcal{K}_X \cdot F = -2$. It follows from the assumption on $F$ being singular that $F$ must be reducible. Let $F = \sum F_i$ be the decomposition into irreducible components. Now $g(F)=0$ implies $g(F_i)=0$ for all $i$. 

We apply the same argument as above to the component $F_i$ of $F$: after averaging over $G$ we deduce that $GF_i$ is in the $G$-extremal ray $R$ and $ \mathcal K _X \cdot F_i <0$.
Since $-2 = \mathcal K _X \cdot F = \sum \mathcal K _X \cdot F_i$, we may conclude that $F = F_1 +F_2$ and $F_i^2=-1$. 
The desired result follows.
\end{proof}
\section*{$G$-minimal models of surfaces}
Let $X$ be a surface with $G$-action such that $\mathcal{K}_X$ is not nef, i.e., $\overline{NE}(X)_{\mathcal{K}_X<0}$ is nonempty. 
\begin{lemma}
There exists a $G$-extremal ray $R$ such that $\mathcal{K}_X \cdot R<0$. 
\end{lemma}
\begin{proof}
Let $[C] \in \overline{NE}(X)_{\mathcal{K}_X<0}\neq \emptyset$ and consider $[GC] \in \overline{NE}(X)^G$. The $G$-orbit or $G$-average of a $\mathcal{K}_X$-negative effective curve is again $\mathcal{K}_X$-negative. It follows that $\overline{NE}(X)^G_{\mathcal{K}_X<0}$ is nonempty. 
Let $L$ be a $G$-invariant ample line bundle on $X$.
By the cone theorem,
for any $\varepsilon>0$
\[
\overline{NE}(X)^G = \overline{NE}(X)^G _{(\mathcal{K}_X+\varepsilon L) \geq 0} + \sum_{\text{finite}} \mathbb R_{\geq 0} G[C_i]. 
\]
where $\mathcal K _X \cdot C_i < 0$ for all $i$.
Since $\overline{NE}(X)^G_{\mathcal{K}_X<0}$ is nonempty, we may choose $\varepsilon>0$ such that $\overline{NE}(X)^G \neq \overline{NE}(X)^G _{(\mathcal{K}_X+\varepsilon L) \geq 0}$.
If the ray $R_1 = \mathbb R_{\geq 0} G[C_1]$ is not extremal in $\overline{NE}(X)^G$, then its generator $G[C_1]$ can be decomposed as a sum of elements of $\overline{NE}(X)^G$ not contained in $R_1$. It follows that 
\[
\overline{NE}(X)^G = \overline{NE}(X)^G _{(\mathcal{K}_X+\varepsilon L) \geq 0} + \underset{\text{finite}}{\sum_{i\neq 1}} \mathbb R_{\geq 0} G[C_i],
\]
i.e., the ray $R_1$ is superfluous in the formula.
By assumption $\overline{NE}(X)^G \neq \overline{NE}(X)^G _{(\mathcal{K}_X+\varepsilon L) \geq 0}$ and we may therefore not remove all rays $R_i$ from the formula and at least one ray $R_i = \mathbb R_{\geq 0} G[C_i]$ is $G$-extremal.
\end{proof}
We apply the equivariant contraction theorem to $X$:
In case 1 we obtain from $X$ a new surface $Z$ by blowing down a $G$-orbit of disjoint (-1)-curves. There is a canonically defined  holomorphic $G$-action on $Z$ such that the blow-down is equivariant. If $K_Z$ is not nef, we repeat the procedure which will stop after a finite number of steps. In case 2 we obtain an equivariant conic bundle structure on $X$. In case 3 we conclude that $X$ is a Del Pezzo surface with $G$-action. We call the $G$-surface obtained from $X$ at the end of this procedure a \emph{$G$-minimal model of $X$} \index{minimal model!$G$-minimal model}.  

As a special case, we consider a rational surface $X$ with $G$-action. 
Since the canonical bundle $\mathcal{K}_X$ of a rational surface $X$ is never nef (cf. Theorem VI.2.1 in \cite{BPV}), a $G$-minimal model of $X$ is an equivariant conic bundle over $Z$ or a Del Pezzo surface with $G$-action. Note that the base curve $Z$ must be rational: if $Z$ is not rational, one finds nonzero holomorphic one-forms on $Z$. Pulling these back to $X$ gives rise to nonzero holomorphic one-forms on the rational surface $X$, a contradiction.

This proves the well-known classification of $G$-minimal models of rational surfaces (cf. \cite{maninminimal}, \cite{isk}).
Although this classification does classically not rely on Mori theory, the proof given above is based on Mori's approach. We therefore refer to an equivariant reduction $Y \to Y_\mathrm{min}$ as an \emph{equivariant Mori reduction}\index{equivariant Mori reduction}.

In the following chapters we will apply the equivariant minimal model program to quotients of K3-surfaces by nonsymplectic automorphisms.
%
%
%
%
\chapter{Centralizers of antisymplectic involutions}\label{chapterlarge}
This chapter is dedicated to a rough classification of K3-surfaces with anti\-sym\-plec\-tic involutions centralized by large groups of symplectic transformations (Theorem \ref{roughclassi}). 

We consider a K3-surface $X$ with an action of a finite group $G
\times C_2 < \mathrm{Aut}(X)$ and assume that the action of $G$ is by
symplectic transformations whereas $C_2$ is generated by an
antisymplectic involution $\sigma$ centralizing $G$. Furthermore, we assume that $\mathrm{Fix}_X(\sigma) \neq \emptyset$. Let $\pi: X \to X/ \sigma = Y$ denote the quotient map. The quotient surface $Y$ is a smooth rational $G$-surface to which we apply the equivariant minimal model program developed in the previous chapter. A $G$-minimal model of $Y$ can either be a Del Pezzo surface or an equivariant conic bundle over $\mathbb P_1$. In the later case, the possibilities for $G$ are limited by the classification of finite groups with an effective action on $\mathbb P_1$
\begin{remark}\label{autP_1}
The classification of finite subgroups of $\mathrm{SU}(2, \mathbb C)$ (or $\mathrm{SO}(3, \mathbb R)$) yields the following list of finite groups with an effective action on $\mathbb P_1$:
\begin{samepage}
\nomenclature{$D_{2n}$}{the dihedral group of order $2n$}
\nomenclature{$T_{12}$}{the tetrahedral group}
\nomenclature{$O_{24}$}{the octahedral group}
\nomenclature{$I_{60}$}{the icosahedral group}
\begin{itemize}
\item
cyclic groups $C_n$
\item
dihedral groups $D_{2n}$
\item
the tetrahedral group $T_{12} \cong A_4$
\item
the octahedral group $O_{24} \cong S_4$
\item
the icosahedral group $I_{60} \cong A_5$
\end{itemize}
\end{samepage}
If $G$ is any finite group acting on a space $X$,
we refer to the number of elements in an orbit $G.x = \{ g.x \, | \, g \in G\}$ as the \emph{length of the $G$-orbit $G.x$}\index{length of an orbit}.
Note that the length of a $T_{12}$-orbit in $\mathbb P_1$ is at least four, the length of an $O_{24}$-orbit in $\mathbb P_1$ is at least six, and the length of an $I_{60}$-orbit in $\mathbb P_1$ is at least twelve.
\end{remark}
\begin{lemma}\label{conicbundle}
If a $G$-minimal model $Y_\mathrm{min}$ of $Y$ is an equivariant conic bundle, then $|G| \leq 96$.
\end{lemma}
\begin{proof}
Let $\varphi: Y_\mathrm{min} \to \mathbb{P}_1$ be an equivariant conic bundle structure on $Y_\mathrm{min}$. By definition, the general fiber of $\varphi$ is isomorphic to $\mathbb{P}_1$. We consider the induced action of $G$ on the base $\mathbb P_1$. If this action is effective, then $G$ is among the groups specified in the remark above. Since the maximal order of an element in $G$ is eight (cf. Remark \ref{order of symp aut}), it follows that the order $G$ is bounded by 60.

If the action of $G$ on the base $\mathbb P_1$ is not effective, every element $n$ of the ineffectivity $N < G$ has two fixed points in the general fiber. This gives rise to a positive-dimensional $n$-fixed point set in $Y_\mathrm{min}$ and $Y$. A symplectic automorphism however has only isolated fixed points. It follows that the action of $n$ on $X$ coincides with the action of $\sigma $ on $\pi^{-1}(\mathrm{Fix}_Y(N))$. In particular, the order of $n$ is two. Since $N$ acts effectively on the general fiber, it follows that $N$ is isomorphic to either $C_2$ or $C_2 \times C_2$. 

If $G/N$ is isomorphic to the icosahedral group $I_{60}= A_5$, then $G$ fits into the exact sequence $ 1 \to N \to G \to A_5 \to 1$ for $N=C_2$ or $C_2 \times C_2$. Let $\eta$ be an element of order five inside $A_5$. One can find an element $\xi$ of order five in $G$ which is mapped to $\eta$. Since neither $C_2$ nor $C_2 \times C_2$ has automorphisms of order five it follows that $\xi$ centralizes the normal subgroup $N$. In particular, there is a subgroup $C_2 \times C_5 \cong C_{10}$ in $G$ which is contrary to the assumption that $G$ is a group of symplectic transformations and therefore its elements have order at most eight.

If $G/N$ is cyclic or dihedral, we again use the fact that the order of elements in $G$ is bounded by $8$ and conclude $|G/N| \leq 16$. It follows that the maximal possible order of $G/N$ is 24. Using $|N| \leq 4$ we obtain $|G| \leq 96$.
\end{proof}
If $| G | >96$, the lemma above allows us to restrict our classification to the case where a $G$-minimal model $Y_\mathrm{min}$ of $Y$ is a Del Pezzo surface. The next section is devoted to a brief introduction to Del Pezzo surfaces and their automorphisms groups. 
\section{Del Pezzo surfaces}
A \emph{Del Pezzo surface}\index{Del Pezzo surface} is a smooth surface $Z$ such that the
anticanonical bundle $\mathcal K_Z^{-1} = \mathcal O _Z(-K_Z)$ is ample. 
The self-intersection number of the
canonical divisor\index{canonical divisor}\nomenclature{$K_X$}{the canonical divisor of $X$} $d:= K_Z^2$ is referred to as the \emph{degree} \index{degree of a Del Pezzo surface} of the Del
Pezzo surface and  $1 \leq d
\leq 9$ (cf. Theorem 24.3 in \cite{manin}). 
\begin{example}
Let $Z = \{f_3=0\} \subset \mathbb P_3$ be a smooth cubic surface. The anticanonical bundle $\mathcal K_Z^{-1}$ of $Z$ is given by the restriction of the hyperplane bundle $\mathcal O_{\mathbb P_3}(1)$ to $Z$ and therefore ample. 
\end{example}
As a consequence of the adjunction formula, an
irreducible curve with negative self-inter\-sec\-tion on a Del Pezzo
surface is a (-1)-curve. The following theorem (cf. Theorem 24.4 in \cite{manin}) gives a classification
of Del Pezzo surfaces according to their degree.
\begin{theorem}\label{classiDelPezzo}
Let $Z$ be a Del Pezzo surface of degree $d$. 
\begin{itemize}
\item
If $d=9$, then $Z$ is isomorphic to $\mathbb P_2$.
\item
If $d=8$, then $Z$ is isomorphic to either $\mathbb P_1 \times \mathbb
P_1$ or the blow-up of $\mathbb P_2$ in one point.
\item
If $1 \leq d \leq 7$, then $Z$ is isomorphic to the blow-up of
$\mathbb P_2$ in $9-d$ points in general position, i.e., no three
points lie on one line and no six points lie on one conic. 
\end{itemize}
\end{theorem}
In our later considerations of Del Pezzo surfaces Table \ref{minus one curves} below
(cf. Theorem 26.2 in \cite{manin})
specifying the number of (-1)-curves on a Del Pezzo surface of degree
$d$ will be very useful.
\begin{table}[H]
\centering
\begin{tabular}{l|ccccccc}
degree $d$ & 1&2&3&4&5&6&7 \\
\hline
number of (-1)-curves & 240 & 56 & 27 & 16 & 10 & 6 & 3
\end{tabular}
\caption{(-1)-curves on Del Pezzo surfaces}\label{minus one curves}
\end{table}
\begin{example}
Let $Z$ be a Del Pezzo surface of degree 5. It follows from the theorem above that $Z$ is isomorphic to the blow-up of $\mathbb P_2$ in four points $p_1, \dots, p_4$ 
 in general position. We denote by $E_i$ the preimage of $p_i$ is $Z$. Let $L_{ij}$ denote the line in $\mathbb P_2$ joining $p_i$ and $p_j$ and note that there are precisely six lines of this type. The proper transform of $L_{ij}$ is a (-1)-curve in $Z$. We have thereby specified all ten (-1)-curves in $Z$. Their incidence graph is known as the \emph{Petersen graph} \index{Petersen graph}. 
\end{example}
The following theorem summarizes properties of the
anticanonical map, i.e., the map associated to the linear system $|-K_Z|$ of
the anticanonical divisor (Theorem 24.5 in \cite{manin} and Theorem
8.3.2 in \cite{dolgachev})
\begin{theorem}\label{antican models of del pezzo}
Let $Z$ be a Del Pezzo surface of degree $d$. If $d \geq 3 $, then
$\mathcal K_Z^{-1}$ is very ample and the anticanonical map is a holomorphic
embedding of $Z$ into $\mathbb P_d$ such that the image of $Z$ in $\mathbb P_d$
is of degree $d$.

If $d=2$, then the anticanonical map is a
holomorphic degree two cover $\varphi: Z \to \mathbb P_2$ branched
along a smooth quartic curve. 

If $d=1$, then the linear system
$|-K_Z|$ has exactly one base point $p$. Let $ Z' \to Z$ be the blow-up
of $p$. Then the pull-back of $-K_Z$ to $Z'$ defines an elliptic
fibration $f: Z' \to \mathbb P_1$. The linear system $|-2K_Z|$ defines a finite map of degree two onto a quadric cone $Q$ in $\mathbb P_3$. Its branch locus is given by the intersection of $Q$ with a cubic surface. 
\end{theorem}
Our understanding of Del Pezzo surfaces as surfaces obtained by
blowing-up points in $\mathbb P_ 2$ in general position or as degree
$d$ subvarieties of $\mathbb P_d$ enables us the decide whether
certain finite groups $G$ can occur as subgroups of the automorphisms
group $\mathrm{Aut}(Z)$ of a Del Pezzo surface $Z$. 
\begin{example}\label{DelPezzoC3C7}
Consider the semi-direct product $G= C_3 \ltimes C_7$ where the action of
$C_3$ on $C_7$ is defined by the embedding of $C_3$ into
$\mathrm{Aut}(C_7) \cong C_6$. The group $G$ is a maximal subgroup of the simple group $L_2(7)$ which is discussed below.
Let $Z$ be a  Del Pezzo surface of degree $d$ with an effective action of $G$. 
Since $G$ does not admit a two-dimensional representation, it follows
that $G$ does not have fixed points in $Z$. In particular, $d \neq
1$. For the same reason, $Z$ is not the blow-up of $\mathbb P_2$ in
one or two points.  
Since there is no nontrivial homomorphisms $G \to C_2$ and no
injective homomorphism $ G \to \mathrm{PSL}(2, \mathbb C)$ it follows
that $G \not\hookrightarrow \mathrm{Aut}(\mathbb P_1 \times \mathbb
P_1) = (\mathrm{PSL}_2(\mathbb C) \times \mathrm{PSL}_2(\mathbb C)) \rtimes C_2$.
\end{example}
In many cases it can be useful to consider possible actions of a
finite group $G$ on the union of (-1)-curves on a Del Pezzo surfaces. 
\begin{example}\label{DelPezzoL2(7)}
We consider $G= L_2(7)$, the simple group of order 168. Its maximal
subgroups are $C_3 \ltimes C_7$ and $S_4$. Assume $G$ acts effectively on a Del
Pezzo surface $Z$ of degree $d$. Since $L_2(7)$ does not
stabilize any smooth rational curve, the $G$-orbit of a (-1)-curve $E
\subset Z$ consists of 7, 8, 14, 24 or more curves. It now follows from
Table \ref{minus one curves} that $d \neq 3, 5, 6$. 

If $d=4$, then the
union of (-1)-curves on $Z$ would consist of two $G$-orbits of length 8. In
particular, $\mathrm{Stab}_G(E)\cong C_3 \ltimes C_7$ for any
(-1)-curve $E \subset Z$. 
Blowing down $E$ to a point $p \in Z'$ induces an action of $C_3 \ltimes C_7$ on $Z'$ fixing $p$. Since $C_3 \ltimes C_7$ does not admit a two-dimensional representation, it follows that the normal subgroup $C_7$ acts trivially on $Z'$ and therefore on $Z$. This is a contradiction.

Using the result of the previous example, it follows that $Z$ is
either a Del Pezzo surface of degree 2 or isomorphic to $\mathbb P_2
$. Both cases will play a role in our discussion of K3-surfaces with
an action of $L_2(7)$.
\end{example}
\begin{example}
Let be the Del Pezzo surface obtained by blowing up one point $p$ in $\mathbb P_2$. Then its automorphims group is the subgroup of $\mathrm{Aut}(\mathbb P_2)$ fixing the point $p$. Similarly, if $Z$ is the Del Pezzo surface obtained by blowing up two points $p,q$ in $\mathbb P_2$, then $\mathrm{Aut}(Z) = G \rtimes C_2$ where $G$ is the subgroup of $\mathrm{Aut}(\mathbb P_2)$ fixing the two points $p, q$ and $C_2$ acts by switching the exceptional curves $E_p,E_q$.
\end{example}
In the previous chapter we have shown that Del Pezzo surfaces can occur as equivariant minimal models. It should be remarked that the blow-up of $\mathbb P_2$ in one or two points is never equivariantly minimal: Let $Z$ be the surface obtained by blowing up one or two points in $\mathbb P_2$. Then $Z$ contains an $\mathrm{Aut}(Z)$-invariant (-1)-curve, namely the curve $E_p$ in the first case and the proper transform of the line joining $p$ and $q$ in the second case. This curve can always be blown down equivariantly. Using the language of equivariant Mori theory introduced in the previous chapter, the $\mathrm{Aut}(Z)$-invariant (-1)-curve spans a $\mathrm{Aut}(Z)$-extremal ray $R$ of the cone of invariant curves $\overline{NE}(X)^{\mathrm{Aut}(Z)}$ with $\mathcal K_Z \cdot R <0$. Its contraction defines an $\mathrm{Aut}(Z)$-equivariant map to $\mathbb P_2$. In particular, $Z$ is not equivariantly minimal. 
\begin{remark}
A complete classification of automorphisms groups of Del Pezzo surfaces
can be found in \cite{dolgachev}. 
\end{remark}
\section{Branch curves and Mori fibers} \label{branch curves mori fibers}
We return to the initial setup where $X$ is a K3-surface with an action of $G \times \langle \sigma \rangle$ and $\pi: X \to X/\sigma =Y$ denotes the quotient map, and fix an equivariant Mori reduction $M: Y \to Y_\mathrm{min}$.

A rational curve $E \subset Y$ is called a \emph{Mori fiber} \index{Mori fiber} if it is contracted in some step of the equivariant Mori reduction $Y \to Y_\mathrm{min}$. The set of all Mori fibers is denoted by $\mathcal E$. Its cardinality $|\mathcal E|$ is denoted by $m$. We let $n$ denote the total number of rational curves in $\mathrm{Fix}_X(\sigma)$. 
\begin{lemma}
The total number $m$ of Mori fibers in $Y$ is bounded by $m \leq n+ 12 - e(Y_\mathrm{min}) \leq n+9$. 
\end{lemma}
\begin{proof}\label{moribound}\index{Euler characteristic formula!}
Recall that $\mathrm{Fix}_X(\sigma)$ is a disjoint union of smooth curves.
We choose a triangulation of $\mathrm{Fix}_X(\sigma)$ and extend it to a triangulation of the surface $X$. The topological Euler characteristic of the double cover is 
\begin{align*}
e(X) = 24       &= 2e(Y) - \sum_{C \subset \mathrm{Fix}_X(\sigma)} e(C)\\
                &= 2e(Y) - \sum_{C \subset \mathrm{Fix}_X(\sigma)} (2-2g(C))\\
                &= 2e(Y) - 2n + \underset{\ g(C) \geq 1}{\sum_{C \subset \mathrm{Fix}_X(\sigma)}} (2g(C)-2)\\
                & \geq 2e(Y) -2n\\
                &= 2(e(Y_\mathrm{min}) +m) -2n
\end{align*}
This yields $m \leq n+12 - e(Y_\mathrm{min})$, and $e(Y_{min}) \geq 3$ completes the proof of the lemma.
\end{proof}
Let $R:= \mathrm{Fix}_X(\sigma) \subset X$ denote the ramification
locus of $\pi$ and let $B:=\pi(R) \subset Y$ be its branch locus. In the following, we repeatedly use the fact that for a finite proper surjective holomorphic map of complex manifolds (spaces) $\pi: X \to Y$ of degree $d$, the intersection number of pullback divisors fulfills $(\pi^*D_1 \cdot \pi^*D_2) =  d(D_1 \cdot D_2)$.
\begin{lemma}\label{preimage of E in X}
Let $E \in \mathcal{E}$ be a Mori fiber such that $E \not\subset B$ and $|E\cap B| \geq 2$ or $E \cdot B \geq 3$.
 Then $E^2 = -1$ and  $\pi^{-1}(E)$
is a smooth rational curve in $X$. Furthermore, $E \cdot B = |E \cap B| =2$. 
\end{lemma}
\begin{proof}

Let $k < 0$ denote self-intersection number of $E$. By the remark above, the divisor
$\pi^{-1}(E) = \pi^* E$ has self-inter\-section $2k$. Assume that $\pi^{-1}(E)$ is
reducible and let $\tilde E_1, \tilde E_2$ denote its irreducible components.  They are rational
and therefore, by adjunction on the K3-surface $X$, have self-intersection number $-2$.  Write
$$
0 > 2k = (\pi^{-1}(E))^2 = \tilde E_1^2 + \tilde E_2^2 + 2 (\tilde E_1 \cdot \tilde E_2) = -4 + 2 (\tilde E_1 \cdot \tilde E_2).
$$
Since $\tilde E_1$ and $\tilde E_2$ intersect at points in the preimage of $E \cap B$, we obtain $\tilde E_1 \cdot \tilde E_2 \geq 2$, a contradiction. It follows that $\pi^{-1}(E)$ is irreducible. Consequently, $k=-1$ and $\pi^{-1}(E)$ is a smooth rational curve with two $\sigma$-fixed points .
\end{proof}
\begin{remark}\label{self-int of Mori-fibers}
Let $E \in \mathcal E$ be a Mori fiber. 
\begin{itemize}
 \item 
If $E \subset B$, then $E$ is the image of a rational curve in $X$ and $E^2 = -4$. (cf. Corollary \ref{minusfour} below). 
\item
If $E \not\subset B$ and $\pi^{-1}(E)$ is irreducible, then $2E^2 = (\pi^{-1}(E))^2 <0$. Adjunction on $X$ implies that $(\pi^{-1}(E))^2=-2$ and that $\pi^{-1}(E)$ is a smooth rational curve in $X$. 
The action of $\sigma$ has two fixed points on $\pi^{-1}(E)$ and 
the restricted degree two map $\pi|_{\pi^{-1}(E)}: \pi^{-1}(E) \to E$ is necessarily branched, i.e., $E \cap B \neq \emptyset$.
\item
If $E \not\subset B$ and $\pi^{-1}(E)= \tilde E_1 + \tilde E_2$ is reducible, then 
\[
2E^2 = \underset{\geq -2}{\underbrace{\tilde E_1^2}} + \underset{\geq 0}{\underbrace{2(\tilde E_1 \cdot \tilde E_2)}} + \underset{\geq -2}{\underbrace{\tilde E_2^2}} \geq -4.
\]
In particular, $E^2 \in \{-1,-2\}$.
\begin{itemize}
\item
If $E^2=-1$, then $\tilde E_1 \cdot  \tilde E_2 =1$ and $E\cap B \neq \emptyset$. 
\item
If $E^2=-2$, then $\tilde E_1\cdot \tilde E_2 =0$ and $E\cap B = \emptyset$. 
\end{itemize}
\end{itemize}
In summary, a Mori fiber $E \not\subset B$ has self-intersection -1 if and only if $E\cap B \neq \emptyset$ and self-intersection -2 if and only if $E\cap B = \emptyset$. A Mori fiber $E$ has self-intersection -4 if and only if $E \subset B$.

More generally, any (-1)-curve $E$ on $Y$ meets $B$ in either one or two points. If $|E \cap B |=1$, then $\pi^{-1}(E) = E_1  \cup E_2$ is reducible. If $|E \cap B |=2$, then $\pi^{-1}(E)$ is irreducible and meets $\mathrm{Fix}_X(\sigma)=R = \pi^{-1}(B)$ in two points.
\end{remark}
\begin {proposition} \label {at most two}
Every Mori fiber $E \in \mathcal{E}$, $ E \not\subset B$ meets the branch locus $B$ in at
most two points. If $E$ and $B$ are tangent at $p$, then 
$E\cap B = \{p\}$ and $(E \cdot B)_p =2$.
\end {proposition}
\begin{proof} 
Let $E \in \mathcal E$, $ E \not\subset B$ and assume $|E\cap B| \geq 2$ or $E \cdot B \geq
3$. Then by the lemma above, $\tilde E = \pi^{-1}(E)$ is a smooth
rational curve in $X$. Since $\tilde E \not\subset \mathrm{Fix}_X(\sigma)$, the involution $\sigma$ has exactly two fixed points on $\tilde E$ showing $|E\cap B| = 2$. It remains to show that the intersection is transversal.

To see this, let $N_{\tilde{E}}$ denote the normal bundle of $\tilde{E}$ in $X$. We
consider the induced action of $\sigma$ on $N_{\tilde{E}}$ by a bundle
automorphism. Using an
equivariant tubular neighbourhood theorem we may equivariantly identify a
neighbourhood of $\tilde E$ in $X$ with $N_{\tilde E}$ via a 
$C^{\infty}$-diffeomorphism. The $\sigma$-fixed point curves
intersecting $\tilde{E}$ map to curves of $\sigma$-fixed points in
$N_{\tilde{E}}$ intersecting the zero-section and vice versa. 
Let $D$ be a curve of $\sigma$-fixed point in $N_{\tilde{E}}$. If $D$ is
not a fiber of $N_{\tilde E}$, it follows that $\sigma$ stabilizes all 
fibers intersecting
$D$ and the induced action of $\sigma$ on the base must be trivial, a contradiction.
It follows that the $\sigma$-fixed point curves correspond to fibers of 
$N_{\tilde{E}}$, and $E$ and $B$ meet transversally.

By negation of the implication above, if $E$ and $B$ are tangent at $p$, then $|E \cap B |=1$ and $E \cdot B=2$.
\end{proof}
\subsection{Rational branch curves}
In this section we find conditions on $G$, in particular conditions on the order of $G$, guaranteeing the absence of rational curves in $\mathrm{Fix}_X(\sigma)$. 
\begin{lemma}\label{selfintbranch}
Let $\pi:X \to Y$ be a cyclic degree two cover of surfaces and let $C \subset X$ be a smooth curve contained in the ramification locus of $\pi$. Then the image of $C$ in $Y$ has self-intersection $(\pi(C))^2= 2 C^2$.
\end{lemma}
\begin{proof}
We recall that the intersection of pullback divisors fulfills $\pi^*D_1 \cdot  \pi^*D_2 =  2(D_1 \cdot D_2)$. In the setup of the lemma, $(\pi^* \pi(C))^2 = 2 (\pi(C))^2$. Now $\pi^* \pi(C) \sim 2C$ implies the desired result. 
\end{proof}
Note that the lemma above can also be proved by considering the normal bundle $N_C$ of $C$ and the induced action of $\sigma$ on it. The normal bundle $N_{\pi(C)}$ is isomorphic to $N_C^2$. Since the self-intersection of a curve is the degree of the normal bundle restricted to the curve, the formula follows. 
\begin{corollary}\label{minusfour}
Let $X$ be a K3-surface and let $\pi: X \to Y$ be a cyclic degree two cover.
Then a rational branch curve of $\pi $ has self-intersection -4.
\end{corollary}
\begin{proof}
Let $C$ be a rational curve on the K3-surface $X$. Then by adjunction $C^2 =-2$ and the image $\pi(C)$ in $Y$ is a (-4)-curve by Lemma \ref{selfintbranch} above.
\end{proof}
On a Del Pezzo surface a curve with negative self-intersection necessarily has self-intersection -1. So if $Y_\mathrm{min}$ is a Del Pezzo surface, all rational branch curves of $\pi$, which have self-intersection -4 by Corollary \ref{minusfour}, need to be modified by the Mori reduction when passing to $Y_{\mathrm{min}}$ and therefore have nonempty intersection with the union of Mori fibers. 

An important tool in the study of rational branch curves is provided by the following lemma which describes the behaviour of self-intersection numbers under monoidal transformations.
\begin{lemma}\label{selfintblowdown}
Let $\tilde X$ and $X$ be smooth projective surfaces and let $b: \tilde X \to X$ be the blow-down of a (-1)-curve $ E \subset \tilde X$. For a curve $B \subset \tilde X$ having no common component with $ E$ the self-intersection of its image in $X$ is given by
\[
b( B)^2 =  B^2 + ( E \cdot B)^2.
\]
\end{lemma}
\begin{proof}
We may choose an ample divisor $H$ in $X$ with $p \not\in \mathrm{supp}(H)$ and $D$ linearly equivalent to $ b(B) +H$ such that $p \not\in \mathrm{supp}(D)$. 
Since $b$ is biholomorphic away from $p$, we know 
\[
(b(B) +H)^2 = D^2 =(b^*D)^2 =  (b^*((b(B) +H))^2.
\]
Using $(b^*H)^2 = H^2$ and $b^*(b(B))\cdot b^*H = b(B) \cdot H$ we find 
$b(B)^2 = (b^*B)^2$.
Now $b^*B = B + \mu E$ where $\mu$ denotes the multiplicity of the point $p \in b(B)$. This multiplicity equals the intersection multiplicity $E \cdot  B$. Therefore,
\[
b( B)^2  = (b^*B)^2= (B + \mu  E)^2 =  B^2 + 2 \mu^2 -\mu^2=  B^2 + \mu^2.
\]
and the lemma follows.
\end{proof}
We denote by $\mathcal{C}$ the set of rational branch curves of
$\pi$. The total number $|\mathcal{C}|$ of these curves is denoted
by $n$. The union of all Mori fibers not contained in the branch locus $B$ is denoted by $\bigcup E_i$.

Let $\mathcal{C}_{\geq k}= \{ C \in \mathcal C \, | \, |C \cap  \bigcup E_i | \geq k\} $ be the set of those rational branch curves $C$ which meet $\bigcup E_i$ in at least $k$ distinct points and let $| \mathcal{C}_{\geq k}| = r_k$. 
We let
$\mathcal{E}_{\geq k}$ denote the set of Mori fibers $E \not\subset B$ which
intersect some $C$ in $\mathcal{C}_{\geq k}$ and define 
\[
P_k = \{ (p,E) \, |\, p \in C \cap E,\, E \in \mathcal{E}_{\geq k},\, C
\in \mathcal{C}_{\geq k}\} \subseteq Y \times \mathcal{E}_{\geq k}
\]
and the projection map $\mathrm{pr}_k: P_k \to \mathcal{E}_{\geq k}$ mapping
$(p,E)$ to $E$. This map is surjective by definition of
$\mathcal{E}_{\geq k}$ and its fibers consist of $\leq 2$ points by
Proposition \ref{at most two}. Using $|P_k| \geq k r_k$ we see
\begin{equation}\label{boundforE_k}
|\mathcal{E}_{\geq k}| \geq \frac{k}{2}r_k.
\end{equation} 
Let $N$ be the largest positive integer such that $\mathcal{C}_{\geq N} =
\mathcal{C}$, i.e., each rational ramification curve is intersected at least $N$
times by Mori fibers. 
A curve $C \in \mathcal{C}$ which is intersected precisely $N$ times by
Mori fibers is referred to as a \emph{minimizing curve}. In the
following, let $C$ be a minimizing curve and let $ H =
\mathrm{Stab}_G(C) < G$ be the stabilizer of $C$ in $G$.
\begin{remark}
The index of
$H$ in $G$ is bounded by $n= r_N$.
\end{remark} 
\subsection*{Bounds for $n$}
A smooth rational curve
on a K3-surface has self-intersection -2 and all curves in
$\mathrm{Fix}_X(\sigma)$ are disjoint. Therefore, the rational curves
in  $\mathrm{Fix}_X(\sigma)$ generate a sublattice of
$\mathrm{Pic}(X)$ of signature $(0,n)$. It follows immediately that $n
\leq 19$.

A sharper bound $n \leq 16$ for the number of disjoint (-2)-curves on a K3-surface has been
obtained by Nikulin \cite{NikulinKummer} and the following optimal bound in our setup is due to Zhang \cite{ZhangInvolutions}, Theorem 3.
\begin{proposition}
The total number of connected curves in the fixed point set of an antisymplectic
involution on a K3-surface is bounded by 10.
\end{proposition}
\begin{corollary}\label{atmostten}
The number $n$ of rational curves in $\mathrm{Fix}_X(\sigma)$ is at
most 10. If $n=10$, then $\mathrm{Fix}_X(\sigma)$ is a union of
rational curves. 
\end{corollary}
In the following, we use Zhang's bound $n \leq 10$. Note, however, that all results can likewise be obtained by using the weakest bound $n \leq 19$. 

For $N \geq 4 $ Zhang's bound can be sharpened using the notion of
Mori fibers and minimizing curves.
\begin{lemma}\label{boundforn}
$\frac{N}{2}n  \leq  n +12 - e(Y_\mathrm{min}) \leq n+9$.
\end{lemma}
\begin{proof}
Using Lemma \ref{moribound} and inequality \eqref{boundforE_k}
$
\frac{N}{2}n = \frac{N}{2}r_N \leq |\mathcal{E}_{\geq N}|\leq
|\mathcal{E}| \leq  n +12 - e(Y_\mathrm{min})  \leq n+9
$.
\end{proof}
In the following we consider the stabilizer $H$ of a minimizing curve
$C$ and using the above bounds for $n$, we obtain bounds for $|G|$.
\subsection*{A bound for $|G|$}
\begin{proposition}
Let $X$ be a K3-surface with an action of a finite group $G \times \langle \sigma
\rangle$ such that $ G < \mathrm{Aut}_\mathrm{symp}(X)$ and $\sigma$ is
an antisymplectic involution with fixed points. If $| G | > 108$,
then $\mathrm{Fix}_X(\sigma)$ contains no rational curves.
\end{proposition}
\begin{proof}
 Assume that $\mathrm{Fix}_X(\sigma)$ contains rational curves and
 consider a minimizing curve $C \subset B$ and its stabilizer $\mathrm{Stab}_G(C) =:H$. Since a symplectic automorphism on $X$ does not admit a one-dimensional set of fixed points, it follows that the action of $H$ on $C$ is effective and $H$ is among the groups discussed in Remark \ref{autP_1}. 
We recall the possible lengths of $H$-orbits in $C$: the length of an orbit of a dihedral group is at least two, the length of a $T_{12}$-orbit in $\mathbb P_1$ is at least four, the length of an $O_{24}$-orbit in $\mathbb P_1$ is at least six, and the length of an $I_{60}$-orbit in $\mathbb P_1$ is at least twelve.

 Let $Y_\mathrm{min}$ be a $G$-minimal model of $X/\sigma =Y$. Recall that by Lemma \ref{conicbundle} $Y_\mathrm{min}$ is a Del Pezzo surface. Each rational branch curve is a (-4)-curve in $Y$. Since its image in $Y_\mathrm{min}$ has self-intersection $\geq -1$, it must intersect Mori fibers.
\begin{itemize}
 \item 
If $N=1$, i.e., the rational curve $C$ meets the union of Mori fibers in exactly one point $p$, then $p$ is a fixed point of the $H$-action on $C$. In particular, $H$ is a cyclic group $C_k$. By Remark \ref{order of symp aut} $k \leq 8 $. Since the index of $H$ in $G$ is bounded by $n \leq 10$, it follows that $|G | \leq 80$.
\item
If $N=2$, then $H$ is either a cyclic or a dihedral group. By Proposition 3.10 in \cite{mukai} the maximal order of a dihedral group of symplectic automorphisms on a K3-surface is 12. We first assume $H \cong D_{2m}$ and that the $G$-orbit $G.C$ of the rational branch curve $C$ has the maximal length $n=|G.C|=10$, i.e., $B = G \cdot C$. 
Each curve in $G \cdot C$ meets the union of Mori fibers in precisely two points forming an $D_{2m}$-orbit. If a Mori fiber $E_C$ meets the curve $C$ twice, then it follows from Proposition \ref{at most two} that $E$ meets no other curve in $B$. The contraction of $E$ transforms $C$ into a singular curve of self-intersection zero. The Del Pezzo surface $Y_\mathrm{min}$ does however not admit a curve of this type. It follows, that $E$ meets a Mori fiber $E'$ which is contracted in a later step of the Mori reduction and meets no other Mori fiber than $E'$.
The described configuration $G E \cup G E'$ requires a total number of at least 20 Mori fibers and therefore contradicts Lemma  \ref{moribound}. If $C$ meets two distinct Mori fibers $E_1, E_2$, each of these two can meet at most one further curve in $B$. The contraction of $E_1$ and $E_2$ transforms $C$ into a (-2)-curve. As above, the existence of further Mori fibers meeting $E_i$ follows. Again, by invariance, the total number of Mori fibers exceeds 20, a contradiction. It follows that either $H$ is cyclic or $ |G.C|\leq 9$. Both imply $|G| \leq 108$. 
\item
If $N=3$, let $S=\{p_1,p_2,p_3\}$ be the points
of intersection of $C$ with the union $\bigcup E_i$ of Mori fibers. The set $S$ is
$H$-invariant. It follows that $H$ is either trivial or isomorphic to $C_2$, $C_3$
or $D_6$ and that $|G| \leq 60$
\item
If $N = 4$, it follows from Lemma \ref{boundforn} that $n \leq 9 $. Now $|H| \leq 12$ implies $|G| \leq 108$. The bound for the
order of $H$ is attained by the tetrahedral group $T_{12}$. If the group $G$
does not contain a tetrahedral group, then $|H| \leq 8$ and $|G| \leq
72$.
\item
If $N=5$, the largest possible group acting on $C$ such
that there is an invariant subset of car\-di\-na\-li\-ty 5 is the dihedral
group $D_{10}$. Since \ref{boundforn} implies $n \leq 6$, we conclude $|G| \leq 60$.
\item 
If $N=6$, then  $n \leq 4 $ and $|H| \leq 24$ implies $|G| \leq
96$. This bound is attained if and only if $H \cong O_{24}$. If there is no
octahedral group in $G$, then $|H| \leq 12$ and $|G|\leq 48$.
\item
If $N\geq 12$, then $n=1$ and $H =G$. The maximal order 60 is attained by
the icosahedral group. 
\item
If $6 < N <12$,
we combine $n \leq 4$ and $|H| \leq 24$ to obtain
$|G| \leq 96$.
If $H$ is not the octahedral group, then  $|H| \leq 16$ and $|G| \leq
64$.
\end{itemize}
The case by case discussion shows that the existence of a rational curve in $B$ implies $| G|\leq 108$ and the proposition follows.
\end{proof}
\begin{remark}
If the group $G$ under consideration does not contain certain subgroups (such as large dihedral groups or $T_{12}$, $O_{24}$ or $I_{60}$) then the condition $|G| >108$ in the proposition above can be improved and non-existence of rational ramification curves also follows for smaller $G$.
\end{remark}
\subsection{Elliptic branch curves}
The aim of this section is to find conditions on the order of $G$ which allow us to exclude elliptic curves in $\mathrm{Fix}_X(\sigma)$. We prove:
\begin{proposition}\label{elliptic branch}
Let $X$ be a K3-surface with an action of a finite group $G \times \langle \sigma
\rangle$ such that $ G < \mathrm{Aut}_\mathrm{symp}(X)$ and $\sigma$ is
an antisymplectic involution with fixed points. If $| G | > 108$, 
then $\mathrm{Fix}_X(\sigma)$ contains neither rational nor elliptic ramification curves.
\end{proposition}
\begin{proof}
By the previous proposition $\mathrm{Fix}_X(\sigma)$ contains no rational curves. It follows from Nikulin's description of $\mathrm{Fix}_X(\sigma)$ (cf. Theorem \ref{FixSigma}) that it is either a single curve of genus $g \geq 1$ or the disjoint union of two elliptic curves.

Let $T \subset B$ be an elliptic branch curve and let $H:= \mathrm{Stab}_G(T)$. If $H \neq G$, then $H$ has index two in $G$. The action of $H$ on $T$ is effective.
The automorphism group $\mathrm{Aut}(T)$ of $T$ is a semidirect product $L \ltimes T$, where $L$ is a
linear cyclic group of order at most 6. We consider the projection $\mathrm{pr}_L: \mathrm{Aut}(T) \to
L$  and let $\lambda \in \mathrm{Pr}_L(H)$ be a generating root of unity. We consider $T$ as a quotient $\mathbb C / \Gamma$ and choose $h \in H$ with $h(z) = \lambda z + \omega$ and $t \in T$ such that $\omega +(1-\lambda)t = 0$. After conjugation with the translation $z \mapsto z+t$ the group $H < \mathrm{Aut}(T)$ inherits the semidirect
product structure of $\mathrm{Aut}(T)$, i.e.,
\[
H = (H \cap L) \ltimes  (H \cap T) .
\]
We refer to this decomposition as the \emph{normal form} of $H$.
By Lemma \ref{conicbundle} a $G$-minimal model of $Y$ is a Del Pezzo surface and therefore does not admit elliptic curves with self-intersection zero. It follows that $T$ meets the union $\bigcup E_i$ of Mori fibers. Let $E$ be a Mori fiber meeting $T$. By Proposition \ref{at most two} $|T \cap  E |  \in \{1,2\}$. The stabilzer of $E$ in $H$ is denoted by $\mathrm{Stab}_H(E)$. Since the total number of Mori fibers is bounded by 9 (cf. Lemma \ref{moribound}), the index of $\mathrm{Stab}_H(E)$ in $H$ is bounded by 9.

If $T \cap E =\{p\}$, then $\mathrm{Stab}_H(E)$ is a cyclic group of order less than or equal to six. It follows that $|G | \leq 6\cdot 9 \cdot 2 = 108$.

If $T \cap E =\{p_1, p_2\}$, then $B \cap E = T \cap E$ and the stabilizer $\mathrm{Stab}_G(E)$ of $E$ in $G$ is contained in $H$. If both points $p_1,p_2$ are fixed by $\mathrm{Stab}_G(E)$, then $|\mathrm{Stab}_G(E)| \leq 6$. If $p_1,p_2$ form a $\mathrm{Stab}_G(E)$-orbit, then in the normal form $|\mathrm{Stab}_G(E) \cap T | =2$. It follows that $\mathrm{Stab}_G(E)$ is either $C_2$ or $D_4= C_2 \times C_2$. The index of $\mathrm{Stab}_G(E)$ in $G$ is bounded by 9 and $|G| \leq 54$.

In summary, the existence of an elliptic curve in $B$ implies $|G| \leq 108$ and the proposition follows. 
\end{proof}
\section{Rough classification}
With the preparations of the previous sections we may now turn to a classification result for K3-surfaces with antisymplectic involution centralized by a large group.
\begin{theorem}\label{roughclassi}
  Let $X$ be a K3-surface with a symplectic action of $G$ centralized by an
antisymplectic involution $\sigma$ such that
$\mathrm{Fix}(\sigma)\neq \emptyset$. If $|G|>96$, then $Y$ is a $G$-minimal
Del Pezzo surface and there are no rational or elliptic curves in $\mathrm{Fix}(\sigma)$. In
particular, $\mathrm{Fix}(\sigma)$ is a single smooth curve $C$ with $ g(C)\geq
3$ and $\pi(C) \sim -2K_Y$, where $K_Y$ denotes the canonical divisor on $Y$.
\end{theorem}
\begin{proof}
The group $G$ is a subgroup of one of the eleven groups on Mukai's list \cite{mukai} (cf. Theorem \ref{mukaithm} and Table \ref{TableMukai}). 
The orders of these Mukai groups are 48, 72, 120, 168, 192, 288, 360, 384, 960. None of these groups can have a subgroup $G$ with $96 < |G| < 120$.
In particular, the order of $G$ is at least 120.

We may therefore apply the results of the previous two sections and conclude that $\pi: X \to Y$ is branched along a single smooth curve $C$ of general type. Its genus $g(C)$ must be $\geq 3$ by Hurwitz' formula.
It remains to show that $Y$ is $G$-minimal.

Assume the contrary and let $E\subset Y$ be a Mori fiber with $E^2 = -1$. As before we let $B \subset Y$ denote the branch locus of $\pi: X \to Y$. By Remark \ref{self-int of Mori-fibers} $E \cap B \neq \emptyset$.
It follows that $|E \cap B| \in \{1,2\}$. Let $\mathrm{Stab}_G(E)$ denote the stabilizer of $E$ in $G$. 

If $\pi^{-1}(E)$ is reducible its two irreducible components meet transversally in one point corresponding to $\{p\} =E \cap B$. The curve $E$ is tangent to $B$ at $p$ and we consider the linearization of the action of $\mathrm{Stab}_G(E)$ at $p$. If the action of $\mathrm{Stab}_G(E)$ on $E$ is not effective, the linearization of the ineffectivity $I < \mathrm{Stab}_G(E)$ yields a trivial action of $I$ on the tangent line of $B$ at $p$. It follows that the action of $I$ is trivial in a neighbourhood of $\pi^{-1}(p) \in R =\pi^{-1}(B)$. This is contrary to the assumption that $G$ acts symplectically on $X$.  
Consequently, the action of  $\mathrm{Stab}_G(E)$ on $E$ is effective and in particular, $\mathrm{Stab}_G(E)$ is a cyclic group. 

If $\pi^{-1}(E)$ is irreducible, then it is a smooth rational curve with an effective action of $\mathrm{Stab}_G(E)$. It follows that $\mathrm{Stab}_G(E)$ is either cyclic or dihedral. The largest dihedral group with a symplectic action on a K3-surface is $D_{12}$ (Proposition 3.10 in \cite{mukai}).

We conclude that the order of  $\mathrm{Stab}_G(E)$ is bounded 12 and the index of $G_E$ in $G$ is $> 9 $. By Lemma \ref{moribound} the total number $m$ of Mori fibers however satifies $m \leq 9$. This contradiction shows that $Y$ is $G$-minimal and, in particular, a Del Pezzo surface. 
\end{proof}
\begin{remark}\label{stab of minus one curve}
Let $X$ be a K3-surface with a symplectic action of $G$ centralized by an
antisymplectic involution $\sigma $ with $\mathrm{Fix}_X(\sigma) \neq \emptyset$ and let $E$ be a (-1)-curve on $Y = X/\sigma$. Then the argument above can be applied to see that the stabilizer of $E$ in $G$ is cyclic or dihedral and therefore has order at most 12.
\end{remark}
In the following chapter, the classification above is applied and extended to the case where $G$ is a maximal group of symplectic transformations on a K3-surface. 
%
%
%
%
\chapter{Mukai groups centralized by antisymplectic involutions}\label{chaptermukai}
In this chapter we consider K3-surfaces with a symplectic action of one of the eleven groups from Mukai's list (Table \ref{TableMukai}) and assume that it is centralized by an antisymplectic involution. We prove the following classification result.
\begin{theorem}\label{thm mukai times invol}
Let $G$ be a Mukai group acting on a K3-surface $X$ by symplectic transformations and $\sigma $ be an antisymplectic involution on $X$ centralizing $G$ with $\mathrm{Fix}_X(\sigma) \neq \emptyset$. Then the pair $(X,G)$ is in Table \ref{Mukai times invol} below. In particular, for groups $G$ numbered 4-8 on Mukai's list, there does not exist a K3-surface with an action of $G \times C_2$ with the properties above. 
\renewcommand{\baselinestretch}{1.4}
\begin{table}[h]
\centering
\begin{tabular}{l|l|l|l}
  & $G$ & $|G|$ & \textbf{K3-surface} $X$ \\ \hline 
1a & $L_2(7)$ & 168 & $\{x_1^3x_2+x_2^3x_3+x_3^3x_1+x_4^4 =0\} \subset \mathbb P_3$\\ \hline
1b & $L_2(7)$ & 168 & Double cover of $\mathbb P_2$ branched along \\
& & & 		$\{x_1^5x_2+x_3^5x_1+x_2^5x_3-5x_1^2 x_2^2 x_3^2 =0\}$\\ \hline
2 & $A_6$ & 360 & Double cover of $\mathbb P_2$ branched along \\
& & & 		$\{10 x_1^3x_2^3+ 9 x_1^5x_3 + 9 x_2^3x_3^3-45 x_1^2 x_2^2 x_3^2-135 x_1 x_2 x_3^4 + 27 			x_3^6  =0\}$\\ \hline
3a & $S_5$ & 120 & $\{\sum_{i=1}^5 x_i = \sum_{i=1}^6 x_1^2 = \sum_{i=1}^5 x_i^3=0\} \subset \mathbb P_5$\\ \hline
3b & $S_5$ & 120 & Double cover of $\mathbb P_2$ branched along \\
& & & 		$\{ F_{S_5} =0\}$\\ \hline
9 & $N_{72}$ & 72 & $\{ x_1^3+ x_2 ^3 + x_3^3 +x_4^3= x_1x_2 + x_3x_4+ x_5^2 = 0 \} \subset \mathbb P_4$ \\ \hline
10 & $M_9$ & 72 & Double cover of $\mathbb P_2$ branched along \\
& & & 		$\{x_1^6+x_2^6 +x_3^6 -10(x_1^3x_2^3 + x_2^3x_3^3 +x_3^3x_1^3) =0\}$\\ \hline
11a & $ T_{48}$ & 48 & Double cover of $\mathbb P_2$ branched along \\
& & & 		$\{x_1x_2(x_1^4-x_2^4)+ x_3^6 =0\}$\\ \hline
11b & $ T_{48}$ & 48 & Double cover of $\{  x_0x_1(x_0^4-x_1^4)+ x_2^3+x_3^2=0 \} \subset \mathbb P(1,1,2,3)$\\
& & 		& branched along $\{x_2=0\} $
\end{tabular}
\caption{K3-surfaces with $G \times C_2$-symmetry}\label{Mukai times invol}
\end{table}
\renewcommand{\baselinestretch}{1.1}\\
The polynomial $F_{S_5}$ in case 3b) is given by 
\begin {align*}
&2(x^4yz+xy^4z+xyz^4)-2(x^4y^2+x^4z^2+x^2y^4+x^2z^4+y^4z^2+y^2z^4)\\
+&2(x^3y^3+x^3z^3+y^3z^3)+x^3y^2z+x^3yz^2+x^2y^3z+x^2yz^3+xy^3z^2+xy^2z^3-6x^2y^2z^2.
\end{align*}
\end{theorem}
\begin{remark}
The examples 1a, 3a, 9, 10, and 11a appaer in Mukai's list,
the remaining cases 1b, 2, 3b, and 11b provide additional examples of K3-surfaces with maximal symplectic symmetry.
\end{remark}
For the proof of this theorem we consider each group separately and apply the following general strategy.

For a K3-surface $X$ with $G \times C_2$-symmetry we consider the quotient $Y= X/C_2$ and a $G$-minimal model $Y_\mathrm{min}$ of the rational surface $Y$. We show that $Y_\mathrm{min}$ is a Del Pezzo surface and investigate which Del Pezzo surfaces admit an action of the group $G$.

It is then essential to study the branch locus $B$ of the covering $X \to Y$. As a first step, we exclude rational and elliptic curves in $B$. In order to exclude rational branch curves, we study their images in $Y_\mathrm{min}$ and their intersection with the union of Mori fibers. 

We then deduce that $B$ consists of a single curve of genus $\geq 2$ with an effective action of the group $G$. The possible genera of $B$ are restricted by the nature of the group $G$ and the Riemann-Hurwitz formula for the quotient of $B$ by an appropriate normal subgroup $N$ of $G$. The equations of $B$ or $X$ given in Table \ref{Mukai times invol} are derived using invariant theory.

Throughout the remainder of this chapter, the Euler characteristic formula\index{Euler characteristic formula!}
\[
 24 = e(X) = 2 e(Y_\mathrm{min}) + 2 m -2n +\underset{\text{branch curve exists}}{\underset{\text{if non-rational}}{\underbrace{(2g-2)}}}
\]
is exploited various times.
Here $m$ denotes the total number of Mori contractions of the reduction $Y \to Y_ \mathrm{min}$, the total number of rational branch curves is denoted by $n$ and $g$ is the genus of a non-rational branch curve.

All classification results are up to \emph{equivariant equivalence}:
\begin{definition}\label{equivariantequivalence}
Let $(X_1, \sigma_1)$ and $(X_2, \sigma_2)$ be K3-surfaces with antisymplectic involution and let $G$ be a finite group acting on $X_1$ and $X_2$ by 
\[
 \alpha_i: G \to \mathrm{Aut}_\mathrm{symp}(X_i),
\]
such that $\alpha_i(g) \circ \sigma_i = \sigma_i \circ \alpha_i(g)$ for $i =1,2$ and all $g \in G$.
Then the surfaces $(X_1, \sigma_1)$ and $(X_2, \sigma_2)$ are considered \emph{equivariantly equivalent}\index{equivariantly equivalent} if there exist a biholomorphic map $\varphi: X_1 \to X_2$ and a group automorphism $\psi \in \mathrm{Aut}(G)$ such that
\[
 \alpha_2(g) \varphi (x) = \varphi ( \alpha_1(\psi(g))x) \quad \text{and} \quad \sigma_2(\varphi (x)) = \varphi ( \sigma_1(x)).
\]
for all $x \in X_1$ and all $g \in G$.

More generally, two surfaces $Y_1$ and $Y_2$, without additional structure such as a symplectic form or an involution, with actions of a finite group $G$
\[
 \alpha_i: G \to \mathrm{Aut}(Y_i)
\]
are considered equivariantly equivalent if there exist a biholomorphic map $\varphi: Y_1 \to Y_2$ and a group automorphism $\psi \in \mathrm{Aut}(G)$ such that
\[
 \alpha_2(g) \varphi (y) = \varphi ( \alpha_1(\psi(g))y)
\]
for all $y \in Y_1$ and all $g \in G$.
\end{definition}
This notion differs from the notion of equivalence in representation theory. Two non-equivalent linear represenations of a group $G$ can induce equivalent actions on the projective plane if they differ by an outer automorphism of the group.
\begin{remark}
If two K3-surfaces $(X_1, \sigma_1)$ and $(X_2, \sigma_2)$ are $G$-equivariantly equivalent, then the quotient surfaces $X_i / \sigma_i$ are equivariantly equivalent with respect to the induced action of $G$. 

Conversely, let $Y$ be a rational surface with two action of a finite group $G$ which are equivalent in the above sense and let $\varphi \in \mathrm{Aut}(Y)$ be the isomorphims identifying these two actions. We consider a smooth $G$-invariant curve $B$ linearly equivalent to $-2K_Y$ and the K3-surfaces $X_B$ and $X_{\varphi(B)}$ obtained as double covers branched along $B$ and $\varphi(B)$ equipped with their respective antisymplectic covering involution. 

Note that $X_B$ and $X_{\varphi(B)}$ are constructed as subsets of the anticanonical line bundle where the involution $\sigma$ is canonically defined. The induced biholomorphic map $\varphi_X: X_B \to X_{\varphi(B)}$ fulfills $\sigma \circ \varphi_X = \varphi_X \circ \sigma$ by construction.

If all elements of the group $G$ can be lifted to symplectic transformations on $X_B$ and $X_{\varphi(B)}$, then the central degree two extensions $E$ of $G$ acting on $X_B$, $X_{\varphi(B)}$, respectively, split as $E = E_\mathrm{symp} \times C_2$ with $E_\mathrm{symp} =G$. In this case the group $G$ acts by symplectic transformations on 
$X_B$ and $X_{\varphi(B)}$ and these are $G$-equivariantly equivalent in strong sense introduced above.
This follows from the assumption that the corresponding $G$-actions on the base $Y$ are equivalent and the fact that for each $g \in G \subset \mathrm{Aut}(Y)$ there is only one choice of symplectic lifting $\tilde g \in \mathrm{Aut}(X_B)$ and $\mathrm{Aut}(X_{\varphi(B)})$.
\end{remark}
In the following sections we will go through the lists of Mukai groups and for each group we prove the classification claimed in Theorem \ref{thm mukai times invol}.
\section{The group $L_2(7)$}\label{mukaiL2(7)}\index{Mukai group! $L_2(7)$}
Let $G \cong L_2(7)$ be the finite simple group of order 168. If $G$ acts on a K3-surface $X$, then the kernels of the homomorphism $G \to \mathrm{Aut}(X)$ and the homomorphism $ G \to \Omega^2(X)$ are trivial and the action is effective and symplectic. Let $\sigma$ be an antisymplectic involution on $X$ centralizing $G$. Since $G$ has an element of order seven which is known to have exactly three fixed points $p_1, p_2,p_3$ and $\sigma$ acts on this set of three points, we know that $\mathrm{Fix}_X(\sigma) \neq \emptyset$. By Theorem \ref{roughclassi}, the K3-surface $X$ is a double cover of a Del Pezzo surface $Y$. Our study of Del Pezzo surfaces with an action of $L_2(7)$ in Example \ref{DelPezzoL2(7)} has revealed that $Y$ is either $\mathbb P_2$ or a Del Pezzo surface of degree 2. In the first case, $\pi: X \to Y$ is branched along a curve of genus 10, in the second case $\pi$ is branched along a curve of genus 3. 
Section \ref{168} in the next chapter is devoted to an inspection of K3-surfaces with an action of $L_2(7) \times C_2$ and a precise classification result in the setup above will be obtained. The pair $(X,G)$ is equivariantly isomorphic to either the surface 1a) or 1b).
\section{The group $A_6$}\label{A6Valentiner}\index{Mukai group! $A_6$}
Let $G \cong A_6$ be the alternating group degree 6. It is a simple group and if it acts on a K3-surface $X$, then this action effective and symplectic. Let $\sigma$ be an antisymplectic involution on $X$ centralizing $G$ and assume that $\mathrm{Fix}_X(\sigma) \neq \emptyset$.  By Theorem \ref{roughclassi}, the K3-surface $X$ is a double cover of a Del Pezzo surface $Y$ with an effective action of $A_6$. 
\begin{lemma}
The Del Pezzo surface $Y$ is isomorphic to $\mathbb P_2$ with a uniquely determined action of $A_6$ given by a nontrivial central extension $V=3.A_6$ of degree three known as \emph{Valentiner's group}\index{Valentiner's group}. 
\end{lemma}
\begin{proof}
We go through the list of Del Pezzo surfaces.
\begin{itemize}
 \item 
If $Y$ has degree one, then $| -K_Y|$ has precisely one base point which would have to be an $A_6$-fixed point. This is contrary to the fact that $A_6$ has no faithful two-dimensional representation. 
\item
We recall that the stabilizer of a (-1)-curve $E$ in $Y$ is either cyclic or dihedral (Remark \ref{stab of minus one curve}). In particular, its order is at most 12 and therefore its index in $A_6$ is at least 30. 
Using Table \ref{minus one curves} we see that $Y$ can not be a Del Pezzo surface of degree 2,3,4,5,6.
\item
Since the blow-up of $\mathbb P_2$ in one point is never $G$-minimal,
it remains to exclude $Y \cong \mathbb P_1 \times \mathbb P_1$. Assume there is an action of $A_6$ on $\mathbb P_1 \times \mathbb P_1$. Since $A_6$ has no subgroups of index two, it follows that $A_6 < \mathrm{PSL}(2, \mathbb C) \times  \mathrm{PSL}(2, \mathbb C)$ and both canonical projections are $A_6$-equivariant. Since $A_6$ has neither an effective action on $\mathbb P_1$ nor nontrivial normal subgroups of ineffectivity, it follows that $A_6$ acts trivially on $Y$. 
\end{itemize}
It follows that $Y \cong \mathbb P_2$. The action of $A_6$ on $\mathbb P_2$ is given by a degree three central extension of $A_6$. Since $A_6$ has no faithful three-dimensional representation, this extension is nontrivial and isomorphic the unique nontrivial degree three extension $V=3.A_6$ known as Valentiner's group. Up to equivariant equivalence, there is a unique action of $A_6$ on $\mathbb P_2$. This follows from the classification of finite subgroup of $\mathrm{SL}_3(\mathbb C)$ (cf. \cite{blichfeldtbook}, \cite{blichfeldt}, and \cite{YauYu}) and can also be derived as follows: An action of $A_6$ on $\mathbb P_2$ is given by a threedimensional projective representation. We wish to show that any two actions induced by $\rho_1, \rho_2$ are equivalent. We restrict the projective representations $\rho_1$ and $\rho_2$ to the subgroup $A_5$. The restricted representations are linear and after a change of coordinates $\rho_1(A_5) = \rho_2(A_5) \subset \mathrm{SL}_3(\mathbb C)$. We fix a subgroup $A_4$ in $A_5$ and consider its normalizer $N$ in $A_6$. The groups $N$ and $A_4$ generate the full group $A_6$ and it suffices to prove that $\rho_1(N) = \rho_2(N)$. This is shown by considering an explicit three-dimensional representation of $A_4 < A_5$ and the normalizer $\mathcal N$ of $A_4$ inside $\mathrm{PSL}_3(\mathbb C)$. The group $A_4$ has index two in $\mathcal N$ and therefore $\mathcal N = \rho_1(N)= \rho_2(N)$.. 
\end{proof}
The covering $X \to Y$ is branched along an invariant curve $C$ of degree six. This curve is defined by an invariant polynomial $F_{A_6}$ of degree six, which is unique by Molien's formula. Its explicit equation is derived in \cite{crass}. In appropriately chosen coordinates,
\[
F_{A_6}(x_1,x_2,x_3) = 10 x_1^3x_2^3+ 9 x_1^5x_3 + 9 x_2^3x_3^3-45 x_1^2 x_2^2 x_3^2-135 x_1 x_2 x_3^4 + 27 x_3^6.
\]
If a K3-surface with $A_6 \times C_2$-symmetry exists, then it must be the double cover of $\mathbb P_2$ branched along $\{F_{A_6}=0\}$.

The action of $A_6$ on $\mathbb P_2$ induces an action of a central degree two extension of $E$ on the double cover branched along $\{F_{A_6}=0\}$,
\[ 
\{\mathrm{id}\} \to C_2 \to E \to A_6 \to \{\mathrm{id}\}.
 \]

 Let $E_\mathrm{symp} \neq E$ be the normal subgroup of symplectic automorphisms in $E$. Since $A_6$ is simple, it follows that $E_\mathrm{symp}$ is mapped surjectively to $A_6$ and $E_\mathrm{symp} \cong A_6$. In particular, the group $E$ splits as $E_\mathrm{symp} \times C_2$ where $C_2$ is generated by the antisymplectic covering involution. This proves the existence of a unique K3-surface with $A_6 \times C_2$-symmetry. We refer to this K3-surface as the \emph{Valentiner surface}\index{Valentiner surface}.
\section{The group $S_5$}\index{Mukai group! $S_5$}
In this section we study K3-surfaces with an symplectic action of the symmetric group $S_5$ centralized by an antisymplectic involution. 

Let $X$ be a K3-surface with a symplectic action of $G = S_5$ and let $\sigma$ denote an antisymplectic involution centralizing $G$. We assume that $\mathrm{Fix}_X(\sigma) \neq \emptyset$. 
We may apply Theorem \ref{roughclassi} which yields that $X/ \sigma =Y$ is a $G$-minimal Del Pezzo surface and $\pi: X \to Y$ is branched along a smooth connected curve $B$ of genus 
\[
g(B) = 13-e(Y).
\]
We will see in the following that only very few Del Pezzo surfaces admit an effective action of $S_5$ or a smooth $S_5$-invariant curve of appropriate genus.
\begin{lemma}
The degree $d(Y)$ of the Del Pezzo surface $Y$ is either three or five. 
\end{lemma}
\begin{proof}
We prove the statement by excluding Del Pezzo surfaces of degree $\neq 3,5$.
\begin{itemize}
\item 
Assume $Y \cong \mathbb P_2$. Then $G =S_5$ is acting effectively on $\mathbb P_2$, i.e., $S_5 \hookrightarrow \mathrm{PSL}_3(\mathbb C)$. Let $\tilde G$ denote the preimage of $G$ in $\mathrm{SL}_3(\mathbb C)$. Since $A_5$ has no nontrivial central extension of degree three, it follows that the preimage of $A_5 < S_5$ in $\tilde G$ splits as $\tilde A_5 = A_5 \times C_3$. It has index two in $\tilde G$ and therefore is a normal subgroup of $\tilde G$. Let $g \in S_5$ be any transposition and pick $\tilde g$ in its preimage with $\tilde g^2 = \mathrm{id}$. Now $\tilde g$ and $ A_5$ generate a copy of $S_5$ in $\mathrm{SL}_3(\mathbb C)$. The action of $S_5$ is given by a three-dimensional representation.
The irreducible representations of $S_5$ have dimensions $1,4,5$ or $6$ and it follows that there is no faithful three-dimensional represenation of $S_5$ and therefore no effective $S_5$-action on $\mathbb P_2$. 
\item
Assume that $Y$ is isomorphic to $\mathbb P_1 \times \mathbb P_1$. We investigate the action of $S_5 = A_5 \rtimes C_2$ and note that $A_5$ is a simple group. The automorphism group $\mathrm{Aut}(Y)$ is given by 
\[
 (\mathrm{PSL}_2( \mathbb C) \times  \mathrm{PSL}_2(\mathbb C)) \rtimes C_2.
\]
It follows that $A_5 < \mathrm{PSL}_2( \mathbb C) \times  \mathrm{PSL}_2( \mathbb C)$, and the action of $A_5$ respects the product structure, i.e, the canonical projections onto the factors are $A_5$-equivariant.
If $A_5$ acts trivially on one of the factors, then the generator $\tau$ of the outer $C_2$ stabilizes this factor because $A_5$ must act nontrivially on the second factor. It follows that $S_5$ stablizes the second factor which is impossible since there is no effective action of $S_5$ on $\mathbb P_1$. It follows that $A_5$ acts effectively on both factors and $\tau$ exchanges them.
We consider an element $\lambda$ of order five in $A_5$ and chose coordinates on $\mathbb P_1 \times \mathbb P_1$ such that $\lambda$ acts by
\[
 ([z_1:z_2],[w_1:w_2]) \mapsto ([\xi z_1:z_2],[\xi ^a w_1:w_2])
\]
for some $a \in \{1,2,3,4\}$ and $\xi^5 =1$. The automorphism $\lambda$ has four fixed points 
\begin{align*}
p_1 = ([1:0],[1:0]),\,\,
p_2 = ([1:0],[0:1]),\, \,
p_3 = ([0:1],[1:0]),\,\,
p_4 = ([0:1],[0:1]).\,\,
\end{align*}
Since it lifts to a symplectic automorphism on the K3-surface $X$ with four fixed points, all fixed points must lie on the branch curve.
The branch curve $B \subset Y$ is a smooth invariant curve linearly equivalent to $-2K_Y$ and is therefore given by an $S_5$-semi-invariant polynomial $f$ of bidegree $(4,4)$. Since $f$ must be invariant with respect to the subgroup $A_5$, it is a linear combination of $\lambda$-invariant monomials of bidegree $(4,4)$.
For each choice of $a$ one lists all $\lambda$-invariant monomials of bidegree $(4,4)$. For $a=1$ these are
\[
z_1 z_2^3 w_1^4,\,\, z_1^2 z_2^2 w_1^3 w_2 ,\,\, z_1^3 z_2 w_1^2 w_2^2,\,\,  z_1^4 w_1 w_2^3 , \,\,z_2^4 w_2^4.
\]
Since $f$ must vanish at $p_1 \dots p_4$, one sees that $f$ may not contain $z_2^4 w_2^4$. The remaining monomials have a common component $z_1 w_1$ such that $f$ factorizes and $C$ must be reducible, a contradiction.
The same argument can be carried out for each choice of $a$. It follows that the action of $S_5$ on  $\mathbb P_1 \times \mathbb P_1$ does not admit irreducble curves of bidegree $(4,4)$. This eliminates the case $Y \cong \mathbb P_1 \times \mathbb P_1$.
\item
Again using the fact that the largest subgroup of $S_5$ which can stabilize a (-1)-curve in $Y$ is the group $D_{12}$ of index 10,
it follows that the number of (-1)-curves in a $G$-orbit is at least 10. 
A Del Pezzo surface of degree six has six (-1)-curves and therefore $d(Y) \neq 6$.
A Del Pezzo surface of degree four contains sixteen (-1)-curves. Since 16 does not divide the the order of $S_5$, the set of these curves is not a single $G$-orbit. As it cannot be the union of $G$-orbits either, we conclude $d(Y) \neq 4$.
\item
If $d(Y) =2$, then the anticanonical map defines an $\mathrm{Aut}(Y)$-equivariant double cover of $\mathbb P_2$. The induced action of $S_5$ on $\mathbb P_2$ would have to be effective and therefore we obtain a contradiction as in the case $Y \cong \mathbb P_2$.
\item
If $d(Y)=1$ then the anticanonical system $|-K_Y|$ is known to have precisely one base point which has to be fixed point of the action of $S_5$. Since $S_5$ has no faithful two-dimensional representation, this is a contradiction. 
\end{itemize}
Since we have considered all possible $G$-minimal Del Pezzo surfaces the proof of the lemma is completed. 
\end{proof}
\subsection{Double covers of Del Pezzo surfaces of degree three}
The following example of a K3-surface $X$ with an action of $S_5 \times C_2$ such that $X/\sigma$ is a Del Pezzo surface of degree three can be found in Mukai's list \cite{mukai} (cf. also Table \ref{TableMukai}). 
\begin{example}\label{MukaiS5}
Let $X$ be the K3-surface in $\mathbb P_5$ given by 
\[
\sum_{i=1}^5 x_i = \sum_{i=1}^6 x_1^2 = \sum_{i=1}^5 x_i^3=0
 \]
 and let $S_5$ act on $\mathbb P_5$ by permuting the first five variables and by the character $\mathrm{sgn}$ on the last variable. This induces an action on $X$. 

The commutator subgroup $ S_5' = A_5 < S_5$ acts by symplectic transformations. In order to show that the full group acts symplectically, consider the transposition $\tau = (12) \in S_5$ acting on $\mathbb P_5$ by $[x_1:x_2:x_3:x_4:x_5:x_6] \mapsto [x_2:x_1:x_3:x_4:x_5:-x_6]$. One checks that the induced involution on $X$ has isolated fixed points and is therefore symplectic. It follows that $S_5 < \mathrm{Aut}_\mathrm{symp}(X)$. 

Let $\sigma : \mathbb P_5 \to \mathbb P_5$ be the involution $[x_1:x_2:x_3:x_4:x_5:x_6] \mapsto [x_1:x_2:x_3:x_4:x_5:-x_6]$. This defines an involution on $X$ with a positive-dimensional set of fixed point $\{x_6=0 \} \cap X$. Therefore $\sigma$ is an antisymplectic involution on $X$ which centralizes the action of $S_5$. 

The quotient $Y$ of $X$ by $\sigma$ is given by restricting then rational map $[x_1:x_2:x_3:x_4:x_5:x_6] \mapsto [x_1:x_2:x_3:x_4:x_5]$ to $X$. The surface $Y$ is given by 
\[
\{ \sum_{i=1}^5 y_i= \sum_{i=1}^5 y_i^3 =0\} \subset \mathbb P_4.
\]
and is isomorphic to the Clebsch diagonal surface  
$\{z_1^2 z_2 + z_1 z_3^2 + z_3 z_4^2 + z_4 z_2^2 = 0\} \subset \mathbb P_3$ (cf. Theorem 10.3.10 in \cite{dolgachev}), a Del Pezzo surface of degree three. The branch set $B$ is given by $\{ \sum_{i=1}^5 y_1^2=0\} \cap Y \subset \mathbb P_4$.
\end{example}
By the following proposition, the example above is the unique K3-surface with $S_5 \times C_2$-symmetry such that $X/\sigma$ is a Del Pezzo surface of degree three. 
\begin{proposition}\label{S5 on degree three}
Let $X$ be a K3-surface with a symplectic action of the group $S_5$ centralized by an antisymplectic involution $\sigma$. If $Y=X/\sigma$ is a Del Pezzo surface of degree three, then $X$ is equivariantly isomorphic to Mukai's $S_5$-example $\{\sum_{i=1}^5 x_i = \sum_{i=1}^6 x_1^2 = \sum_{i=1}^5 x_i^3=0\} \subset \mathbb P_5$. 
\end{proposition}
\begin{proof}
We consider the $\mathrm{Aut}(Y)$-equivariant embedding of the Del Pezzo surface $Y$ into $\mathbb P_3$ given by the anticanonical map. Any automorphism of $Y$ induced by an automorphism of the ambient projective space. 

It follows from the representation and invariant theory of the group $S_5$ that a Del Pezzo surface of degree three with an effective action of the group $S_5$ is equivariantly isomorphic the Clebsch cubic $\{z_1^2 z_2 + z_1 z_3^2 + z_3 z_4^2 + z_4 z_2^2 = 0\} \subset \mathbb P_3$ (cf. Theorems 10.3.9 and 10.3.10, Table 10.3 in \cite{dolgachev}).

The ramification curve $B\subset Y$ is linearly equivalent to $-2K_Y$. We show that $B$ is given by intersecting $Y$ with a quadric in $\mathbb P_3$. 

Applying the formula 
\[
h^0(Y, \mathcal{O}(-rK_Y))= 1+ \frac{1}{2}r(r+1)d(Y)
\]
(cf. e.g. Lemma 8.3.1 in \cite{dolgachev})
to $d=d(Y)=3$ and $r=2$ we obtain $h^0(Y,\mathcal{O}( -2K_Y))= 10$. This is also the dimension of the space of sections of $\mathcal O_{\mathbb P_3} (2)$ in $\mathbb P_3$ (homogeneous polynomials of degree two in four variables). It follows that the restriction map
\[
H^0(\mathbb P_3, \mathcal O (2)) \to H^0(Y, \mathcal O(-2K_Y))
\]
is surjective and $B = Y \cap Q$ for some quadric $Q = \{f=0\}$ in $\mathbb P_3$.

Since $B$ is an $S_5$-invariant curve in $Y$, it follows that for each $g \in S_5$ the intersection of $gQ = \{ f \circ g^{-1} =0\}$ with Y coincides with $B$. It follows that $f|_Y$ is a multiple of $(f\circ g^{-1})|_Y$, i.e., there exists a constant $c \in \mathbb C^*$ such that $(f \circ g^{-1}) - cf$ vanishes identically on $Y$. Since $Y$ is irreducible, this implies $f \circ g^{-1} = cf$. It follows that the polynomial $f$ is an $S_5$- semi-invariant and therefore invariant with respect to the commutator subgroup $A_5$.

We have previously noted that after a suitable linear change of coordinates the surface $Y$ is given by $\{ \sum_{i=1}^5 y_i= \sum_{i=1}^5 y_i^3 =0\} \subset \mathbb P_4$ where $S_5$ acts  by permutation. The action of any transposition on an $S_5$-semi-invariant polynomial is given by multiplication by $\pm 1$. It follows that in the coordinates $[y_1:\dots:y_5]$ the semi-invariant polynomial $f$ is given by
\[
a \sum_{i=1}^5 y_i^2 + b(\sum_{i=1}^5 y_i)^2 =0
\]
for some $a,b \in \mathbb C$. Using the fact $Y \subset \{\sum_{i=1}^5 y_i=0\}$ it follows that $B$ is given by intersecting $Y$ with $ \{\sum_{i=1}^5 y_i^2 =0\}$ and $X$ is Mukai's $S_5$-example discussed in Example \ref{MukaiS5}.
\end{proof}
\subsection{Double covers of Del Pezzo surfaces of degree five}
A second class of candidates of K3-surfaces with $S_5 \times C_2$-symmetry is given by double covers of Del Pezzo surfaces of degree five. 

Any two Del Pezzo surfaces of degree five are isomorphic and the automorphisms group of a Del Pezzo surface $Y$ of degree five is $S_5$. The ten (-1)-curves on $Y$ form a graph known as the \emph{Petersen graph}. The Petersen graph has $S_5$-symmetry and every symmetry of the abstract graph is induced by a unique automorphism of the surface $Y$. 

The following proposition classifies K3-surfaces with $S_5 \times C_2$-symmetry which are double covers of Del Pezzo surfaces of degree five.
\begin{proposition}
Let $X$ be a K3-surface with a symplectic action of the group $S_5$ centralized by an antisymplectic involution $\sigma$. If $Y=X/\sigma$ is a Del Pezzo surface of degree five, then $X$ is equivariantly isomorphic to the minimal desingularization of the double cover of $\mathbb P_2$ branched along the sextic
\begin{align*}
\{&2(x^4yz+xy^4z+xyz^4)
-2(x^4y^2+x^4z^2+x^2y^4+x^2z^4+y^4z^2+y^2z^4)
+2(x^3y^3+x^3z^3+y^3z^3)\\
&+x^3y^2z+x^3yz^2+x^2y^3z+x^2yz^3+xy^3z^2+xy^2z^3
-6x^2y^2z^2
=0\}
\end{align*}
\end{proposition}
\begin{proof}
Let $B \subset Y$ denote the branch locus of the covering $X \to Y$. The curve $B$ is smooth, connected, invariant with respect to the full automorphism group of $Y$ and linearly equivalent to $-2K_Y$. 

The Del Pezzo surface $Y$ is the blow-up of $\mathbb P_2$ in four points $p_1,p_2,p_3,p_4$ in general position. We may choose coordinates $[x:y:z]$ on $\mathbb P_2$ such that 
\begin{align*}
 p_1=[1:0:0],\quad
 p_2=[0:1:0],\quad
 p_3=[0:0:1],\quad
 p_4=[1:1:1].
\end{align*}
Let $m: Y \to \mathbb P_2$ be the blow-down map and let $E_i = m^{-1}(p_i)$.
Consider the $S_4$-action on $\mathbb P_2$ permuting the points $\{p_i\}$. The isotropy at the point $p_1$ is isomorphic to $S_3$ and induces an effective $S_3$-action on $E_1$. 

Let $E$ be any (-1)-curve on $Y$. By adjunction $E\cdot B=2$. Since $Y$ contains precisely ten (-1)-curves forming an $S_5$-orbit, the group $H = \mathrm{Stab}_{S_5}(E)$ has order 12 and all stabilizer groups of (-1)-curves in $Y$ are conjugate.
It follows that the group $H$ contains $S_3$, which is acting effectively on $E$, and therefore $H$ is isomorphic to the dihedral group of order 12. The points of intersection $B\cap E$  form an $H$-invariant subset of $E$. Since $H$ has no fixed points in $E$ and precisely one orbit $H.p = \{p,q\}$ consisting of two elements, it follows that $B$ meets $E$ transversally in $p$ and $q$.

In particular, each curve $E_i$
meets $B$ in two points and the image curve $C = m(B)$ has nodes at the four points $p_i$. 
By Lemma \ref{selfintblowdown}, the self-intersection number of $C$ is $20 + 4\cdot 4= 36$, so $C$ is a sextic curve. It is invariant with respect to the action of $S_4$ given by permutation on $p_1, \dots p_4$. For simplicity, we first only consider the action of $S_3$ permuting $p_1,p_2,p_3$ and conclude that $C$ is given by $\{f=\sum a_i f_i =0 \}$ as a linear combination of the following degree six polynomials
\begin{align*}
f_1&=x^6+y^6 +z^6\\
f_2&=x^5y + x^5z+ xy^5 +xz^5+y^5z+yz^6\\
f_3&=x^4yz+xy^4z+xyz^4\\
f_4&=x^4y^2+x^4z^2+x^2y^4+x^2z^4+y^4z^2+y^2z^4\\
f_5&=x^3y^3+x^3z^3+y^3z^3\\
f_6&=x^3y^2z+x^3yz^2+x^2y^3z+x^2yz^3+xy^3z^2+xy^2z^3\\
f_7&=x^2y^2z^2
\end{align*}
The fact that $C$ passes through $p_i$ and is singular at $p_i$ yields $a_1=a_2=0$ and
\[
3a_3+6a_4+3a_5+6a_6+a_7=0.
\]
The two tangent lines of $C$ at the node $p_i$ correspond to the unique $\mathrm{Stab}(E_i)$-orbit of length two in $E_i$. We consider the point $p_3$ and the subgroup $S_3 < S_4$ stabilizing $p_3$. The action of $S_3$ on $E_3$ is given by the linearized $S_3$-action on the set of lines through $p_3$. One checks that in local affine coordinates $(x,y)$ the unique orbit of length two corresponds to the line pair $x^2 -xy +y^2 =0$.
Dehomogenizing $f$ at $p_3$, i.e., setting $z=1$, we obtain the local equation $f_\mathrm{dehom}$ of $C$ at $p_3$. 
The polynomial $f_\mathrm{dehom}$ modulo terms of order three or higher must be a multiple of $x^2 -xy +y^2$.
Therefore $a_3 = -a_4$.

Next we consider the intersection of $C$ with the line $L_{34} = \{x=y\}$ joining $p_3$ and $p_4$. We know that $f|_{L_{34}}$ vanishes of order two at $p_3$ and $p_4$ and at one or two further points on $L_{34}$. 

Let $\widetilde L_{34}$ denote the proper transform of $L_{34}$ inside the Del Pezzo surface $Y$. The curve $\widetilde L_{34}$ is a (-1)-curve, hence its stabilizer $\mathrm{Stab}_G(\widetilde L_{34})$ is isomorphic to $D_{12} = S_3 \times C_2$. The factor $C_2$ acts trivially on $\widetilde L_{34}$. Since the intersection of $\widetilde L_{34}$ with $B$ is $\mathrm{Stab}_G(\widetilde L_{34})$ invariant, it follows that $\widetilde L_{34} \cap B$ is the unique $S_3$-orbit a length two in $\widetilde L_{34}$.  

We wish to transfer our determination of the unique $S_3$-orbit of length two in $E_3$ above to the curve $\widetilde L_{34}$ using an automorphism of $Y$ mapping $E_3$ to $\widetilde L_{34}$. Consider the automorphism $\varphi$ of $Y$ induced by the birational map of $\mathbb P_2$ given by
\[
[x:y:z] \mapsto [x(z-y):z(x-y):xz] 
\]
(cf. Theorem 10.2.2 in \cite{dolgachev}) and let $\psi$ be the automorphism of $Y$ induced by the permutation of the points $p_2$ and $p_3$ in $\mathbb P_2$. Then $\psi \circ \varphi$ is an automorphism of $Y$ mapping $E_3$ to $\widetilde L_{34}$. If $[X:Y]$ denote homogeneous coordinates on $E_3$ induced by the affine coordinates $(x,y)$ in a neighbourhood of $p_3$, then a point $[X:Y] \in E_3$ is mapped to the point corresponding to $[X:X: X-Y] \in L_{34} \subset \mathbb P_2$. It was derived above that the unique $S_3$-orbit of length two in $E_3$ is given by $X^2-XY +Y^2$ and it follows that the unique $S_3$-orbit of length two in $\widetilde L_{34}$ corresponds to the points $[x:x:z] \in \mathbb P_2$ fulfilling $x^2 -xz +z^2 =0$. 

Therefore, $f|_{L_{34}}$ is a multiple of polynomial given by $x^2(x-z)^2(x^2 -xz +z^2)$.
Comparing coefficients with $f(x:x:z)$ yields
\begin{align*}
2a_3+2a_6 &= 2a_5 +2a_6\\
2a_4 + a_5&= 2a_4+a_3\\
8a_4+4a_5 &= 2a_4+2a_6 +a_7\\
-6a_4 -3a_5 &= 2a_5 +2a_6.
\end{align*}
We conclude $a_3=a_5=2= -a_4$, $a_6=1$, and $a_7=-6$. So if $X$ as in the lemma exists, it is the double cover of $Y$ branched along the proper transform of $\{f=0\}$ in $Y$, where 
\begin{align*}
f(x,y,z)=&2(x^4yz+xy^4z+xyz^4)\\
&-2(x^4y^2+x^4z^2+x^2y^4+x^2z^4+y^4z^2+y^2z^4)\\
&+2(x^3y^3+x^3z^3+y^3z^3)\\
&+x^3y^2z+x^3yz^2+x^2y^3z+x^2yz^3+xy^3z^2+xy^2z^3\\
&-6x^2y^2z^2.
\end{align*}
In order to prove existence, let $X$ be the minimal desingularisation of the double cover of $\mathbb P_2$ branched along $\{f=0\}$. Then $X$ is the double cover of the Del Pezzo surface $Y$ of degree five branched along the proper transform $D$ of $\{f=0\}$ in $Y$.
Since all automorphisms of $Y$ are induced by explicit biholomorphic or birational transformation of $\mathbb P_2$ 
one can check by direct computations that $D$ is in fact invariant with respect to the action of $\mathrm{Aut}(Y) = S_5$. The covering involution $\sigma$ is antisymplectic. 

On $X$ there is an action of a central extension $E$ of $S_5$ by $C_2$. Let $E_\mathrm{symp}$ be the subgroup of symplectic automorphisms in $E$. Since $E$ contains the antisymplectic covering involution $E_\mathrm{symp} \neq E$. The image $N$ of $E_\mathrm{symp}$ in $S_5$ is normal and therefore either $N \cong S_5$ or $N \cong A_5$.

If $N \cong A_5$ and $| E_\mathrm{symp}| =60$, then  $E_\mathrm{symp} \cong A_5$. Lifting any transposition from $S_5$ to an element $g$ of order two in $E$, the group generated by $g$ and $E_\mathrm{symp}$ inside $E$ is isomorphic to $S_5$. It follows that $E$ splits as $S_5 \times C_2$ and $E / E_\mathrm{symp} \cong C_2 \times C_2$. This is a contradiction.  

If $N \cong A_5$ and $| E_\mathrm{symp}| =120$, then $ E = E_\mathrm{symp} \times C_2$, where the outer $C_2$ is generated by the antisymplectic covering involution $\sigma$, and $E / C_2 = S_5$ implies that $E_\mathrm{symp} \cong S_5$. This is contradictory to the assumption $N \cong A_5$.

In the last remaining case $N \cong S_5$. Since $E_\mathrm{symp} \neq E$, also $E_\mathrm{symp} \cong S_5$ and $E$ splits as $E_\mathrm{symp} \times C_2$. 
It follows that the action of $S_5$ on $Y$ induces an symplectic action of $S_5$ on the double cover $X$ centralized by the antisymplectic covering involution. This completes the proof of the proposition.
\end{proof}
\subsection{Conclusion}
We summarize our results of the previous subsections in the following theorem.
\begin{theorem}
Let $X$ be a K3-surface with a symplectic action of the group $S_5$ centralized by an antisymplectic involution $\sigma$ with $\mathrm{Fix}_X(\sigma) \neq \emptyset$. Then $X$ is equivariantly isomorphic to either Mukai's $S_5$-example or the minimal desingularization of the double cover of $\mathbb P_2$ branched along the sextic
\begin{align*}
\{&F_{S_5}(x_1,x_2,x_3)=\\
&2(x^4yz+xy^4z+xyz^4)-2(x^4y^2+x^4z^2+x^2y^4+x^2z^4+y^4z^2+y^2z^4)+2(x^3y^3+x^3z^3+y^3z^3)\\
&+x^3y^2z+x^3yz^2+x^2y^3z+x^2yz^3+xy^3z^2+xy^2z^3-6x^2y^2z^2=0\}.
\end{align*}
\end{theorem}
\section{The group $M_{20} = C_2^4 \rtimes A_5$}\index{Mukai group! $M_{20}$}
\begin{proposition}
There does not exist a K3-surface with a symplectic action of $M_{20}$ centralized by an antisymplectic involution $\sigma$ with $\mathrm{Fix}_X(\sigma) \neq \emptyset$.
\end{proposition}
\begin{proof}
Assume that a K3-surface $X$ with these properties exists. Applying Theorem \ref{roughclassi} we see that $X \to Y$ is branched along a single $M_{20}$-invariant smooth curve $C$ on the Del Pezzo surface $Y$. The curve $C$ is neither rational nor elliptic. By Hurtwitz' formula, 
\[
| \mathrm{Aut} (C) | \leq 84(g(C)-1),
\]
the genus of $C$ must be at least twelve. Since $C$ is linearly equivalent to $-2K_Y$, the adjunction formula 
\[
2g(C)-2 = (K_Y,C) + C^2 = 2K_Y^2  
\]
implies $\mathrm{deg}(Y) = K_Y^2 \geq 11$. This is a contradiction since the degree of a Del Pezzo surface is at most nine.
\end{proof}
\section{The group $F_{384} = C_2^4 \rtimes S_4$}\index{Mukai group! $F_{384}$}
Before we prove non-existence of K3-surfaces with $F_{384} \times C_2$-symmetry, we note the following useful fact about $S_4$-actions on Riemann surfaces.
\begin{lemma}\label{S4 not on g=1,2}
The group $S_4$ does not admit an effective action on a Riemann surface of genus one or two.
\end{lemma}
\begin{proof}
The automorphism group of a Riemann surface $T$ of genus one is of the form $\mathrm{Aut}(T)= L \ltimes T$ for $L \in \{C_2, C_4, C_6\}$. We have seen before (cf. Proof of Proposition \ref{elliptic branch}) that any subgroup $H$ of $\mathrm{Aut}(T)$ can be put into the form $H = (H \cap L) \ltimes  (H \cap T)$. The nontrivial normal subgroups of $S_4$ are $A_4$ and $C_2 \times C_2$. Since $A_4$ is not Abelian and the quotient of $S_4$ by $S_4 \cap T = C_2 \times C_2$ is not cyclic, we
conclude that $S_4$ is not a subgroup of $ \mathrm{Aut}(T)$.

Assume that $S_4$ acts effectively on a Riemann surface $H$ of genus two. Note that $H$ is hyperelliptic and the quotient of $H$ by the hyperelliptic involution is branched at six points. Since $S_4$ has no normal subgroup of order two, the induced action of $S_4$ on the quotient $\mathbb P_1$ is effective and therefore has precisely one orbit consisting of six points. The isotropy subgroup at these points is isomorphic to $C_4$. The isotropy group at the corresponding points in $H$ must be isomorphic to $C_4 \times C_2$. Since this group is not cyclic, it cannot act effectively with fixed points on a Riemann surface and we obtain a contradiction.
\end{proof}
\begin{proposition}
There does not exists a K3-surface with a symplectic action of $F_{384}$ centralized by an antisymplectic involution $\sigma$ with $\mathrm{Fix}_X(\sigma) \neq \emptyset$.
\end{proposition}
\begin{proof}
As above, assume that a K3-surface $X$ with these properties exists and apply Theorem \ref{roughclassi} to see that $X \to Y$ is branched along a single $F_{384}$-invariant smooth curve $C$ on the Del Pezzo surface $Y$. It follows from Hurwitz' formula that the genus of $C$ is at least 6.

We use the realization of $F_{384}$ as a semi-direct product $C_4^2 \rtimes S_4$ (cf. \cite{mukai}) and consider the quotient $Q$ of the branch curve $C$ by the normal subgroup $N = C_4^2$. On $Q$ there is the induced action of $S_4$. It follows from the lemma above that $Q$ is either rational or $g(Q) >2$. In the second case, if we apply the Riemann-Hurwitz formula to the covering $ C \to Q$, then 
\[
e(C) = 16 e(Q) - \text{branch point contributions} \leq -64
\]
and $g(C) \geq 33$. This contradicts the adjunction formula on the Del Pezzo surface $Y$ and implies that $Q$ is a rational curve.

It follows from adjunction that $K_Y^2 = g(C)-1$. Therefore, the degree of the Del Pezzo surface $Y$ is at least five.
We consider the action of $F_{384}$ on the configuration of (-1)-curves on $Y$ and recall that the order of a stabilizer of a (-1)-curve in $Y$ is at most twelve (cf. Remark \ref{stab of minus one curve}) and therefore has index greater than or equal to $32$ in $G$. It follows that $Y$ is either $\mathbb P_1 \times \mathbb P_1$ or $\mathbb P_2$. In the first case, the canonical projections of $\mathbb P_1 \times \mathbb P_1$ are equivariant with respect to a subgroup of index two in $F_{384}$ and thereby contradict Lemma \ref{conicbundle}. Consequently, $Y \cong \mathbb P_2$. In particular, $g(C) =10$ and $e(C) = -18$. It follows that the branch point contribution of the covering $C \to Q$ must be 50. Since isotropy groups must be cyclic, the only possible isotropy subgroups of $N = C_4^2$ at a point in $C$ are $C_2$ and $C_4$ and have index four or eight. The full branch point contribution must therefore be a multiple of four. This contradiction yields the non-existence claimed. 
\end{proof}
\section{The group $A_{4,4} = C_2^4 \rtimes A_{3,3}$}\index{Mukai group! $A_{4,4}$}
By $S_{p,q}$ for $p+q =n$ we denote a subgroup $S_p \times S_q$ of $S_n$ preserving a partition of the set $\{1,\dots, n\}$ into two subsets of cardinality $p$ and $q$. The intersection of $A_n$ with $S_{p.q}$ is denoted by $A_{p,q}$.
\begin{proposition}
There does not exists a K3-surface with a symplectic action of $A_{4,4}$ centralized by an antisymplectic involution $\sigma$ with $\mathrm{Fix}_X(\sigma) \neq \emptyset$.
\end{proposition}
\begin{proof}
We again assume that a K3-surface with these properties exists. Applying Theorem \ref{roughclassi} we see that $X \to Y$ is branched along a single $A_{4,4}$-invariant smooth curve $C$ on the Del Pezzo surface $Y$. The group $A_{4,4}$ is a semi-direct product $C_2^4 \rtimes A_{3,3}$ (see e.g. \cite{mukai}). We consider the quotient $Q$ of $C$ by the normal subgroup $N \cong C_2^4$. On $Q$ there is an action of $A_{3,3}$.
Since $A_{3,3}$ contains the subgroup $C_3 \times C_3$, which does not act on a rational curve, it follows that $Q$ not rational. We apply the Riemann-Hurwitz formula to the covering $C \to Q$.

If $Q$ is elliptic, then $2g(C) -2$ equals the branch point contribution of the covering $C \to Q$. As above, isotropy groups must be cyclic and the maximal possible isotropy group of the $C_2^4$-action on $C$ is $C_2$ and has index eight in $C_2^4$. Consequently, the branch point contribution at each branch point is eight.
Recall that any group $H$ acting on the torus $Q$ can be put into the form $H = (H \cap L) \ltimes  (H \cap Q)$ for $L \in \{C_2, C_4, C_6\}$. Since $Q$ acts freely, 
the action of $C_3 \times C_3 < A_{3,3}$ on the elliptic curve $Q$ has orbits of length greater than or equal to three. Therefore, the total branch point contribution must be greater than or equal to $24$. In particular, $g(C) = \mathrm{deg}(Y) +1 \geq 13$ contrary to $\mathrm{deg}(Y) \leq 9$.

If $g(Q) \geq 2$, then $g(C) \geq 17$ which is also contrary to $\mathrm{deg}(Y) \leq 9$
\end{proof}
\section{The groups $T_{192} = (Q_8 * Q_8) \rtimes S_3$ and $H_{192} = C_2^4  \rtimes D_{12}$}\index{Mukai group! $T_{192}$}
By $Q_8$\nomenclature{$Q_8$}{the quaternion group} we denote the quaternion group $\{+1,-1, +I,-I, +J,-J,+K,-K\}$ where $I^2= J^2= K^2 = IJK = -1$. 
The central product $ Q_8 * Q_8 $ is defined as the quotient of $Q_8  \times Q_8$ by the central involution $(-1, -1)$, i.e., $Q_8 * Q_8  = (Q_8  \times Q_8) / (-1,-1)$.

Note that both groups $T_{192}$ and $H_{192}$ are semi-direct products $C_2^3 \rtimes S_4$  (cf. \cite {mukai}).
\begin{proposition}
For $G =T_{192}$ or $G = H_{192}$ there does not exists a K3-surface with a symplectic action of $G$ centralized by an antisymplectic involution $\sigma$ with $\mathrm{Fix}_X(\sigma) \neq \emptyset$.
\end{proposition}
\begin{proof}
Assume that a K3-surface with these properties exists. Applying Theorem \ref{roughclassi} we see that $X \to Y$ is branched along a single $G$-invariant smooth curve $C$ on the Del Pezzo surface $Y$. The genus of $C$ is at least four by Hurwitz' formula and therefore $\mathrm{deg}(Y) \geq 3$. We consider the quotient $Q$ of $C$ by the normal subgroup $N = C_2^3$. By Lemma \ref{S4 not on g=1,2} the quotient $Q$ is either rational or $g(Q) >2$. In the second case $g(C) \geq 19$ and we obtain a contradiction to $\mathrm{deg}(Y) = g(C)-1 \leq 9$. It follows that $Q$ is a rational curve.

We consider the action of $G$ on the Del Pezzo surface $Y$ of degree $\geq 3$, in particular the induced action on its configuration of (-1)-curves. By Remark \ref{stab of minus one curve} the stabilizer of a (-1)-curve in $Y$ has index $\geq 16$ in $G$ and we may immediately exclude the cases $\mathrm{deg}(Y) = 3,5,6,7$.
The automorphism group of a Del Pezzo surface of degree four is $C_2^4 \rtimes \Gamma$ for $\Gamma \in \{C_2, C_4, S_3, D_{10}\}$ (cf. \cite{dolgachev}). In particular, the maximal possible order is 160 and therefore $\mathrm{deg}(Y) \neq 4$.

Assume that $Y \cong \mathbb P_1 \times \mathbb P_1$. The canonical projection $\pi_{1,2}: Y \to \mathbb P_1$ is equivariant with respect to a subgroup $H$ of $G$ of index at most two. It follows that $H$ fits into the exact sequences
\begin{align*}
 \{\mathrm{id}\} \to I_1 \to H \overset{(\pi_1)_*}{\to} H_1\to  \{\mathrm{id}\} \\
 \{\mathrm{id}\} \to I_2 \to H \overset{(\pi_2)_*}{\to} H_2\to  \{\mathrm{id}\}
\end{align*}
where $I_i \cong C_2 \times C_2$ is the ineffectivity of the induced $H$-action on the base and $H_i \cong S_4$ (cf. proof of Lemma \ref{conicbundle}).
Since the action of $G$ on $\mathbb P_1 \times \mathbb P_1$ is effective by assumption, it follows that $I_2$ acts effectively on $\pi_1(\mathbb P_1 \times \mathbb P_1)$. We find a set of four points in $\pi_1(\mathbb P_1 \times \mathbb P_1)$ with nontrivial isotropy with respect to $I_2 \cong C_2 \times C_2$. Since $I_2$ is a normal subgroup of $H$, this set is $H$-invariant. The action of $H_1 \cong S_4$ on $\pi_1(\mathbb P_1 \times \mathbb P_1)$ does however not admit invariant sets of cardinality four since the minimal $S_4$-orbit in $\mathbb P_1$ has length six. 

We conclude that $Y$ must be isomorphic to $\mathbb P_2$. It follows that $g(C) =10$. Return to the covering $C \to Q$,
\[
 -18 = e(C) = 8\cdot e(Q) - \text{branch point contributions}.
\]
Since $Q$ is rational, the
branch point contribution must $34$. The possible isotropy of $N = C_2^3$ at a point in $C$ is $C_2$ and the full branch point contribution must be divisible by four. This contradiction yields the desired non-existence. 
\end{proof}
\section{The group $N_{72} = C_3^2 \rtimes D_8$}\index{Mukai group! $N_{72}$}\label{N72}
We let $X$ be a K3-surface with a symplectic action of $G=N_{72}$ centralized by an antisymplectic involution $\sigma$ with $\mathrm{Fix}_X(\sigma) \neq \emptyset$.
Note that in this case we may not apply Theorem \ref{roughclassi} and therefore begin by excluding that a $G$-minimal model of $Y=X/\sigma$ is an equivariant conic bundle.
\begin{lemma}\label{N72notconicbundle}
A $G$-minimal model of $Y$ is a Del Pezzo surface.
\end{lemma}
\begin{proof}
Assume the contrary and let $Y_\mathrm{min}$ be an equivariant conic bundle and a $G$-minimal model of $Y$.  We consider the induced action of $G$ on the base $B=\mathbb P_1$
and denote by $I \lhd G$ the ineffectivity of the $G$-action on $B$. Arguing as in the proof of Lemma \ref{conicbundle}, we see that $I$ is trivial or isomorphic to either $C_2$ or $C_2 \times C_2$. In all cases the quotient $G/I$ contains the subgroup $C_3 \times C_3$, which has no effective action on the rational curve $B$.
\end{proof}
As we will see, only very few Del Pezzo surfaces admit an effective action of the group $N_{72}$. We will explicitly use the group structure of $N_{72}= C_3^2 \rtimes D_8$: the action of $D_8 = C_2 \ltimes (C_2 \times C_2) = \langle \alpha \rangle  \ltimes ( \langle \beta \rangle \times \langle \gamma \rangle) = \mathrm{Aut}(C_3 \times C_3)$ on $C_3 \times C_3$ is given by
\[
\alpha(a,b) = (b,a), \quad \beta(a,b)=(a^2,b), \quad \gamma(a,b) = (a,b^2).
\]
As a first step we show:
\begin{lemma}
The degree of a Del Pezzo surface $Y_\mathrm{min}$ is at most four.
\end{lemma}
\begin{proof}
We exclude Del Pezzo surface of degree $\geq 5$.
\begin{itemize}
 \item 
A Del Pezzo surface of degree five has automorphims group $S_5$ and  $N_{72} \nless S_5$.
\item
The automorphism group of a Del Pezzo surface of degree six is $(\mathbb C^* )^2 \rtimes (S_3 \times C_2)$ (cf. Theorem 10.2.1 in \cite{dolgachev}). Assume that $N_{72} = C_3^2 \rtimes D_8$ is contained in this group and consider the intersection $ A = N_{72} \cap (\mathbb C^* )^2$. The quotient of $N_{72}$ by $A$ has order at most 12 and may not contain a copy of $C_3^2$. Therefore, the order of $A$ is at least six and $A$ contains a copy of $C_3$. If $|A| =6$, then $A = C_6 = C_3 \times C_2$ and $C_2$ is central in $N_{72}$. Using the group structure of $N_{72}$ specified above one finds that there is no copy of $C_2$ in $N_{72}$ centralizing $C_3 \times C_3$ and therefore $C_2$ cannot be contained in the centre of $N_{72}$.
For every choice of $C_3$ inside $C_3 \times C_3$ there is precisely one element in $\{\alpha, \beta, \gamma\}$ acting trivially on it and the centralizer of $C_3$ inside $D_8$ is isomorphic to $C_2$. If $|A|  >6$, then the centralizer of $C_3$ in $D_8$ has order greater then 2, a contradiction.
\item
A Del Pezzo surface of degree seven is obtained by blowing-up to points $p, q$ in $\mathbb P_2$. As was mentioned before, such a surface is never $G$-minimal.
\item
If $G$ acts on $\mathbb P_1 \times \mathbb P_1$, then the canonical projections are equivariant with respect to a subgroup $H$ of index two in $G$. We consider one of these projections. The action of $H$ induces an effective action of $H/I$ on the base $\mathbb P_1$. The group $I$ is either trivial or isomorphic to $C_2$ or $C_2 \times C_2$. In all case we find an effective action of $C_3 ^2$ on the base, a contradiction.  
\item
It remains to exclude $\mathbb P_2$.
If $N_{72}$ acts on $\mathbb P_2$ we consider its embedding into $\mathrm{PSL}_3(\mathbb C)$, in particular the realization of the subgroup $C_3^2 = \langle a \rangle \times \langle b \rangle$ and its lifting to $\mathrm{SL}_3(\mathbb C)$. 

We fix a preimage $\tilde a$ of $a$ inside $\mathrm{SL}_3(\mathbb C)$ and may assume that $\tilde a$ is diagonal. Since the action of $a$ on $\mathbb P_2$ is induced by a symplectic action on $X$, it follows that $a$ does not have a positive-dimensional set of fixed point. In appropiately chosen coordinates
\[
\tilde a=
 \begin{pmatrix}
1 & 0 & 0 \\
0 & \xi & 0\\
0&0& \xi^2
 \end{pmatrix}, 
\]
where $\xi$ is third root of unity.
As a next step, we want to specify a preimage $\tilde b$ of $b$ inside $\mathrm{SL}_3(\mathbb C)$. Since $a$ and $b$ commute in $\mathrm{PSL}_3(\mathbb C)$, we know that
\[
\tilde a \tilde b \tilde a^{-1} \tilde b^{-1} = \xi^k \mathrm{id}_{\mathbb C^3}
\]
for $k \in \{0, 1,2\}$.
Note that $\tilde b$ is not diagonal in the coordinates chosen above since this would 
 give rise to $C_3^2 $-fixed points in $\mathbb P_2$. As these correspond to $C_3^2$-fixed points on the double cover $X \to Y$ and a symplectic action of $C_3^2 \nless \mathrm{SL}_2(\mathbb C)$ on a K3-surface does not admit fixed points, this is a contradiction. An explicit calculation yields
\begin{align*}
\tilde b = \tilde b_1 =
 \begin{pmatrix}
0 & 0 & * \\
* & 0 & 0\\
0 & * & 0 
 \end{pmatrix} \quad \text{or} \quad
\tilde b = \tilde b_2= 
 \begin{pmatrix}
0 & * & 0 \\
0 & 0 & *\\
* & 0 & 0 
 \end{pmatrix}.
\end{align*}
We can introduce a change of coordinates commuting with $\tilde a$ such that
\begin{align*}
\tilde b = \tilde b_1 =
 \begin{pmatrix}
0 & 0 & 1 \\
1 & 0 & 0\\
0 & 1 & 0 
 \end{pmatrix} \quad \text{or} \quad
\tilde b = \tilde b_2= 
 \begin{pmatrix}
0 & 1 & 0 \\
0 & 0 & 1\\
1&0&0 
 \end{pmatrix}.
\end{align*}
Since $\tilde b_1^2 = \tilde b_2$, the two choices above correspond to the two choices of generators $b$ and $b^2$ of $\langle b \rangle$. We pick $\tilde b = \tilde b_1$.

The action of $D_8$ on $C_3^2 $ is specified above and the element $\beta \in D_8$ acts on $C_3^2 $ by $a \to a^2$ and $b \to b$. There is no element $T \in \mathrm{SL}_3(\mathbb C)$ such that (projectively) $T \tilde a T^{-1} = \tilde a^2$ and $T \tilde b T^{-1} = \tilde b$. It follows that there is no action of $N_{72}$ on $\mathbb P_2$.
\end{itemize}
This completes the proof of the lemma.
\end{proof}
As a next step, we study the possibility of rational curves in $\mathrm{Fix}_X(\sigma)$.
\begin{lemma}
There are no rational curves in $\mathrm{Fix}_X(\sigma)$.
\end{lemma}
\begin{proof}
 Let $n$ denote the total number of rational curves in $\mathrm{Fix}_X(\sigma)$ and recall $n \leq 10$. If $n \neq 0$, let $C$ be a rational curve in the image of  $\mathrm{Fix}_X(\sigma)$ in $Y$ and let $H = \mathrm{Stab}_G(C)$ be its stabilzer. The index of $H$ in $G$ is at most nine, therefore the order of $H$ is at least eight. The action of $H$ on $C$ is effective. 

First note that $G$ does not contain $S_4 = O_{24}$ as a subgroup. If this were the case, consider the intersection $S_4 \cap C_3^2$ and the quotient $S_4 \to S_4 / (S_4 \cap C_3^2) < D_8$. Since the only nontrivial normal subgroups of $S_4$ are $A_4$ and $C_2 \times C_2$, this leads to a contradiction.

Consequently, the order of $H$ is at most twelve. In particular, $n \geq 6$. Since $C_8 \nless G$, the group $H$ is not cyclic and any $H$-orbit on $C$ consists of at least two points.

It follows from $C^2 = -4$ that $C$ must meet the union of Mori fibers and the union of Mori fibers meets the curve $C$ in at least two points. Recalling that each Mori fibers meets the branch locus $B$ in at most two points we see that at least $n$ Mori fibers meeting $B$ are required. However, no configuration of $n$ Mori fibers is sufficient to transform the curve $C$ into a curve on a Del Pezzo surface and further Mori fibers are required. By invariance, the total number $m$ of Mori fibers must be at least $2n$.

Combining the Euler-characteristic formula 
\[
24 = 2e(Y_\mathrm{min}) +2m - 2n + \underset{\text{branch curve exists}}{\underset{\text{if non-rational }}{\underbrace{2g-2}}}
\]
with our observation $\mathrm{deg}(Y_\mathrm{min}) \leq 4$, i.e., $e(Y_\mathrm{min}) \geq 8$ we see that $n \leq 4$. However, it was shown above, that if $n \neq 0$, then $n \geq 6$. It follows that $n =0$.
\end{proof}
\begin{proposition}
 The quotient surface $Y$ is $G$-minimal and isomorphic to the Fermat cubic $\{x_1^3 + x_2^3 +x_3^3 +x_4^3 =0\} \subset \mathbb P_3$. Up to equivalence, there is a unique action of $G$ on $Y$ and the branch locus of $X \to Y$ is given by $\{x_1x_2 + x_3x_4 =0\}$. In particular, $X$ is equivariantly isomorphic to Mukai's $N_{72}$-example.
\end{proposition}
\begin{proof}
We first show that the total number $m$ of Mori fibers equals zero. By the Euler-characteristic formula above, the number $m$ is bounded by four. Using the fact that the maximal order of a stabilizer group of a Mori fiber is twelve (cf. proof of Theorem \ref{roughclassi}) we see that $Y$ must be $G$-minimal.

In order to conclude that $Y$ is the Fermat cubic we consult Dolgachev's lists of automorphisms groups of Del Pezzo surfaces of degree less than or equal to four (\cite{dolgachev} Section 10.2.2; Tables 10.3; 10,4; and 10.5):
It follows immediately from the order of $G$ that $Y$ is not of degree two or four. If $G$ were a subgroup of an automorphism group of a Del Pezzo surface of degree one, it would contain a central copy of $C_3$. The group structure of $N_{72}$ does however not allow this.
 After excluding the cases $\mathrm{deg}(Y) \in \{1,2,4\}$ the result now follows from the uniqueness of the cubic surface in $\mathbb P_3$ with an action of $N_{72}$ (cf. Appendix \ref{N72appendix}).
The action of $G$ on $Y$ is induced by a four-dimensional (projective) representation of $G$ and the branch curve $C \subset Y$ is the intersection of $Y$ with an invariant quadric (compare proof of Proposition \ref{S5 on degree three}).

In the Appendix \ref{N72appendix} it is shown that there is a uniquely determined action of $N_{72}$ on $\mathbb P_3$ and a unique invariant quadric hypersurface $\{x_1x_2 + x_3x_4 =0\}$. In particular, the branch curve in $Y$ is defined by $\{x_1x_2 + x_3x_4 =0\} \cap Y$. 

Mukai's $N_{72}$-example is defined by $\{ x_1^3+ x_2 ^3 + x_3^3 +x_4^3= x_1x_2 + x_3x_4+ x_5^2 = 0 \} \subset \mathbb P_4$. An anti-symplectic involution centralizing the action of $N_{72}$ is given by the map $ x_5 \mapsto -x_5$. The quotient of Mukai's example by this involution is the Fermat cubic and the fixed point set of the involution is given by $\{x_1x_2 + x_3x_4= 0 \}$.
\end{proof}
\section{The group $M_9 =C_3^2 \rtimes Q_8$}\index{Mukai group! $M_9$}\label{M9}
Let $G = M_9$ and let $X$ be a K3-surface with a symplectic $G$-action centralized by the antisymplectic involution $\sigma$ such that $\mathrm{Fix}_X(\sigma) \neq \emptyset$.
We proceed in analogy to the case $G=N_{72}$ above. Arguing precisely as in the proof of Lemma \ref{N72notconicbundle} one shows.
\begin{lemma}
A $G$-minimal model of $Y$ is a Del Pezzo surface.
\end{lemma}
We may exclude rational branch curves without studying configurations of Mori fibers.
\begin{lemma}\label{subgroupsM9}
There are no rational curves in $\mathrm{Fix}_X(\sigma)$.
\end{lemma}
\begin{proof}
Let $n$ be the total number of rational curves in $\mathrm{Fix}_X(\sigma)$. Assume $n \neq 0$, let $C$ be a rational curve in the image of $\mathrm{Fix}_X(\sigma)$ in $Y$ and let $H <G$ be its stabilizer. The action of $H$ on $C$ is effective. We go through the list of finite groups with an effective action on a rational curve. 

Since $M_9$ is a group of symplectic transformations on a K3-surface, its element have order at most eight.
Clearly, $A_6 \nless M_9$ and $D_{10}, \, D_{14}, \, D_{16} \nless M_9$. If $S_4 < M_9 = C_3^2 \rtimes Q_8$, then $S_4 \cap C_3^2$ is a normal subgroup of $S_4$ and it is therefore trivial. Now $ S_4 = S_4 / (S_4 \cap C_3^2) <  M_9 / C_3^2 =  Q_8$ yields a contradiction. The same argument can be carried out for $A_4$, $D_8$ and $C_8$.  If $D_{12} < M_9 = C_3^2 \rtimes Q_8$, then either $D_{12} \cap C_3^2 = C_3$ and $C_2 \times C_2 = D_{12} / C_3 < M_9 / C_3^2 =  Q_8$ or $D_{12} \cap C_3^2 = \{ \mathrm{id} \}$ and $D_{12} < Q_8$, both are impossible.

It follows that the subgroups of $M_9$ admitting an effective action on a rational curve have index greater than or equal to twelve. Therefore  $n \geq 12$, contrary to the bound $n \leq 10$ obtained in Corollary \ref{atmostten}.
\end{proof}
\begin{proposition}
 The quotient surface $Y$ is $G$-minimal and isomorphic to $\mathbb P_2$. Up to equivalence, there is a unique action of $G$ on $Y$ and the branch locus of $X \to Y$ is given by
$\{x_1^6 + x_2^6+ x_3^6-10( x_1^3x_2^3 + x_2^3x_3^3+ x_3^3x_1^3 ) =0\}$. In particular, $X$ is equivariantly isomorphic to Mukai's $M_9$-example.
\end{proposition}
\begin{proof}
We first check that $Y$ is $G$-minimal. Again, we proceed as in the proof of Theorem \ref{roughclassi}   and Lemma \ref{subgroupsM9} above to see that the largest possible stablizer group of a Mori fiber is $D_6 < G$. If $Y$ is not $G$-minimal, this implies that the total number of Mori fibers is $ \geq 12$, contradicting $m \leq 9$.

Note that $X \to Y$ is not branched along one or two elliptic curves as this would imply $e(Y) =12$ and contradict the fact that $Y$ is a Del Pezzo surface. 

Let $D$ be the branch curve of $ X \to  Y$ and consider the quotient $Q$ of $D$ by the normal subgroup $N = C_3^2$ in $G$. On $Q$ there is an action of $Q_8$ implying that $Q$ is not rational. 
We show that $Q_8$ does not act on an elliptic curve $Q$. If this were the case, consider
the decomposition $Q_8 = (Q_8 \cap Q) \rtimes (Q_8 \cap L)$ where $(Q_8 \cap L)$ is a nontrivial cyclic group. For any choice of generator of $(Q_8 \cap L)$ the center $\{+1,-1\}$ of $Q_8$ is contained in $(Q_8 \cap L)$.
Let $q: Q_8 \to Q_8 /(Q_8 \cap Q) \cong Q_8 \cap L $ denote the quotient homomorphism. The commutator subgroup $Q_8' = \{+1,-1\}$ must be contained in the kernel of $q$. This contradiction yields that $Q_8$ does not act on an elliptic curve. 
 It follows that the genus of $Q$ is at least two and the genus of $D$ is at least ten. Adjunction on the Del Pezzo surface $Y$ now implies $g=10$ and $Y \cong \mathbb P_2$. 

It is shown in Appendix \ref{M9 on P2} that,
up to natural equivalence, there is a unique action of $M_9$ on the projective plane. In suitably chosen coordinated the generators $a, b$ of $C_3^2$ are represented as
\begin{align*}
\tilde a=
 \begin{pmatrix}
1 & 0 & 0 \\
0 & \xi & 0\\
0&0& \xi^2
 \end{pmatrix}, \quad
\tilde b= 
 \begin{pmatrix}
0 & 1 & 0 \\
0 & 0 & 1\\
1&0&0 
 \end{pmatrix}
\end{align*}
and $I,J \in Q_8$ are represented as
\begin{align*}
\tilde I= \frac{1}{\xi -\xi^2}
 \begin{pmatrix}
1 & 1 & 1 \\
1 & \xi & \xi^2\\
1 & \xi^2& \xi
 \end{pmatrix}, \quad
\tilde J= \frac{1}{\xi -\xi^2}
 \begin{pmatrix}
1 & \xi & \xi \\
\xi^2 & \xi & \xi^2\\
\xi^2 & \xi^2& \xi
 \end{pmatrix}.
\end{align*}
 We study the action of $M_9$ on then space of sextic curves. By restricting our consideration to the subgroup $C_3^2$ first, we see that a polynomial defining an invariant curve must be a linear combination of the following polynomials:
\begin{align*}
f_1 &=  x_1^6 + x_2^6+ x_3^6;\\
f_2 &= x_1^2x_2^2x_3^2;\\
f_3 &= x_1^3x_2^3 + x_1^3x_3^3+ x_2^3x_3^3;\\
f_4 &= x_1^4x_2x_3 + x_1x_2^4x_3+ x_1x_2x_3^4.
\end{align*}
Taking now the additional symmetries into account, we find three $M_9$-invariant sextic curves, namely
\[ 
\{f_1- 10f_3 =x_1^6 + x_2^6+ x_3^6-10( x_1^3x_2^3 + x_2^3x_3^3+ x_3^3x_1^3 ) =0\},
\]
which is the example found by Mukai, and additionally 
\[ 
\{f_a=f_1 + (18-3a) f_2 +2f_3 + af_4  =0\},
\]
where $a$ is a solution of the quadratic equation $a^2-6a+36$, i.e. $a= -6\xi$ or $a= -6\xi^2$. The polynomial $f_a$ is invariant with respect to the action of $M_9$ for $a= -6\xi^2$ and semi-invariant if $a= -6\xi$.

We wish to show that $X$ is not the double cover of $\mathbb P_2$ branched along $\{f_a=0\}$. If this were the case, 
consider the fixed point $p=[0:1:-1]$ of the automorphism $I$ and note that $f_a(p)=0$. So the $\pi^{-1}(p)$ consists of one point $x \in X$ and we linearize the $\langle I \rangle \times \langle \sigma \rangle$ at $x$. In suitably chosen coordintes the action of the symplectic automorphism $I$ of order four is of the form  $(z,w) \mapsto (iz, -iw)$. Since the action of $\sigma$ commutes with $I$, the $\sigma$-quotient of $X$ is locally given by
\[
(z,w) \mapsto (z^2, w) \quad \text{or}\quad (z,w) \mapsto (z, w^2).
\]
It follows that the action of $I$ on $Y$ is locally given by either
\[
\begin {pmatrix}
-1 & 0\\
0 & -i
\end {pmatrix}
\quad \text{or} \quad
\begin {pmatrix}
i & 0\\
0 & -1
\end {pmatrix}.
\]
In particular, the local linearization of $I$ at $p$ has determinant $\neq 1$. 
By a direct computation using the explicit form of $\tilde I$ given above, in particular the facts that $\mathrm{det}(\tilde I) =1$ and $\tilde I v =v$ for $[v]=p$, we obtain a contradiction.

 This completes the proof of the proposition.
\end{proof}
\begin{remark}\label{M9 symplectic}
In the proof of the propostion above we have observed that an element of $\mathrm{SL}_3(\mathbb C)$ does not necessarily lift to a symplectic transformation on the double cover of $\mathbb P_2$ branched along a sextic given by an invariant polynomial.
Mukai's $M_9$-example $X$ is a double cover of $\mathbb P_2$ branched along the sextic curve $\{x_1^6 + x_2^6+ x_3^6-10( x_1^3x_2^3 + x_2^3x_3^3+ x_3^3x_1^3 ) =0\}$ and for this particular example, the action of $M_9$ does lift to a group of symplectic transformation as claimed by Mukai. 

To see this consider the set $\{a,b,I,J\}$ of generators of $M_9$. Since $a$ and $b$ are commutators in $M_9$, they can be lifted to symplectic transformation $\overline a, \overline b$ on $X$. For $I,J$ consider the linearization at the fixed point $[0:1:-1]$ and check that it has determinant one. Since $[0:1:-1]$ is \emph{not} contained in the branch set of the covering, its preimage in $X$ consists of two points $p_1,p_2$. We can lift $I$ ($J$, respectively) to a transformation of $X$ fixing both $p_1,p_2$ and a neighbourhood of $p_1$ is $I$-equivariantly isomorphic to a neighbourhood of $ [0:1:-1] \in \mathbb P_2$. In particular, the action of the lifted element $\overline I$ ($\overline J$, respectively) is symplectic. On $X$ there is the action of a degree two central extension $E$ of $M_9$,
\[
\{\mathrm{id}\} \to C_2 \to E \to M_9 \to \{\mathrm{id}\}.
\]
The elements $\overline a, \overline b, \overline I, \overline J$ generate a subgroup $\tilde M_9$ of $E_\mathrm{symp}$ mapping onto $M_9$. Since $E_\mathrm{symp} \neq E$, the order of $\tilde M_9$ is 72 and it follows that $\tilde M_9$ is isomorphic to $M_9$. In particular $E$ splits as $E_\mathrm{symp} \times C_2$ with  $E_\mathrm{symp}= M_9$.
\end{remark}

\section{The group $T_{48} = Q_8 \rtimes S_3$}\index{Mukai group! $T_{48}$}
We let $X$ be a K3-surface with an action of $T_{48} \times C_2$ where the action of $G = T_{48}$ is symplectic and the generator $\sigma$ of $C_2$ is antisymplectic and has fixed points. The action of $S_3$ on $Q_8$ is given as follows: The element $c$ of order three in $S_3$ acts on $Q_8$ by permuting $I,J,K$ and an element $d$ of order two acts by exchanging $I$ and $J$ and mapping $K$ to $-K$.
\begin{lemma}\label{not two}
A $G$-minimal model $Y_\mathrm{min}$ of $Y$ is either $\mathbb P_2$, a Hirzebruch surface $\Sigma_n$ with $n >2$, or $e(Y_\mathrm{min}) \geq 9$.
\end{lemma}
\begin{proof}
Let us first consider the case where $ Y_\mathrm{min}$ is a Del Pezzo surface and go through the list of possibilities.
\begin{itemize}
 \item 
Let $Y_\mathrm{min} \cong \mathbb P_1 \times \mathbb P_1$. Since $T_{48}$ acts on $Y_\mathrm{min}$, both canonical projections are equivariant with respect to the index two subgroup $G'= Q_8 \rtimes C_3$. 
Since $Q_8$ has no effective action on $\mathbb P_1$, it follows that the subgroup $Z =\{+1,-1\} < Q_8$ acts trivially on the base. Since this holds with respect to both projections, the subgroup $Z$ acts trivially on $Y_\mathrm{min}$, a contradiction. 
\item
Using the group structure of $T_{48}$ one checks that the only nontrivial normal subgroup $N$ of $T_{48}$ such that $N \cap Q_8 \neq Q_8$ is the center $Z= \{+1,-1\}$ of $T_{48}$. It follows that $T_{48}$ is neither a subgroup of $(\mathbb C^*)^2 \rtimes (S_3 \times C_2)$ nor a subgroup of any of the automorphism groups $C_2^4 \rtimes \Gamma$ for $\Gamma \in \{C_2, C_4, S_3, D_{10}\}$ of a Del Pezzo surface of degree four. 
Furthermore,
$T_ {48} \nless S_5$. Thus it follows that $d(Y_\mathrm{min}) \neq 4,5,6$.
\end{itemize}
So if $Y_\mathrm{min}$ is a Del Pezzo surface, then $Y_\mathrm{min} \cong \mathbb P_2$ or $e(Y_\mathrm{min}) \geq 9$ .

Let us now turn to the case where $Y_\mathrm{min}$ is an equivariant conic bundle.
We first show that $Y_\mathrm{min}$ is not a conic bundle with singular fibers. We assume the contrary and let $p: Y_\mathrm{min} \to \mathbb P_1$ be an equivariant conic bundle with singular fibers. 
The center $Z = \{+1,-1\}$ of $G= T_{48}$ acts trivially on the base an has two fixed points in the generic fiber.
Let $C_1$ and $C_2$ denote the two curves of $Z$-fixed points in $Y_\mathrm{min}$. By Lemma \ref{singular fibers of conic bundle} any singular fiber $F$ is the union of two (-1)-curves $F_1,F_2$ meeting transversally in one point. We consider the action of $Z$ on this union of curves. The group $Z$ does not act trivially on either component of $F$ since linearization at a smooth point of $F$ would yield a trivial action of $Z$ on $Y_\mathrm{min}$.
Consequently,
it has either one or three fixed points on $F$. The first is impossible since $C_1$ and $C_2$ intersect $F$ in two points. It follows that $Z$ stabilizes each curve $F_i$. We linearize the action of $Z$ at the point of intersection $F_1 \cap F_2$. The intersection is transversal and the action of $Z$ is by $-\mathrm{Id}$ on $T_{F_1} \oplus T_{F_2}$ contradicting the fact the $Z$ acts trivially on the base. Thus $Y_{min}$ is not a conic bundle with singular fibers.

If $Y_\mathrm{min} \to \mathbb P_1$ is a Hirzebruch surface $\Sigma_n$, then the action of $T_{48}$ induces an effective action of $S_4$ on the base $\mathbb P_1$. 

The action of $T_{48}$ on $\Sigma_n$ stabilizes two disjoint sections $E_\infty$ and $E_0$, the curves of $Z$-fixed points. This is only possible if $E_0^2 = -E_\infty^2 = n$. Removing the exceptional section $E_\infty$ from $\Sigma_n$, we obtain the hyperplane bundle $H^n$ of $\mathbb P_1$. Since $T_{48}$ stabilizes the section $E_0$, we chose this section to be the zero section and conclude that the action of $T_{48}$ on $H^n$ is by bundle automorphisms. 

If $n=2$, then $H^n$ is the anticanonical line bundle of $\mathbb P_1$ and the action of $S_4$ on the base induces an action of $S_4$ on $H^2$ by bundle automorphisms. It follows that $T_{48}$ splits as $S_4 \times C_2$, a contradiction. Thus, if $Y_\mathrm{min}$ is a Hirzebruch surface $\Sigma_n$, then $n \neq 2$. 
\end{proof}
\begin{lemma}\label{no rat T48}
 There are no rational curves in $\mathrm{Fix}_X(\sigma)$.
\end{lemma}
\begin{proof}
We let $n$ denote the total number of rational curves in $\mathrm{Fix}_X(\sigma)$ and assume $n >0$. Recall $n \leq 10$, let $C$ be a rational curve in $B=\pi(\mathrm{Fix}_X(\sigma)) \subset Y$ and let $H = \mathrm{Stab}_G(C) < G$ be its stabilizer group. The action of $H$ on $C$ is effective, the index of $H$ in $G$ is at most 8. Using the quotient homomorphism $T_{48} \to T_{48}/Q_8 = S_3$ one checks that $T_{48}$ does not contain $O_{24}=S_4$ or $T_{12}= A_4$ as a subgroup. It follows that $H$ is a cyclic or a dihedral group.

If $H\in \{C_6, C_8, D_8\}$, then $H$ and all conjugates of $H$ in $G$ contain the center $Z= \{+1,-1\}$ of $G$. It follows that $Z$ has two fixed point on each curve $gC$ for $g \in G$. Since there are six (or eight) distinct curves $gC$ in $Y$, it follows that $Z$ has at least 12 fixed points in $Y$ and in $X$. This contradicts to assumption that $Z < G$ acts symplectically on $X$ and therefore has eight fixed points in the K3-surface $X$.

It remains to study the cases $H = D_{12}$ and $H = D_6$ where $n = 8$ or $n=4$. 

We note that a Hirzebruch surface has precisely one curve with negative self-intersection and only fibers have self-intersection zero. A Del Pezzo surface does not contains curves of self-intersection less than $-1$. The rational branch curves must therefore meet the union of Mori fibers in $Y$.

The total number of Mori fibers is bounded by $n+9$. We study the possible stabilizer subgroups $\mathrm{Stab}_G(E) < G$ of Mori fibers. A Mori fiber $E$ with self-intersection (-1) meets the branch locus $B$ in one or two points and its stabilizer is either cyclic or dihedral. If $\mathrm{Stab}_G(E) \in \{C_4, D_8\}$, then the points of intersection of $E$ and $B$ are fixed points of the center $Z$ of $G$ and we find too many $Z$-fixed points on $X$. 

Assume $n=4$ and let $R_1, \dots R_4$ be the rational curves in $B$. We denote by $\tilde R_i$ their images in $Y_\mathrm{min}$. The total number $m$ of Mori fibers is bounded by 12.
We go through the list of possible configurations:
\begin{itemize}
 \item 
If $m = 4$, there is no invariant configuration of Mori fibers such that the contraction maps the four rational branch curves to a configuration on the Hirzebruch or Del Pezzo surface $Y_\mathrm{min}$.
\item
If $m= 6$, then $\mathrm{Stab}_G(E) = C_8$ and the points of intersection of $E$ and $B$ are $Z$-fixed. Since $Z$ has at most eight fixed points on $B$, it follows that each curve $E$ meets $B$ only once. The images $\tilde R_i$ of the $R_i$ contradict our observations about curves in Del Pezzo and Hirzebruch surfaces.
\item
If $m=8$ and all Mori fibers have self-intersection $-1$, then each Mori fiber meets $\bigcup R_i$ in a $Z$-fixed point.
 Since there at at most eight such points, it follows that each Mori fibers meets $\bigcup R_i$ only once and their contractions does not transform the curves $R_i$ sufficiently. 
\item
If $m=8$ and only four Mori fibers have self-intersection $-1$, we consider the four Mori fibers of the second reduction step. Each of these meets a Mori fiber $E$ of the first step in precisely one point. By invariance, this would have to be a fixed points of the stabilizer $ \mathrm{Stab}_G(E)= D_{12}$, a contradiction.
\item
If $m=12$, then either $e(Y_\mathrm{min}) = 3$ and there exist a branch curve $D_{g=2}$ of genus two or $e(Y_\mathrm{min}) = 4$ and $B = \bigcup R_i$. In the first case, $Y_\mathrm{min}  \cong \mathbb P_2$ and twelve Mori fibers are not sufficient to transform $B = D_{g=2} \cup \bigcup R_i$ into a configuration of curves in the projective plane.
So $Y_\mathrm{min} = \Sigma_n$ for $n > 2$. 

 Since $Z$ has two fixed points in each fiber of $p: \Sigma_n \to \mathbb P_1$ the $Z$-action on $\Sigma_n$ has two disjoint curves of fixed points. As was remarked above, these curves are the exceptional section $E_\infty$ of self-intersection -$n$ and a section $E_0 \sim E_\infty + n F$ of self-intersection $n$. Here $F$ denotes a fiber of $p: \Sigma_n \to \mathbb P_1$. There is no automorphisms of $\Sigma_n$ mapping $E_\infty$ to $E_0$. 

Each rational branch curve $\tilde R_i$ has two $Z$-fixed points. These are exchanged by an element of $\mathrm{Stab}_G(R_i)$ and therefore both lie on either $E_\infty$ or $E_0$, i.e., $\tilde R_i$ cannot have nontrivial intersection with both $E_0$ and $E_\infty$. By invariance all curves $\tilde R_i$ either meet $E_0$ or $E_\infty$ and not both.

Using the fact that $\sum \tilde R_i$ is linearly equivalent to $-2K_{\Sigma_n} \sim 4 E_\infty +(2n +4)F$ we find that $\tilde R_i \cdot E_\infty = 0$ and $n=2$, a contradiction to Lemma \ref{not two}. 
\end{itemize}
We have shown that
all possible configurations in the case $n \neq 4$ lead to a contradiction. We now turn to the case $n=8$ and let $R_1, \dots R_8$ be the rational ramification curves. The total number of Mori fibers is bounded by 16. Note that by invariance, the orbit of a Mori fiber meets $\bigcup R_i$ in at least 16 points or not at all. In particular, Mori fibers meeting $R_i$ come in orbits of length $\geq 8$. As above, we go through the list of possible configurations.
\begin{itemize}
 \item 
If $m =16$, then the set of all Mori fibers consists of two orbits of length eight. If all 16 Mori fibers meet $B$, then each meets $B$ in one point and $R_i$ is mapped to a (-2)-curve in $Y_\mathrm{min}$. If only eight Mori fibers meet $B$, then each of the eight Mori fibers of the second reduction step meets one Mori fiber $E$ of the first reduction step in one point. This point would have to be a $\mathrm{Stab}_G(E)$-fixed point. But if $\mathrm{Stab}_G(E)$ is cyclic, its fixed points coincide with the points $E\cap B$. 
\item
If $m=12$, then the set of all Mori fibers consists of a single $G$-orbit and each curve $R_i$ meets three distinct Mori fibers. Their contraction transforms $R_i$ into a (-1)-curve on $Y_\mathrm{min}$. It follows that $Y_\mathrm{min}$ contains at least eight (-1)-curves and is a Del Pezzo surface of degree $\leq 5$. We have seen above that $d( Y_\mathrm{min}) \neq 4,5$ and therefore $e(Y_\mathrm{min}) \geq 9$. With $m=12$ and $n=8$, this contradicts the Euler characteristic formula $24 = 2 e(Y_\mathrm{min}) +2m -2n +(2g-2)$.
\item
If $m=8$ there is no invariant configuration of Mori fibers such that the contraction maps the eight rational branch curves to a configuration on the Hirzebruch or Del Pezzo surface $Y_\mathrm{min}$
\end{itemize}
This completes the proof of the lemma.
\end{proof}
Since there is an effective action of $T_{48}$ on $\mathrm{Fix}_X(\sigma)$, it is neither an elliptic curve nor the union of two elliptic curves. It follows that $ X \to Y$ is branched along a single $T_{48}$-invariant curve $B$ with $g(B) \geq 2$.
\begin{lemma}
 The genus of $B$ is neither three nor four.
\end{lemma}
\begin{proof}
We consider the quotient $Q = B/Z$ of the curve $B$ by the center $Z$ of $G$ and apply the Euler characteristic formula,  $e(B) = 2 e(Q) - |\mathrm{Fix}_B(Z)|$. On $Q$ there is an effective action of the group $G/Z = (C_2 \times C_2) \rtimes C_3 = S_4$. Using Lemma \ref{S4 not on g=1,2} we see that $e(Q) \in \{2, -4, -6,-8 \dots\}$.

If $g(B)=3$, then $e(B) = -4$ and the only possibility is $Q \cong \mathbb P_1$ and $|\mathrm{Fix}_B(Z)| =8$. In particular, all $Z$-fixed points on $X$ are contained in the curve $B$. Let $A < G$ be the group generated by $I \in Q_8 = \{ \pm 1, \pm I, \pm J, \pm K \}$. The four fixed points of $A$ in $X$ are contained in $\mathrm{Fix}_X(Z) = \mathrm{Fix}_B(Z)$ and the quotient group $A/Z \cong C_2$ has four fixed points in $Q$. This is a contradiction.

If $g(B)=4$, then $e(B) = -6$ and the only possibility is $Q \cong \mathbb P_1$ and $|\mathrm{Fix}_B(Z)| =10$. This contradicts the fact that $Z$ has at most eight fixed points in $B$ since it has precisely eight fixed points in $X$.
\end{proof}
In Lemma \ref{not two} we have reduced the classification to the cases $e(Y_\mathrm{min}) \in \{3,4, 9, 10, 11\}$. In the following, we will exclude the cases $e(Y_\mathrm{min}) \in \{4, 9, 10,\}$ and describe the remaining cases more precisely. Recall that the maximal possible stabilizer subgroup of a Mori fiber is $D_{12}$, in particular, $m = 0$ or $ m \geq 4$.
\begin{lemma}
 If  $e(Y_\mathrm{min}) =3$, then $Y_\mathrm{min} = Y = \mathbb P_2$ and $X \to Y$ is branched along the curve $\{x_1x_2(x_1^4-x_2^4) + x_3^6=0\}$. In particular, $Y$ is equivariantly isomorphic to Mukai's $T_{48}$-example.
\end{lemma}
\begin{proof}
Let $ M: Y \to \mathbb P_2$ denote a Mori reduction of $Y$ and let $B \subset Y$ be the branch curve of the covering $X \to Y$. 
If $Y = Y_\mathrm{min}$, then $B= M(B)$ is a smooth sextic curve.
If $Y \neq Y_\mathrm{min}$, then the Euler characteristic formula with $m \in \{4,6,8\}$ shows that $g(B) \in \{2,4,6\}$. The case $m=6$, $g(B) =4$ has been excluded by the previous lemma.

If $m=4$, then the stabilizer group of each Mori fiber is $D_{12}$ and each Mori fiber meets $B$ in two points. Furthermore, since in this case $g(B) =6$, the self-intersection of $\mathrm{Fix}_X(\sigma)$ in $X$ equals ten and therefore $B^2 = 20$. The image $M(B)$ of $B$ in $Y_\mathrm{min}$ has self-intersection $20 + 4 \cdot 4 = 36$ and follows to be an irreducible singular sextic.

If $m=8$, then $g(B) =2$ and $B^2 =4$. Since the self-intersection number $M(B)^2$ must be a square, one checks that all possible invariant configurations of Mori fibers yield $M(B)^2 =36$ and involve Mori fibers meeting $B$ is two points. In particular, $M(B)$ is a singular sextic.

We study the action of $T_{48}$ on the projective plane. As a first step, we may choose coordinates on $\mathbb P_2$ such that the automorphism $-1 \in Q_8 < T_{48}$ is represented as
\[
\widetilde {-1}=
 \begin{pmatrix}
-1 & 0 & 0 \\
0 & -1 & 0\\
0&0& 1
 \end{pmatrix}.
\]
We denote by $V$ to the $-1$-eigenspace of this operator.
For each element $I, J , K$ there is a unique choice $\widetilde I, \widetilde J, \widetilde K$ in $\mathrm{SL}_3(\mathbb C)$ such that $\widetilde I ^2  =  \widetilde J^2 =  \widetilde K^2 = \widetilde {-1}$. One checks $\widetilde I \widetilde J\widetilde K = \widetilde {-1}$. 
Therefore $\widetilde I, \widetilde J, \widetilde K$ generate a subgroup of $\mathrm{SL}_3(\mathbb C)$ isomorphic to $Q_8$.
By construction $\widetilde I, \widetilde J, \widetilde K$ stabilze the vector space $V$. Up to isomorphisms, there is a unique faithful 2-dimensional representation of $Q_8$ and it follows that $I,J,K$ are represented as
\begin{align*}
\widetilde I = 
 \begin{pmatrix}
-i & 0 & 0 \\
0 & i & 0\\
0&0& 1
 \end{pmatrix}, \quad
\widetilde J= 
 \begin{pmatrix}
0 & -1 & 0 \\
1 & 0 & 0\\
0 & 0  & 1
 \end{pmatrix}, \quad
\widetilde K= 
 \begin{pmatrix}
0 & i & 0 \\
i & 0 & 0\\
0 & 0  & 1
 \end{pmatrix}. 
\end{align*}
We recall that the action of $S_3$ on $Q_8$ is given as follows: The element $c$ of order three in $S_3$ acts on $Q_8$ by permuting $I,J,K$ and an element $d$ of order two acts by exchanging $I$ and $J$ and mapping $K$ to $-K$.
With $\mu = \sqrt{\frac{i}{2}}$ and $ \nu = \frac{i}{\sqrt{2}}$ it follows that the elements $c$ and $d$ are represented as
\begin{align*}
\widetilde c = 
 \begin{pmatrix}
-i\mu & i\mu& 0 \\
\mu & \mu & 0\\
0&0& 1
 \end{pmatrix}, \quad
\widetilde d= 
 \begin{pmatrix}
-i\nu & -\nu & 0 \\
\nu & i\nu & 0\\
0 & 0  & -1
 \end{pmatrix}.
\end{align*}
In particular, there is a unique action of $T_{48}$ on $\mathbb P_2$. In the following, we denote by $[x_1:x_2:x_3]$ homogeneous coordintes such that the action of $T_{48}$ is as above. Using the explicit form of the $T_{48}$-action and the fact that the commutator subgroup of $T_{48}$ is $Q_8 \rtimes C_3$ one can check that any invariant curve of degree six is of the form 
\[
C_\lambda = \{ x_1x_2(x_1^4-x_2^4) + \lambda x_3^6 =0\}
\]
In order to avoid this calculation, one can also argue that the polynomial $x_1x_2(x_1^4-x_2^4)$ is the lowest order invariant of the octahedral group $S_4 \cong T_{48} / Z$. 

The curve $C_\lambda$ is smooth and it follows that $Y = Y_\mathrm{min}$.
We may adjust the coordinates equivariantly such that $\lambda =1$ and find that our surface $X$ is precisely Mukai's $T_{48}$-example.
\end{proof}
\begin{remark}\label{T48 symplectic}
As claimed by Mukai, the action of $T_{48}$ on $\mathbb P_2$ does indeed lift to a symplectic action of $T_{48}$ on the double cover of $\mathbb P_2$ branched along the invariant curve $\{ x_1x_2(x_1^4-x_2^4) +  x_3^6 =0\}$. The elements of the commutator subgroup can be lifted to symplectic transformation on the double cover $X$. 

The remaining generator $d$ 
is an involution fixing the point $[0:0:1]$. Any involution $\tau$ with a fixed point $p$ outside the branch locus can be lifted to a symplectic involution on the double cover $X$ as follows:

The linearized action of $\tau$ at $p$ has determinant $\pm 1$. We consider the lifting $\tilde \tau$ of $\tau$ fixing both points in the preimage of $p$. Its linearization coincides with the linearization on the base and therefore also has determinant $\pm 1$. 
In particular, $\tilde \tau$ is an involution. 
It follows that either $\tilde \tau$ or the second choice of a lifting $\sigma \tilde \tau$ acts symplectically on $X$. 

The group generated by all lifted automorphisms is either isomorphic to $T_{48}$ or to the full central extension $E$
\[
 \{\mathrm{id}\} \to C_2 \to E \to T_{48} \to  \{\mathrm{id}\} 
\]
acting on the double cover.
Since $E_\mathrm{symp} \neq E$ the later is impossible it follows that $E$ splits as $E_\mathrm{symp} \times C_2$ with $E_\mathrm{symp} = T_{48}$.
\end{remark}
Finally, we return to the remaining possibilities $e(Y_\mathrm{min}) \in \{4, 9, 10, 11\}$.
\begin{lemma}
$e(Y_\mathrm{min}) \not\in \{4,9,10\}$.
\end{lemma}
\begin{proof}
Recalling that the genus of the branch curve $B$ is neither three nor four and that $m$ is either zero or $\geq 4$, we may exclude $e(Y_\mathrm{min}) =9,10$ using the Euler characteristic formula $12 = e(Y_\mathrm{min}) +m +g-1$. It remains to consider the case $Y_\mathrm{min} = \Sigma _n$ with $n >2$ and we claim that this is impossible.

Let $M= Y \to Y_\mathrm{min} = \Sigma_n$ denote a (possibly trivial) Mori reduction of $Y$. The image $M(B)$ of $B$ in $\Sigma _n$ is linearly equivalent to $-2K_{\Sigma _n}$. Now $M(B) \cdot E_\infty = 2(2-n) < 0$ and it follows that $M(B)$ contains the rational curve $E_\infty$. This is a contradiction since $B$ does not contain any rational curves by Lemma \ref{no rat T48}.
\end{proof}
In the last remaining case, i.e., $e(Y_\mathrm{min})=11$, the quotient surface $Y$ is a $G$-minimal Del Pezzo surface of degree 1. Consulting \cite{dolgachev}, Table 10.5, we find that $Y$ is a hypersurface in weighted projective space $\mathbb P(1,1,2,3)$ defined by the degree six equation
\[
 x_0x_1(x_0^4-x_1^4)+ x_2^3+x_3^2.
\]
This follows from the invariant theory of the group $S_4 \cong T_{48}/Z$ and fact that $Y$ is a double cover of a quadric cone $Q$ in $\mathbb P_3$ branched along the intersection of $Q$ with a cubic hypersurface (cf. Theorem \ref{antican models of del pezzo}).

The linear system of the anticanonical divisor $K_Y$ has precisely one base point $p$. In coordinates $[x_0:x_1:x_2:x_3]$ this point is given as $[0:0:1:i]$. It is fixed by the action of $T_{48}$. The linearization of $T_{48}$ at $p$ is given by the unique faithful 2-dimensional represention of $T_{48}$. This represention has implicitly been discussed above as a subrepresentation $V$ of the three-dimensional representation of $T_{48}$. It follows that there is a unique action of $T_{48}$ on $Y$. The branch curve $B$ is linearly equivalent to $-2K_Y$, i.e., $B = \{s=0\}$ for a section $s \in \Gamma(Y, \mathcal O(-2K_Y))$ which is either invariant or semi-invariant.

By an adjunction formula for hypersurfaces in weighted projective space $\mathcal O (-2K_Y)) = \mathcal O _Y(2)$. The four-dimensional space of sections $\Gamma(Y, \mathcal O(-2K_Y))$ is generated by the weighted homogeneous polynomials $x_0^2, x_1^2, x_0x_1, x_2$. 
We consider the map $ Y \to \mathbb P (\Gamma(Y, \mathcal O(-2K_Y))^*)$ associated to $|-2K_Y|$. Since this map is equivariant with respect to $\mathrm{Aut}(Y)$, the fixed point $p$ is mapped to a fixed point in $\mathbb P (\Gamma(Y, \mathcal O(-2K_Y))^*)$. It follows that the section corresponding to the homogeneous polynomial $x_2$ is invariant or semi-invariant with respect to $T_{48}$. It is the only section of $\mathcal O(-2K_Y)$ with this property since the representation of $T_{48}$ on the span of $x_0^2, x_1^2, x_0x_1$ is irreducible.

The curve $B \subset Y$ defined by $s=0$ is connected and has arithmetic genus $2$. Since $T_{48}$ acts effectively on $B$ and does not act on $\mathbb P_1$ or a torus, it follows that $B$ is nonsingular.

It remains to check that the action of $T_{48}$ on $Y$ lifts to a group of symplectic transformation on the double cover $X$ branched along $B$. First note that $B$ does not contain the base point $p$. 
For $I,J,K, c \in T_{48}$ we we choose liftings $ \overline I, \overline J, \overline K , \overline c \in \mathrm{Aut}(X)$ fixing both points in $\pi^{-1}(p)= \{p_1,p_2\}$. The linearization of $ \overline I, \overline J, \overline K , \overline c$ at $p_1$ is the same as the linearization at $p$ and in particular has determinant one. 
By the general considerations in Remark \ref{T48 symplectic} the involution $d$ can be lifted to a symplectic involution on $X$.
The symplectic liftings of $I,J,K,c,d$ generate a subgroup $\tilde G$ of $\mathrm{Aut}(X)$ which is isomorphic to either $T_{48}$ or to the central degree two extension of $T_{48}$ acting on $X$. 
In analogy to Remarks \ref{M9 symplectic} and \ref{T48 symplectic} we conclude that $\tilde G \cong T_{48}$ and the action of $T_{48}$ on $Y$ induces a symplectic action of $T_{48}$ on the double cover $X$.

This completes the classification of K3-surfaces with $T_{48} \times C_2$-symmetry. We have shown:
\begin{theorem}
 Let $X$ be a K3-surface with a symplectic action of the group $T_{48}$ centralized by an antisymplectic involution $\sigma$ with $\mathrm{Fix}_X(\sigma) \neq \emptyset$. Then $X$ is equivariantly isomorphic either to Mukai's $T_{48}$-example or to the double cover of 
\[
 \{x_0x_1(x_0^4-x_1^4)+ x_2^3+x_3^2=0\} \subset \mathbb P(1,1,2,3)
\]
branched along $\{x_2=0\}$
\end{theorem}
\begin{remark}
The automorphism group of the Del Pezzo surface $Y = \{x_0x_1(x_0^4-x_1^4)+ x_2^3+x_3^2=0\} \subset \mathbb P(1,1,2,3)$ is the trivial central extension $C_3 \times T_{48}$.  By contruction, the curve $B=\{s=0\}$ is invariant with respect to the full automorphism group. The double cover $X$ of $Y$ branched along $B$ carries the action of a finite group $\tilde G$ of order $2 \cdot 3 \cdot 48 = 288$ containing $T_{48} < \tilde G_\mathrm{symp}$. Since $T_{48}$ is a maximal group of symplectic transformations, we find $T_{48} = \tilde G_\mathrm{symp}$ and therefore
\[
 \{\mathrm{id}\} \to T_{48} \to \tilde G  \to C_6 \to \{\mathrm{id}\}.
\]
In analogy to the proof of Claim 2.1 in \cite{OZ168}, one can check that 288 is the maximal order of a finite group $H$ acting on a K3-surface with $T_{48} < H_\mathrm{symp}$. It follows that $\tilde G$ is maximal finite subgroup of $\mathrm{Aut}(X)$. For an arbitrary finite group $H$ acting on a K3-surface with 
$\{\mathrm{id}\} \to T_{48} \to H \to C_6 \to \{\mathrm{id}\}$,
there need however not exist an involution in $H$ centralizing $T_{48}$.
\end{remark}
%
%
%
%
\chapter{K3-surfaces with an antisymplectic involution centralizing $C_3 \ltimes C_7$}\label{chapterC3C7}
In this chapter it is illustrated that a classification of K3-surfaces with antisymplectic involution $\sigma$ can be carried out even even if the centralizer $G$ of $\sigma$ inside the group of symplectic transformations is relatively small, i.e., well below the bound 96 obtianed in Theorem \ref{roughclassi}, and not among the maximal groups of symplectic transformations. We consider the group $G = C_3 \ltimes C_7$, which is a subgroup of $L_2(7)$. The principles presented in Chapter \ref{chapterlarge} can be transferred to this group $G$ and yield a description of K3-surfaces with $G \times \langle \sigma \rangle$-symmetry. Using this, we deduce the classification K3-surfaces with an action of $L_2(7) \times C_2$ announced in Section \ref{mukaiL2(7)}. The results presented in this chapter have appeared in \cite{FKH1}.

To begin with, we present a family of K3-surfaces with $G \times \langle \sigma \rangle$-symmetry.
\begin{example}
We consider the action of $G$ on $\mathbb P_2$ given by one of its three-dimensional representations.
After a suitable change of coordinates, the action of the commutator subgroup $G'=C_7 < G$ is given by 
\[
[z_0:z_1:z_2] \mapsto [ \lambda z_0: \lambda^2 z_1: \lambda^4 z_2]
\]
for $\lambda = \mathrm{exp}(\frac{2\pi i}{7})$ and $C_3$ is generated by the permutation
\[
[z_0:z_1:z_2] \mapsto [z_2:z_0:z_1].
\]
The vector space of $G$-invariant homogeneous polynomials of degree six is the span of $P_1 = z_0^2 z_1^2 z_2^2$ and $P_2 = z_0^5 z_1 + z_2 ^5 z_0 + z_1^5 z_2$. 

The family $\mathbb{P}(V)$ of curves defined by polynomials in $V$
contains exactly four singular curves, namely the curve defined by
$z_0^2z_1^2z_2^2$ and those defined by 
$3z_0^2z_1^2z_2^2 -\zeta ^k(z_0^5z_1+z_2^5z_0+z_1^5z_2)$, where
$\zeta $ is a nontrivial third root of unity, $k=1,2,3$.
We let $\Sigma = \mathbb P(V) \backslash \{z_0^2z_1^2z_2^2=0\}$.

The double cover of $\mathbb P_2$ branched along a curve $C \in \Sigma$ is a K3-surface (singular K3-surface if $C$ is singular) with an action of $G \times C_2$ where $C_2$ acts nonsymplectically. It follows that $\Sigma$ parametrizes a family of K3-surface with $G \times C_2$-symmetry.
\end{example}
\begin{remark}
Let us consider the cyclic group $\Gamma$ of order
three generated by the transformation 
$[z_0:z_1:z_2]\mapsto [z_0:\zeta z_1: \zeta ^2z_2]$ and its induced action on the space $\Sigma$. One finds that the three irreducible singular $G$-invariant curves form a $\Gamma$-orbit. Furthermore, if two curves $C_1, C_2 \in  \Sigma$ are equivalent with respect to the action of $\Gamma$, then the corresponding K3-surfaces are equivariantly isomorphic (see Section \ref{EquiEqui} for a detailed discussion).
\end{remark}
\begin{remark}
The singular curve $C_\text{sing} \subset \mathbb P_2$ defined by
$3z_0^2z_1^2z_2^2 -(z_0^5z_1+z_2^5z_0+z_1^5z_2)$ has exactly seven
singular points $p_1, \dots p_7$ forming an $G$-orbit. Since they are in
general position (cf. Proposition \ref{blowdownseven}), the blow up of $\mathbb{P}_2$ in these points defines a Del Pezzo surface $Y_\text{Klein}$ of degree two with an action of $G$. It is seen to be the double cover of $\mathbb{P}_2$ branched along Klein's quartic curve\index{Klein's curve!}
$$
C_\text{Klein}:=\{z_0z_1^3+z_1z_2^3+z_2z_0^3=0\}.
$$
The proper transform $B$ of $C_\text{sing}$ in $Y_\text{Klein}$ is a
smooth $G$-invariant curve. It is a normalization of $C_\text{sing}$ and has genus three by the genus formula. The curve $B$ coincides with the
preimage of $C_{\text{Klein}}$ in $Y_\text{Klein}$.
The minimal resolution $\tilde X_\text{sing}$ of the singular surface $X_\text{sing}$ defined as the double cover of $\mathbb P_2$ branched along $C_\text{sing}$ is a K3-surface with an action of $G$. 
 By construction, it is the double cover of
$Y_\text{Klein}$ branched along $B$. In particular, $\tilde X
_\text{sing}$ is the degree four cyclic cover of $\mathbb P_2$ branched along $C_\text{Klein}$ and known as the Klein-Mukai-surface $X_{\text{KM}}$ (cf. Example \ref{L2(7)example}). 
\end{remark}
\begin{notation}
In the following, the notion of ``$G \times C_2$-symmetry'' abbreviates a symplectic action of $G$ centralized an antisymplectic action of $C_2$.
\end{notation}
In this chapter we will show that the space $\mathcal M = \Sigma / \Gamma$ parametrizes K3-surfaces with $G \times C_2$-symmetry up to equivariant equivalence. More precisely, we prove:
\begin{theorem}\label{mainthmc3c7}
The K3-surfaces with a symplectic action of $G = C_3 \ltimes C_7$ centralized by an antisymplectic involution $\sigma$ are parametrized by the space $\mathcal M = \Sigma / \Gamma$ of equivalence classes of sextic branch curves in $\mathbb P_2$. The Klein-Mukai-surface occurs as the minimal desingularization of the double cover branched along the unique singular curve in $\mathcal M$.  
\end{theorem}
Inside the family $\mathcal M$ one finds two K3-surfaces with a symplectic action of the larger group $L_2(7)$ centralized by an antisymplectic involution.
\begin{theorem}\label{L2(7) times invol}
There are exactly two K3-surfaces with an action of the group $L_2(7)$ centralized by an antisymplectic involution. These are the Klein-Mukai-surface $X_\mathrm{KM}$ and the double cover of $\mathbb P_2$ branched along the curve $\mathrm{Hess}(C_\text{Klein})=\{z_0^5z_1+z_2^5z_0+z_1^5z_2-5z_0^2 z_1^2 z_2^2=0\}$.
\end{theorem}
\section{Branch curves and Mori fibers}
Let $X$ be a K3 surface with an symplectic action of $G=C_3 \ltimes C_7$ centralized by the antisymplectic involution $\sigma$. We consider the quotient $\pi: X \to X/\sigma =Y$. Since the action of $G'$ has precisely three fixed points in $X$ and $\sigma$ acts on this point set, we know that $\mathrm{Fix}_X(\sigma)$ is not empty. It follows that $Y$ is a smooth rational surface with an effective action of the group $G$ to which we apply the equivariant minimal model program. The following lemma excludes the possibility that a $G$-minimal model is a conic bundle. The argument resembles that in the proof of Lemma \ref{conicbundle}.
\begin{lemma}\label{C3C7 conic bundle}
A $G$-minimal model of $Y$ is a Del Pezzo surface.
\end{lemma}
\begin{proof}
Assume the contrary and let $Y_\mathrm{min} \to \mathbb P_1$ be a $G$-equivariant conic bundle. Since $G$ has no effective action on the base, there must be a nontrivial normal subgroup acting trivially on the base. This subgroup must be $G'$. The action of $G'$ on the generic fiber has two fixed points and gives rise to a positive-dimensional $G'$-fixed point set in $Y_\mathrm{min}$ and $Y$. Since the action of $G'$ on $Y$ is induced by a symplectic action of $G'$ on $X$, this is a contradiction.
\end{proof}
\begin{remark}\label{P1xP1C3C7}
Since $G$ has no subgroup of index two,
 the above proof also shows that $Y_\mathrm{min} \not\cong \mathbb P_1 \times \mathbb P_1$. 
\end{remark}
In analogy to the procedure of the previous chapter we exclude rational and elliptic ramification curves and show that $\pi$ is branched along a single curve of genus greater than or equal to three.
\begin{proposition}
The set $\mathrm{Fix}_X(\sigma)$ consists of a single curve $C$ and $g(C) \geq 3$.
\end{proposition}
\begin{proof}
We let $\{x_1,x_2,x_3\} = \mathrm{Fix}_X(G')$. Since $G$ has no faithful two-dimensional representation, it has no fixed points in $X$ an therefore acts transitively on $\{x_1,x_2,x_3\}$. It follows that the central involution $\sigma$, which fixes at at least one point $x_i$, fixes all three points by invariance. Now $\{x_1,x_2,x_3\} \subset \mathrm{Fix}_X(\sigma)$ implies that $G'$ has precisely three fixed points in $Y$. Let $C_i$ denote the connected component of $\mathrm{Fix}_X(\sigma)$ containing $x_i$. Since $G$ acts on the set $\{C_1,C_2,C_3\}$, it follows that either $C_1=C_2=C_3$ or no two of them coincide. 

In the later case, it follows from Theorem \ref{FixSigma} that at least two curves $C_1,C_2$ are rational. The action of $G'$ on a rational curves $C_i$ has two fixed points. We therefore find at least five $G'$-fixed points in $X$ contradicting $|\mathrm{Fix}_X(G')|=3$.

It follows that all three points $x_1,x_2,x_3$ lie on one $G$-invariant connected component $C$ of $\mathrm{Fix}_X(\sigma)$. The action of $G$ on $C$ is effective and it follows that $C$ is not rational. 

If $g(C)=1$, then an effective action of $G$ on $C$ would force $G'$ to act by translations on $C$, in particular freely, a contradiction.

If $g(C)=2$, then $C$ is hyperelliptic. The quotient $C \to \mathbb P_1$ by the hyperellitic involution is $\mathrm{Aut}(C)$-equivariant and would induce an effective action of $G$ on $\mathbb P_1$, a contradiction.

It follows that $g(C) \geq 3$ and it remains to check that there are no rational ramification curves.

We let $n$ denote the total number of rational curves in $\mathrm{Fix}_X(\sigma)$. Since $G'$ acts freely on the complement of $C$ in $X$, it follows that the number $n$ must be a multiple of seven. Combining this observation with the bound $n \leq 9$ from Corollary \ref{atmostten} we conclude that $n$ is either 0 or 7.

We suppose $n =7$ and let $m$ denote the total number of Mori contractions of a reduction $Y \to Y_\mathrm{min}$. The Euler characeristic formula
\[
13 - g(C) = e(Y_\mathrm{min}) +m -n
\]
with $n=7$, $g(C) \geq 3$ and $e(Y_\mathrm{min})\geq 3$ implies $m \leq 14$. 

Let us first check that no Mori fiber $E$ coincides with a rational branch curve $B$. If this was the case, then all seven rational branch curves coincide with Mori fibers. Rational branch curves have self-intersection -4 by Corollary \ref{minusfour}. Before they may by contracted, they need to be transformed into (-1)-curves by earlier reduction steps. The remaining seven or less Mori contraction are not sufficient to achieve this transformation. It follows that each rational branch curve is mapped to a curve in $Y_\mathrm{min}$ and not to a point.

We now first consider the case $m=14$. The Euler characteristic formula implies $Y_\mathrm{min} \cong \mathbb P_2$ and $g(C)=3$. Using our study of Mori fibers and branch curves in Section \ref{branch curves mori fibers}, in particular Remark \ref{self-int of Mori-fibers} and Proposition \ref{at most two}, we see that no configuration of 14 Mori fibers is such that the images in $Y_\mathrm{min} \cong \mathbb P_2$ of any two rational branch curves have nonempty intersection. It follows that $m \leq 13$.

Let $R_1, \dots, R_7 \subset Y$ denote the rational branch curves. Each curve $R_i$ has self-intersection -4 and therefore has nontrivial intersection with at least one Mori fiber. Let $E_1$ be a Mori fiber meeting $R_1$, let $H \cong C_3$ be the stabilizer of $R_1$ in $G$ and let $I$ be the stabilizer of $E_1$ in $G$. Since $m \leq 13$ the group $I$ is nontrivial. If $I$ does not stabilize $R_1$, then $E_1$ meets the branch locus in at least three points. This is contrary to Proposition \ref{at most two}. It follows that $I=H$. If $E_1$ meets any other rational branch curve $R_2$, then it meets all curves in the $H$-orbit through $R_2$. Since $H$ acts freely on the set  $\{R_2, \dots, R_7\}$, it follows that $E_1$ meets three more branch curves. This is again contradictory to Proposition \ref{at most two}. 

Since $m \leq 13$ it follows that each rational branch curve meets exactly one Mori fiber. Their intersection can be one of the following three types:
\begin{enumerate}
\item
$E_i \cap R_i = \{p_1,p_2\}$ or
\item
$E_i \cap R_i = \{p\}$ and $ (E_i, R_i)_p =2$ or
\item
$E_i \cap R_i = \{p\}$ and $ (E_i, R_i)_p =1$.
\end{enumerate}
In all three cases the contraction of $E_i$ alone does not transform the curve $R_i$ into a curve on a Del Pezzo surface. So further reduction steps are needed and require the existence of Mori fibers $F_i$ disjoint from $\bigcup R_i$. Each $F_i$ is a (-2)-curve meeting $\bigcup E_i$ transversally in one point and the total number of Mori fibers exceeds our bound 13.

This contradiction yields $n=0$ and the proof of the proposition is completed.
\end{proof}
\section{Classification of the quotient surface $Y$}
We now turn to a classification of the quotient surface $Y$. 
\begin {proposition}
The surface $Y$ is either $G$-minimal or the blow up of $\mathbb P_2$ in seven singularities of an irreducible $G$-invariant sextic..
\end {proposition}
\begin {proof}
Since $n=0$, the Euler characteristic formula yields 
$m \leq 7$. The fact that
$G$ acts on the set of Mori fibers implies that 
$m \in \{ 0,3,6, 7\}$.  If $m \in \{3, 6\}$,
then $G'$ stabilizes every Mori fiber, and
consequently it has more then three fixed points, a contradiction. Thus we must only 
consider the case $m =7$.

In this case the set of Mori fibers is a $G$-orbit and it follows that every
Mori fiber has self-inter\-section -1 and therefore has
nonempty intersection with $\pi(C)$ by Remark \ref{self-int of Mori-fibers}.

As before, the Euler characteristic formula
implies that $g(C)=3$ and $Y_\mathrm{min}=\mathbb P_2$ and 
adjunction in $X$ shows that $(\pi(C))^2=8$ in $Y$. The fact that
$\pi(C)$ has nonempty intersection with seven different Mori fibers implies
that its image $D$ in $Y_\mathrm{min}$ has self-intersection either $15 = 8 +7$ or $36 = 8 + 4 \cdot 7$. Since the first is impossible it follows that $E \cdot \pi(C)=2$ for all Mori fibers $E$ and the $G$-invariant irreducible sextic $D$ has seven singular points corresponding to the images of $E$ in $\mathbb P_2$.
\end {proof}
\begin{corollary}\label{sing exam}
If $Y$ is not $G$-minimal, then $X$ is the minimal desingularization of a double cover of $\mathbb{P}_2$ branched along an irreducible $G$-invariant sextic with seven singular points.
\end{corollary}
We conclude this section with a classification of possible
$G$-minimal models of $Y$.
\begin {proposition}
The surface $Y_\mathrm{min}$ is either a Del Pezzo surface of degree two
or $\mathbb P_2$.
\end {proposition}
\begin {proof}
The case $Y_\mathrm{min}=\mathbb P_1\times \mathbb P_1$ is excluded by Example \ref{DelPezzoC3C7} and also by Remark \ref{P1xP1C3C7}.

Thus $Y_\mathrm{min}=Y_d$ 
is a Del Pezzo surface of degree $d=1,\ldots ,9$ which is a blowup
of $\mathbb P_2$ in $9-d$ points.

If $Y_\mathrm{min}=Y_1$ the anticanonical map has exactly one base point.  
This point has to be $G$-fixed and since $G$ has no faithful 
two-dimensional representations, this case does not occur.

It remains to eliminate $d=8,\ldots ,3$. In these
cases the sets $\mathcal S$ of (-1)-curves consist
of 1, 2, 6, 10, 16 or 27 elements, respectively (cf. Table \ref{minus one curves}).  
The $G$-orbits in $\mathcal S$
consist of $1$, $3$, $7$ or $21$ curves and  
there must be orbits of length three or one.
If $G$ stablizes a curve
in $\mathcal S$, then its contraction gives rise to a two-dimensional
representation of $G$ which does not exist. If $G$ has an
orbit consisting of three curves, then $G'$ stabilizes each of the 
curves in this orbit. Thus $G'$ has at least six fixed points in
$Y_\mathrm{min}$ and in $Y$. This contradicts the fact that $| \mathrm{Fix}_Y(G')|=3$. 
\end {proof}
\section{Fine classification - Computation of invariants}
We have reduced the classification of K3-surfaces with $G \times C_2$-symmetry to the study of equivariant double covers of rational surfaces $Y$ branched along a single invariant curve of genus $g \geq 3$. Here $Y$ is either $\mathbb P_2$, the blow-up of $\mathbb P_2$ in seven singular points of an irreducible $G$-invariant sextic, or a Del Pezzo surface of degree two. 
\subsection{The case $Y = Y_\mathrm{min} = \mathbb P_2$}
An effective action of $G$ on $\mathbb P_2$ is given by an injective homomorphisms $G \to \mathrm {PSL}_3(\mathbb C)$. There are two central degree three extension of $G$, the trivial extension and $C_9 \ltimes C_7$. A study of their three-dimensional representation reveals that in both cases the action of $G$ on $\mathbb P_2$ is given by an irreducible representation $G \hookrightarrow \mathrm {SL}_3(\mathbb C)$. There are two isomorphism classes of irreducible 3-dimensional representations.
Since these differ by a group automorphism
and the corresponding actions on $\mathbb{P}_2$ are therefore equivalent, 
we may assume that in appropriately chosen coordinates a generator of $G'$ acts by
\begin{equation}\label{c7action}
[z_0:z_1:z_2]\mapsto [\lambda z_0,\lambda ^2,z_1,\lambda ^4z_2],
\end{equation}
 where $\lambda =\mathrm{exp}{\frac{2\pi i}{7}}$ and a generator of
$C_3$ acts by the cyclic permutation 
 $\tau $ which is defined by 
\begin{equation}\label{tauaction}
[z_0:z_1:z_2]\mapsto [z_2:z_0:z_1].
\end{equation}

A homogeneous polynomial defining an invariant curve 
must be a $G$-semi-invariant with $G'$ acting with eigenvalue one.
The $G'$-invariant monomials of degree six are 
$$ 
\mathbb C[z_0,z_1,z_2]_{(6)}^{G'}=
\mathrm {Span}\{z_0^2z_1^2z_2^2, z_0^5z_1,z_2^5z_0,z_1^5z_2\}\,.
$$
Letting $P_1=z_0^2z_1^2z_2^2 $ and $P_2=z_0^5z_1+z_2^5z_0+z_1^5z_2$,
it follows that 
$$
\mathbb C[z_0,z_1,z_2]_{(6)}^{G}=
\mathrm {Span}\{P_1,P_2\}=:V\,.
$$
There are two $G$-semi-invariants which are not invariant, namely 
$z_0^5z_1+\zeta z_2^5z_0+\zeta ^2z_1^5z_2$ for $\zeta ^3=1$
but $\zeta \not =1$. By direct computation one checks that
the curves defined by these polynomials are smooth and that in both
cases all $\tau $-fixed points in $\mathbb P_2$ lie on them.
Thus, $\tau $ has only three fixed points on the K3-surface
$X$ obtained as a double cover and therefore does not act symplectically (cf. Table \ref{fix points symplectic}). Consequently, 
$G$ does not lift to an action by symplectic transformations
on the K3-surfaces defined
by these two curves. Hence it is enough to consider  
ramified covers $X\to Y=\mathbb P_2$, where the branch curves are defined by invariant polynomials $f \in V$. 

We wish to determine which polynomials 
$P_{\alpha,\beta}=\alpha P_1+\beta P_2$ define singular curves.
Since $\mathrm {Fix}(\tau) = \{ [1:\zeta :\zeta ^2] \, | \, \zeta ^3 =1 \}$, the
curves which contain $\tau $-fixed points are
defined by condition $\alpha +3\zeta \beta =0$.
Let $C_{P_1} = \{ P_1 =0\} $ and let $C_\zeta $ be
the curve defined by $P_{\alpha,\beta}$ for $\alpha +3\zeta \beta =0$.  A direct computation shows that $C_\zeta $ is singular at the point $[1:\zeta :\zeta ^2]$.  We let $\Sigma_\mathrm{reg} $ be
the complement of this set of four curves,  $\Sigma_\mathrm{reg} =   \mathbb P(V) \backslash \{C_{P_1}; \, C_\zeta \, | \, \zeta ^3 =1 \}$. 
\begin{lemma}
A curve $C\in \mathbb P(V)$ is smooth if and only if $C \in \Sigma_\mathrm{reg} $.
\end{lemma}
\begin{proof}
Let $C \in \Sigma_\mathrm{reg} $. Since $\tau $ has no fixed points in $C$ by definition and every subgroup of
order three in $G$ is conjugate to $\langle \tau \rangle $, it follows that any $G$-orbit $G.p$ through a point $p\in C$ has length three or 21. 

The only
subgroup of order seven in $G$ is the commutator group
$G'$.  So the $G$-orbits of length three 
are the orbits of the $G'$-fixed points $[1:0:0], [0:1:0],[0:0:1]$. One checks by direct computation that every $C\in \Sigma_\mathrm{reg} $ is smooth
at these three points.  

An irreducible curve of degree six has at most ten singular points by the genus formula.
Suppose that $C$ is singular at some point $q$. Then it
is singular at each of the 21 points in $G.q$ and $C$ must be reducible.
Considering
the $G$-action on the space of irreducible components of $C$ yields a contradiction
and it follows that $C$ is smooth.
\end{proof}
For any curve $C \in \Sigma_\mathrm{reg}$ the double cover of $\mathbb P_2$ branched along $C$ is a K3-surface $X_C$ with an action of a degree two central extension of $G$. By the following lemma, this action is always of the desired type.
\begin{lemma}\label{G acts sympl}
For every $C \in \Sigma_\mathrm{reg}$ the K3-surface $X_C$ carries an action of the group 
$G \times \langle \sigma \rangle $.
The group $G$ acts by symplectic transformations on $X_C$ and $\sigma$ denotes the covering involution.
\end{lemma}
\begin{proof}
It follows from the group structure of $G$
that the central degree two extension of $G$ acting on $X_C$ splits as $G \times C_2$. The factor $C_2$ is by construction generated by the covering involution $\sigma$. It remains to check that $G$ acts symplectically. As the commutator subgroup $G'$ acts symplectically it is sufficient to check whether $\tau$ lifts to a symplectic automorphism. Consider the $\tau$-fixed point $p=[1:1:1]$ and check that the linearization of $\tau$ at $p$ is in $\mathrm{SL}(2, \mathbb C)$. Since $p$ is not contained in $C$, it follows that the linearization of $\tau$ at a corresponding fixed point in $X_C$ is also in $\mathrm{SL}(2, \mathbb C)$. Consequently, the group $G$ acts by symplectic transformations on $X_C$.
\end{proof}
\subsection{Equivariant equivalence} \label{EquiEqui}
We wish to describe the space of K3-surfaces with $G \times C_2$-symmetry modulo equivariant equivalence.
For this, we study the family of K3-surfaces parametrized by the family of branch curves $\Sigma_\mathrm{reg}$. Consider the cyclic group $\Gamma$ of order three in $\mathrm{PGL}(3, \mathbb C)$ generated by 
\[
[z_0:z_1:z_2]\mapsto [z_0:\zeta z_1: \zeta ^2z_2] 
\]
for $\zeta = \mathrm{exp}(\frac{2 \pi i}{3})$. The group $\Gamma$ acts on $\Sigma_\mathrm{reg}$ and by the following proposition the induced equivalence relation is precisely equivariant equivalence  formulated in Definition \ref{equivariantequivalence}.
\begin{proposition}
 Two K3-surfaces $X_{C_1}$ and $X_{C_2}$ for $C_1,C_2 \in \Sigma_\mathrm{reg}$ are equivariantly equivalent if and only if $C_1 = \gamma C_2$ for some $\gamma \in \Gamma$, i.e., the quotient $\Sigma_\mathrm{reg}/ \Gamma$ parametrizes equivariant equivalence classes of K3-surfaces $X_C$ for $C \in \Sigma_\mathrm{reg}$. 
\end{proposition}
\begin{proof}
If two K3-surfaces $X_{C_1}$ and $X_{C_2}$ for $C_1,C_2 \in \Sigma_\mathrm{reg}$ are equivariantly equivalent, then the isomorphism $X_{C_1} \to X_{C_2}$ induces an automorphism of $\mathbb P_2$ mapping $C_1$ to $C_2$.

Let $C\in \Sigma_\mathrm{reg}$ and for $T\in \mathrm {SL}_3(\mathbb C)$
 assume that $T(C)\in \Sigma_\mathrm{reg} $. We consider the group span
 $S$ of $TGT^{-1}$ and $G$. By Lemma \ref{G acts sympl}, the group $G$ acts by symplectic transformations on $X_C$ and $X_{T(C)}$. We argue precisely as in the proof of this lemma to see that $TGT^{-1}$ also acts symplectically on the K3-surface $X_{T(C)}$. It follows that $S$ is acting as a group of symplectic transformations on this K3-surface. 

If $S=G$, then $T$ normalizes $G$. The normalizer $N$ of $G$ in
$\mathrm {PGL}_3(\mathbb C)$
is the product
$\Gamma \times G$ and
 it follows that $gT$ is contained
 in $\Gamma $ for some $g \in G$ and $T(C) = gT(C) = \gamma C$. 

Note 
that $L_2(7)$ is the only group in Mukai's
list which contains $G$.  Therefore,
$S$ is a subgroup of $L_2(7)$. The group $G$ is a maximal subgroup of $L_2(7)$ and if $S\not=G$, then it follows that
$S=L_2(7)$. Any two subgroups of order 21 in $L_2(7)$ are conjugate. This implies the existence of $s \in S= L_2(7)$ such that $sTGT^{-1}s^{-1} = G$. Now $sT \in N = \Gamma \times G$ can be written as $sT = \gamma g$ for $(\gamma, g) \in  \Gamma \times G$. By assumption, $s$ stabilizes $T(C)$ and $T(C) = sT(C) = \gamma g (C) = \gamma C$.
This completes the proof of the proposition.
\end{proof}
\subsection{The case $Y \neq Y_\mathrm{min}$}
Let us now consider the three singular irreducible curves in our family $\mathbb P (V)$. They are identified by the action of $\Gamma$. Using Corollary \ref{sing exam} we see that if $Y = X/\sigma$ is not $G$-minimal, then, up to equivariant equivalence, the K3-surface $X$ is the minimal desingularization of 
the double cover of $\mathbb P_2$ branched along $C_{\zeta=1} = C_\mathrm{sing}$ and $Y$ is the blow-up of $\mathbb P_2$ in the seven singular points of $C_\mathrm{sing}$. These points are the $G'$-orbit of $[1:1:1]$. In the following propostion we prove that these are in general position and therefore $Y$ is a Del Pezzo surface.
\begin{proposition}\label{blowdownseven}
If $Y$ is not minimal, then it is the Del Pezzo surface of degree two which
arises by blowing up the seven singular points $p_1,\ldots ,p_7$ on the curve
$C_\mathrm {sing}$ in $\mathbb P_2$. The corresponding map $Y \to
\mathbb P_2$ is $G$-equivariant and therefore a Mori reduction of
$Y$. 
\end{proposition}
\begin{proof}
We show that the points $\{p_1, \dots, p_7\} = G'.[1:1:1]$ 
are in general position, i.e., no three lie on one line and 
no six lie on one conic.
It follows from direct computation that no three points in $G'.[1:1:1]$ 
lie on one line. 
If $p_1,\dots p_6$ lie on a conic $Q$, then $g.p_1, \dots , g.p_6$ lie on 
$g.Q$ for every  $g \in G$. Since $\{p_1, \dots, p_7\}$ is a
$G$-invariant set, the conics $Q$ and $g.Q$ intersect in at 
least five points and therefore coincide. It follows that 
$Q$ is an invariant conic meeting $C_\mathrm{sing}$ at its 
seven singularities and $(Q,C_\mathrm{sing}) \geq 14$ implies 
$Q \subset C_\mathrm{sing}$, a contradiction.
\end{proof}
\section{Klein's quartic and the Klein-Mukai surface}\label{KMsurface}
In this section we show that the Del Pezzo surface discussed in Proposition \ref{blowdownseven} above can be realized as the double cover of $\mathbb P_2$ branched along Klein's quartic curve. 
\begin{proposition}\label{DelPezzoYKlein}
A Del Pezzo surface of degree two with an action of $G$ is equivariantly isomorphic to the double cover $Y_\mathrm{Klein}$ of $\mathbb P_2$ branched along Klein's quartic curve.
\end{proposition}
\begin{proof}
Recall that the anticanonical map of a Del Pezzo surface $Y$ of degree two defines a 2:1 map to $\mathbb P_2$. This map is branched along a smooth curve of degree four and equivariant with respect to $\mathrm{Aut}(Y)$. We obtain an action of $G$ on $\mathbb P_2$ stabilizing a smooth quartic. As before, we may choose coordinates such that $G$ is acting as in equations \eqref{c7action} and \eqref{tauaction}. Then
$$
\mathbb C[z_0:z_1:z_2]_{(4)}^{G'}
=\mathrm {Span}\{z_0^3z_2, z_1^3z_0,z_2^3z_1\}\,.
$$ 
is a direct sum of $G$-eigenspaces. The eigenspace
of the eigenvalue $\zeta $ is spanned by the polynomial
$Q_\zeta :=z_0^3z_2+\zeta z_2^3z_1+\zeta ^2z_1^3z_0$ with $\zeta $
being a third root of unity.

In order to take into account equivariant equivalence
we consider the cyclic group 
$\Gamma \subset \mathrm {SL}_3(\mathbb C)$ which 
is generated by the transformation 
$\gamma$, $[z_0:z_1:z_2]\mapsto [z_0:\zeta z_1:\zeta^2z_2]$.
The induced action on $\mathbb C[z_0:z_1:z_2]_{(4)}^{G'}$ is transitive on the $G$-eigenspaces
spanned by the $Q_\zeta $.  Consequently, up to equivariant
equivalence, we may assume that
$Y\to \mathbb P_2$ is branched along Klein's curve $C_\text{Klein}$\index{Klein's curve!}
which is defined by $Q_1$. 
\end{proof}
\begin{corollary}
A Del Pezzo surface of degree two with an action of $G$ is never $G$-minimal. Its Mori reduction $Y_\mathrm{Klein} \to \mathbb P_2$ is precisely the map discussed in Proposition \ref{blowdownseven}.
\end{corollary}
We summarize our observartions in the following proposition.
\begin {proposition}
If $X$ is a K3-surface  with a symplectic $G$-action centralized
 by an antisymplectic involution $\sigma $,  then $Y_{min}=\mathbb P_2$.
In all but one case $X/\sigma =Y=Y_{min}$. In the exceptional case
$Y=Y_\mathrm{Klein}$, the Mori reduction $Y\to Y_{min}$ is the contraction of seven (-1)-curves to the singular 
points of $C_\mathrm{sing}$ and the branch set $B$ of $X\to Y$
is the proper transform of $C_\mathrm {Klein}$ in $Y$.
\end {proposition}
\begin {proof}
It remains to prove that $B$
is the proper transform of $C_\mathrm {Klein}$ in $Y$.
Suppose that the
branch curve of $X\to Y$ is some other curve 
$\widetilde B$ linearly equivalent to $-2K_Y$. Let
$I:=\widetilde B\cap B$ and note that $\vert I\vert \le B \cdot \widetilde B = 4 K_Y ^2 = 8$.
Since $G$ has no fixed points in $B$, it follows that
$\vert I\vert =3$ and that $I$ is a $G$-orbit. 
Thus the
intersection multiplicities at the three points in $\widetilde B\cap B$ are
the same. Since 3 does not divide 8, this is a contradiction.
\end {proof}
In order to complete the proof of Theorem \ref{mainthmc3c7} it remains to show that the action of $G$ on $Y_\mathrm{Klein}$ lifts to a group of symplectic transformation on the K3-surface $X=X_{KM}$ defined as a double cover of $Y_\mathrm{Klein}$ branched along the proper transform of $C_\mathrm{sing}$.

Since $G$ stabilizes $C_\text{Klein}$ and does not admit nontrivial central extensions of degree two, it lifts to a subgroup of $\mathrm{Aut}(Y_\text{Klein})$ and subsequently to a subgroup of $\mathrm{Aut}(X)$. 

The covering involution  
$Y_\text{Klein}\to \mathbb P_2$,
lifts to a holomorphic transformation of $X$ where we also find the
involution defining $X \to Y_\text{Klein}$. These two transformations
generate a group $F$
of order four. The elements of $F$ all have a positive-dimensional
fixed point set. It follows that $F$ acts solely by nonsymplectic transformations and is therefore isomorphic to $C_4$.
The full preimage of $G$ in $\mathrm {Aut}(X)$ splits 
as $G\times C_4$.  

Since the commutator group $G'$ 
automatically acts by symplectic transformations, we must
only check that the lift of the cyclic permutation $\tau $,
$[z_0:z_1:z_2]\mapsto [z_2:z_0:z_1]$, acts symplectically.
As above, this follows from a linearization argument at a $\tau$-fixed point not in $C_\mathrm{Klein}$. 

In conclusion, up to equivalence there
is a unique action of $G$ by symplectic transformations 
on the K3-surface $X_{KM}$. It is centralized by a cyclic 
group of order four which acts faithfully on the symplectic
form. 

The Klein-Mukai-surface is the only surface with $G \times C_2$-symmetry for which $Y \not\cong \mathbb P_2$. As in the introduction of this chapter, we define $\Sigma$ as the complement of $C_{P_1}$ in $\mathbb P(V)$. Then $\Sigma = \Sigma_\mathrm{reg} \cup \{C_\zeta \, | \, \zeta^3 =1\}$. Using this notation the space
\[
\mathcal M = \Sigma / \Gamma
\]
parametrizes the space of K3-surfaces with $G \times C_2$-symmetry up to equivariant equivalence.
This completes the proof of Theorem \ref{mainthmc3c7}.
\section{The group $L_2(7)$ centralized by an antisymplectic involution}\label{168}
We consider the simple group of order 168. This group is $\mathrm{PSL}(2, \mathbb F_7)$ and usually denoted by $L_2(7)$. It contains our group $G = C_3 \ltimes C_7$ as a subgroup. 
Since $L_2(7)$ is a simple group, if it acts on a K3-surface, it automatically acts by symplectic transformations. 

We wish to prove Theorem \ref{L2(7) times invol} stating that there are exactly two K3-surfaces with an action of the group $L_2(7)$ centralized by an antisymplectic involution. These are the Klein-Mukai-surface $X_\mathrm{KM}$\index{Klein-Mukai surface!} and the double cover of $\mathbb P_2$ branched along the curve $\mathrm{Hess}(C_\text{Klein})=\{z_0^5z_1+z_2^5z_0+z_1^5z_2-5z_0^2 z_1^2 z_2^2=0\}$.

We have to check which elements of $\mathcal M$ have the symmetry of the larger group. The Klein-Mukai-surface is known to have $L_2(7) \times C_4$-symmetry (cf. Example \ref{L2(7)example}). 
If $X \neq X_\mathrm{KM}$ has $L_2(7)$-symmetry, then it follows from the considerations of the previous sections that $X$ is an $L_2(7)$-equivariant double cover of $\mathbb P_2$ branched along a smooth $L_2(7)$-invariant sextic curve. I.e., it remains to identify the surfaces with $L_2(7)$-symmetry in the family parametrized by $\Sigma_\mathrm{reg} / \Gamma$. 
\begin{lemma}\label{L2(7)onP2}
The action of $L_2(7)$ on $\mathbb P_2$ is necessarily given by a three-dimensional represention.
\end{lemma}
\begin{proof}
The lemma follows from the fact that the group $L_2(7)$ does not admit nontrivial degree three central extensions. This can be derived from the cohomology group $H^2(L_2(7), \mathbb C^*) \cong C_2$ known as the Schur Multiplier. 
\end{proof}
There are two isomorphism classes of three-dimensional representations and these differ by an outer automorphism. We may therefore consider the particular representation given in Example \ref{L2(7)example}. One checks that the curve $\mathrm{Hess}(C_\mathrm{Klein})$ is $L_2(7)$-invariant.
The maximal possible isotropy group is $C_7$ and each $L_2(7)$-orbit in  $\mathrm{Hess}(C_\mathrm{Klein})$ consists of at least 21 elements. If there was another $L_2(7)$-invariant curve $C$ in $\Sigma_\mathrm{reg}$, then the invariant set $ C \cap \mathrm{Hess}(C_\mathrm{Klein})$ consists of at most 36 points. This is a contradiction and it follows that $\mathrm{Hess}(C_\mathrm{Klein})$ is the only $L_2(7)$-invariant curve in $\Sigma_\mathrm{reg}$.

It remains to check that $L_2(7)$ lifts to a subgroup of $\mathrm{Aut}(X_{\mathrm{Hess}(C_\mathrm{Klein})})$: On $X_{\mathrm{Hess}(C_\mathrm{Klein})}$ we find an action of a central degree two extension $E$ of $L_2(7)$. Since $E \neq E_\mathrm{symp}$ and $L_2(7)$ is simple, the subgroup of symplectic transformations inside $E$ must be isomorphic to $L_2(7)$.

It follows that $X_\mathrm{KM}$ and the double cover of $\mathbb P_2$ branched along $\mathrm{Hess}(C_\mathrm{Klein})$ are the only examples of K3-surfaces with $L_2(7) \times C_2$ symmetry. This completes the proof of Theorem \ref{L2(7) times invol}.
\begin{remark}
 If we consider the quotient $Y_\mathrm{Klein}$ of $X_\mathrm{KM}$ by the antisymplectic involution $\sigma \in C_4$, this surface was seen not to be minimal with respect to the action of $C_3 \ltimes C_7$. It is however $L_2(7)$-minimal as we cannot find a equivariant contraction morphism blowing down an orbit of disjoint (-1)-curves in $Y_\mathrm{Klein}$ .
Such an orbit would have to consists of seven Mori fibers. The only subgroup of index seven is $S_4$. A Mori fiber of self-intersection (-1) does however not admit an action of the group $S_4$ (cf. Proof of Theorem \ref{roughclassi}).
\end{remark}
%
%
%
%
\chapter{The simple group of order 168}\label{chapter non exist}
In this chapter we consider finite groups containing $L_2(7)$, the simple group of order 168, and their actions on K3-surfaces. Based on our considerations about $L_2(7) \times C_2$-actions on K3-surfaces in Section \ref{168} we derive a classification result (Theorem \ref{improve OZ}). This gives a refinement of a lattice-theoretic result due to Oguiso and Zhang \cite{OZ168}. The main part of this chapter is dedicated to proving the non-existence of K3-surfaces with an action of the group $L_2(7) \times C_3$ (Theorem \ref{nonexist}) using equivariant Mori reduction.
\section{Finite groups containing $L_2(7)$}
If $H$ is a finite group acting on a K3-surface and $L_2(7) \lneqq H$, then it follows from Mukai's theorem and the fact that $L_2(7)$ is simple, that $H$ fits into the short exact sequence
\[
1 \to L_2(7) = H_\mathrm{symp} \to H \to C_m \to 1
\]
for some $m \in \mathbb N$. As it is noted by Oguiso and Zhang, Claim 2.1 in \cite{OZ168}, it follows from Proposition 3.4 in \cite{mukai} that $m \in \{1,2,3,4,6\}$. 

The action of $H$ on $L_2(7)$ by conjugation defines a homomorphism 
$H \to  \mathrm{Aut}(L_2(7))$.
 Factorizing by the group of inner automorphism of $L_2(7)$ we obtain a homomorphism
\[
 C_m \cong H/L_2(7) \to \mathrm{Out}(L_2(7)) \cong C_2.
\]
If $H$ is not the nontrivial semidirect product $L_2(7) \rtimes C_2$, this homomorphism has a nontrivial kernel. In particular, we find a cyclic group $C_k < C_m$ centralizing $L_2(7)$. If $k$ is even, we may apply our results on K3-surfaces with $L_2(7) \times C_2$-symmetry from the previous chapter. 

If $m = 3,6$, then $k=3$ or $k=6$. These cases may be excluded as is shown in \cite{OZ168}, Added in proof, Proposition 1. An independent proof of this fact, i.e., the non-existence of K3-surfaces with $L_2(7) \times C_3$ symmetry, using equivariant Mori theory, in particular the classification of $L_2(7)$-minimal models, is given below (Theorem \ref{nonexist}).

We summarize our observations about K3-surfaces with $L_2(7)$-symmetry in the following theorem, which improves the classification result due to Oguiso and Zhang.
\begin{theorem}\label{improve OZ}
Let $H$ be finite group acting on a K3-surface $X$ with $L_2(7) \lneqq H$. Then 
\begin{samepage}
\begin{itemize}
\item
$|H/L_2(7)| \in \{2,4\}$. 
\item
If  $|H/L_2(7)|= 4$, then $H = L_2(7) \times C_4$ and $X \cong X_\mathrm{KM}$. 
\item
If $|H/L_2(7)|= 2$ and $H = L_2(7) \times C_2$, then either $X \cong X_\mathrm{KM}$ or $X \cong X_{\mathrm{Hess}(C_\mathrm{Klein})}$
\end{itemize}
\end{samepage}
\end{theorem}
The first statement follows from the non-existence of K3-surfaces with $L_2(7) \times C_3$-symmetry (Theorem \ref{nonexist} below) and the third statement follows from  Theorem \ref{L2(7) times invol}. The remaining part ist covered in the following lemma (cf. Main Theorem in \cite{OZ168}).
\begin{lemma}\label{OZresult}
If $X$ is a K3-surface with an action of a finite group containing the $L_2(7)$ as a subgroup of index four, then $X$ is the Klein-Mukai-surface.
\end{lemma}
\begin{proof}
We let $X$ be a K3-surface and $H$ be a finite subgroup of $\mathrm{Aut}(X)$ with $L_2(7) < H$ and $|H/L_2(7)|=4$. 

Since $L_2(7)$ is simple and a maximal group of symplectic transformations, 
it coincides with the group of symplectic transformations in $H$. 
In particular, $H / L_2(7) = C_4$ and a group 
$\langle \sigma \rangle$ of
order two is contained in the kernel of the homomorphism $H\to \mathrm {Aut}(L_2(7))$.
It follows that we are in the setting of Theorem \ref{L2(7) times invol} where
$\Lambda :=H/\langle \sigma \rangle$ acts on $Y=X/\sigma $.  If
$X\not =X_{KM}$, then $Y= \mathbb P_2$. This possibility needs to be eliminated.

Let $\tau $ be any element of $\Lambda $ which
is not in $L_2(7)$ and let $\Gamma = C_3 \ltimes C_7 < L_2(7)$. 
Since any two subgroups of order 21 in $L_2(7)$ are conjugate by an element of $L_2(7)$, it follows that there exists
$h\in L_2(7)$ with $(h\tau )\Gamma(h\tau)^{-1}=\Gamma$.  Thus, the normalizer
$N(\Gamma)$ of $\Gamma$ in $\Lambda $ is a group of order 42 which also normalizes the commutator subgroup
$\Gamma'$ and therefore stabilizes its set $F$ of fixed points.

Using coordinates $[z_0:z_1:z_2]$ of $\mathbb{P}_2$ 
as in Theorem \ref{L2(7) times invol} one
checks by direct computation that the only transformations in $\mathrm {Stab}(F)$ which stabilize 
the branch curve $\mathrm {Hess}(C_{\mathrm {Klein}})$
are those in $\Gamma$ itself.  This contradiction shows that
$Y\not =\mathbb P_2$ and therefore $X=X_{KM}$. 
\end{proof}
\section{Non-existence of K3-surfaces with an action of $L_2(7) \times C_3$}
The method of equivariant Mori reduction can be applied to obtain both classification and non-existence results.
In the following, we exemplify a general approach to prove non-existence of K3-surfaces with specified symmetry by considering the group $L_2(7) \times C_3$ and give an independent proof of the following observation of Oguiso and Zhang \cite {OZ168}:
\begin{theorem}\label{nonexist}
There does not exist a K3-surface with an action of $L_2(7) \times C_3$.
\end{theorem}
The remainder of this chapter is dedicated to the proof of this theorem.
\subsection{Global structure}
Let $G \cong L_2(7)$, let $D \cong C_3$, and assume there exists a K3-surface $X$ with a holomorphic action of $G \times
D$. Since $G$ is a simple group and a maximal group of symplectic transformations on a K3-surface, it follows that $G$ acts symplectically whereas the action of $D$ is nonsymplectic. We obtain the following commuting diagram.
\[
\begin{xymatrix}{
X \ar[d]^{\pi} & \hat{X} \ar[d]^{\hat{\pi}} \ar[l]_{b_X}\\
X/D=Y & \hat{Y} \ar[l]_{b_Y}\ar[d]^{M_{\mathrm{red}}}\\
& \hat{Y}_\mathrm{min}=Z
}
\end{xymatrix}
\] 
Here $b_X$ is the blow-up of the isolated $D$-fixed points in $X$.
The singularities of $X/D$ correspond to isolated $D$-fixed points.
Since the linearization of the $D$-action at an isolated fixed point is locally of the form $(z,w) \mapsto (\chi z, \chi w)$ for some nontrivial character $\chi: D \to \mathbb C^*$, each singularity of $X/D$ is resolved by a single blow-up. We let 
$b_Y$ denote the simultanious blow-up of all singularities of $Y$. We fix a $G$-Mori reduction $M_\mathrm{red}: \hat Y \to \hat{Y}_\mathrm{min}=Z$. All maps in
the diagram are $G$-equivariant.  By Theorem \ref{K3quotnonsymp}, the surface $\hat{Y}$ is rational. As conic bundles do not admit an action of $G$ (cf. Lemma \ref{C3C7 conic bundle}), we know that $\hat{Y}_\mathrm{min}$ is
a Del Pezzo surface . The following lemma specifies $Z$.
\begin{lemma}
The Del Pezzo surface $Z$ is either $\mathbb P_2$
or a surface obtained from $\mathbb P_2$ by blowing up 7 points in general
position. In the later case, $Z$ is a $G$-equivariant double cover of $\mathbb P_2$ branched along
Klein's quartic curve. The action of $G$ on $\mathbb P_2$ is given by a three-dimensional representation. 
\end{lemma}
\begin{proof}
The first part of the lemma follows from our observations in Example \ref{DelPezzoL2(7)}, the last part has been discussed in Lemma \ref{L2(7)onP2}. 
If $Z$ is a Del Pezzo surface of degree two, then the anticanonical map realizes it as an equivariant double cover of $\mathbb P_2$ branched along a smooth quartic curve $C$. We choose coordinates on $\mathbb P_2$ such that the action of $G$ is given by the representation $\rho$ of Example \ref{L2(7)example} (or its dual represenation $\rho^*$) and have already seen that Klein's quartic curve
\[
C_\text{Klein}= \{x_1x_2^3 + x_2x_3^3 + x_3x_1^3=0\} \subset \mathbb P_2
\]
is $G$-invariant. If $C \neq C_\text{Klein}$, then $C \cap C_\text{Klein}$ is a $G$-invariant subset of $\mathbb P_2$. Since the maximal cyclic subgroup of $G$ is of order seven, it follows that a $G$-orbit $G.p$ for a point $p \in C \cap C_\text{Klein}$ consists of at least 24 elements. Since $C \cap C_\text{Klein}$ however consists of at most 16 points,
this is a contradiction. Therefore, $C= C_\text{Klein}$ and the lemma follows.
\end{proof}
\subsection*{$D$-fixed points}
The map $\pi$ is in general ramified both at points and along curves.
Let $x$ be an isolated
$D$-fixed point in $X$. As was noted above, the isotropy representation of the nonsymplectic $D$-action at $x$ in local coordinates $(z,w)$ is given by
$(z,w) \mapsto (\chi z, \chi w)$ for some nontrivial character $\chi: D \to \mathbb C^*$. The action of $D$ on the rational
curve $\hat E$ obtained by blowing up $x$ is trivial and therefore $\hat E$ is contained
in the ramification set $\mathrm{Fix}_{\hat{X}}(D)$. Let $\{\hat{E_i}\}$ denote the set of (-1)-curves in $\hat{X}$ obtained from blowing up isolated $D$-fixed points in $X$ and define $E_i = \hat{\pi}(\hat{E_i})$.

If $C$ is a curve of $D$-fixed points in $X$, it follows that $\hat \pi$ is ramified along $b_X^{-1}(C)$. Let $\{\hat {F_j}\}$ denote the set of all ramification curves of type $b_X^{-1}(C)$ and define $F_j= \hat \pi(\hat {F_j})$. The map $\hat \pi$ is a $D$-quotient and ramified along curves 
 $$\mathrm{Fix}_{\hat{X}} (D) = \bigcup \hat{E_i} \cup
  \bigcup \hat{F_j}.$$
\subsection{Mori contractions and $C_7$-fixed points}
Many aspects of the group theory of $G$ can be well understood in term of its generators $\alpha, \beta, \gamma$ of order 7,3,2, respectively. Since the action of $G$ on $\mathbb{P}_2$ is given by a three-dimensional irreducible representation, the action of $G$ on $Z$ is given explicitly in terms of $\alpha, \beta, \gamma$.
We let $S= \langle \alpha  \rangle \cong C_7 < G$ be a cyclic subgroup of order seven in $G$.

The symplectic action of a cyclic group of order seven on an K3-surface has exactly three fixed points.
Since $p_1=[1:0:0]$, $p_2=[0:1:0]$ and $p_3=[0:0:1]$ all lie on $C_\mathrm{Klein} \subset \mathbb P_2$, the action of 
$S$ on $Z$ has exactly three fixed points.

Let $\mathrm{Fix}_{\hat{Y}}(S) =: \{y_1, \dots,y_k\}$ and let
$\mathrm{Fix}_{\hat{X}}(S) =: \{x_1, \dots,x_l\}$. Since blowing-up an $S$-fixed point in $X$ replaces the fixed
point by a rational curve with two $S$-fixed points in $\hat{X}$, we find $3 \leq k \leq l
\leq 6$. 
\begin{lemma}\label{fixSinfixD}
The fixed points of $S$ in $\hat{X}$ are contained in the
$D$-ramification set, i.e., $\mathrm{Fix}_{\hat{X}}(S) \subset
\mathrm{Fix}_{\hat{X}} (D)$.
\end{lemma}
\begin{proof}
Since $D$ centralizes $S$, the action of $D$ stabilizes the $S$-fixed point set.
We first show that $\mathrm{Fix}_{{X}}(S) \subset
\mathrm{Fix}_{{X}} (D)$. Assume the contrary and let $\mathrm{Fix}_X(S)=
\{s_1,s_2,s_3\}$ be a $D$-orbit and $\pi(s_i)=y$. Then $y$ is
a smooth point and fixed by the action of $S$ on $Y$. There exists a neighbourhood of $y$ in $Y$ which is
biholomorphic to a neighbourhood of $b_Y^{-1}(y)= \tilde{y}$ in
$\hat{Y}$. By construction, $\tilde{y}\in \mathrm{Fix}_{\hat{Y}}(S)$. Since $\mathrm{Fix}_{\hat{Y}}(S)$ consists of at least three points, we let $\tilde{\tilde{y}} \neq \tilde{y}$ be an additional $S$-fixed point on
$\hat{Y}$. The fibre $\pi^{-1}(b_Y(\tilde{\tilde{y}}))$ consists of one or three points and is disjoint from $\{s_1,s_2,s_3\}$. Since the point $\tilde{\tilde{y}}$ is a fixed point of $S$, we know that $S \cong C_7$
acts on the fiber $\pi^{-1}(b_Y(\tilde{\tilde{y}}))$ and is seen to fix it pointwise. This is contrary to the fact that $\mathrm{Fix}_X(S)=\{s_1,s_2,s_3\}$. It follows that $\mathrm{Fix}_{{X}}(S) \subset
\mathrm{Fix}_{{X}} (D)$.

It remains to show the corresponding inclusion on $\hat{X}$. If the
points $s_i$ do not coincide with isolated $D$-fixed points, the
statement follows since $b_X$ is equivariant and biholomorphic
outside the isolated $D$-fixed points. 
If $s_i$ is an isolated
$D$-fixed point, we have seen above that the action of $D$ on the blow-up of $s_i$ is trivial. In particular, $\mathrm{Fix}_{\hat{X}}(S) \subset
\mathrm{Fix}_{\hat{X}} (D)$. 
\end{proof}
\subsection*{Excluding the case $|\mathrm{Fix}_{\hat{Y}}(S)|=3$}
\begin{lemma} If $|\mathrm{Fix}_{\hat{Y}}(S)
|=3$, then $\mathrm{Fix}_{\hat{Y}}(S) \cap \bigcup E_i = \emptyset$. 
\end{lemma}
\begin{proof}
Fixed points
of $S$ on a curve $\hat{E_i}$ always come in pairs: If the curve $\hat{E_i}$ contains a fixed point of $S$, then the isotropy
representation of $S$ at the fixed point $b_X(\hat{E_i})$ in $X$ defines an action of the cyclic group $S$ on the rational curve $\hat{E_i}$ with exactly two fixed points.
If $|\mathrm{Fix}_{\hat{Y}}(S)|=|\mathrm{Fix}_{\hat{X}}(S)
|=3$ and $\mathrm{Fix}_{\hat{Y}}(S) \cap \bigcup E_i \neq \emptyset$, then two of the $S$-fixed point lie on the same curve $\hat{ E_i}$ and $|\mathrm{Fix}_{X}(S)| \leq 2$, a contradiction.
\end{proof}
\begin{lemma}\label{fixed points on mori fibers}
If $|\mathrm{Fix}_{\hat{Y}}(S)
|=3$, then the set $\mathrm{Fix}_{\hat{Y}}(S)$ has empty intersection with the exceptional locus of the full equivariant Mori reduction $M_{\mathrm{red}}: \hat{Y} \to Z$. 
\end{lemma}
\begin{proof}
Let $C$ be any exceptional curve of the Mori reduction and assume there is a fixed point of $S$ on $C$. As the point $p$ obtained from blowing down $C$ has to be a fixed point of $S$, it follows that the curve $C$ is $S$-invariant. In particular, we know that the action of $S$ on $C$ has exactly two fixed points. Now blowing down $C$ reduces the number of $S$-fixed point by 1. This contradicts the fact that $|\mathrm{Fix}_{Z}(S)|=3$. 
\end{proof}
\begin{lemma}\label{shapeofSaction}
Let $|\mathrm{Fix}_{\hat{Y}}(S)|=3$ and let $p \in \mathrm{Fix}_Z (S)$. Then there exist local coordinates $(u,v)$ at $p$ and a nontrivial character $\mu: S \to \mathbb{C}^*$ such that the action of $S$ at $p$ is locally given by either
\[
(u,v) \mapsto (\mu^3 u, \mu^{-1} v)\quad \text{or} \quad (u,v) \mapsto (\mu u, \mu^{-3} v).
\]
\end{lemma}
\begin{proof}
On the K3-surface $X$ the action of $S$ at a fixed point is in local coordinates $(z,w)$ given by $(z,w) \mapsto (\mu
z, \mu^{-1}w)$ for some nontrivial character $\mu: S \to \mathbb{C}^*$. Since $\mathrm{Fix}_{\hat{Y}}(S) \cap \bigcup E_i = \emptyset$, the map $b_X$ is biholomorphic in a neighbourhood of the fixed point. 
Recalling that $\mathrm{Fix}_{\hat{X}}(S)$ is contained in the ramification locus of $\hat{\pi}$ (i.e., $p \in \mathrm{Fix}_{\hat{X}}(D)$) the action of $D$ may be linearized at $p$. 
Since $S$ and $D$ commute, the action of $D$ is diagonal in the chosen local coordinates $(z,w)$. We conclude that $\hat \pi$ is locally of the form $(z,w) \mapsto (z^3,w)$ or $(z,w^3)$. The action of $S$ at a fixed point in $\hat{Y}$ is defined by $(\mu^3, \mu^{-1})$ or $(\mu, \mu^{-3})$, respectively.
By the lemma above, the fixpoints of $S$ are not affected by the Mori reduction. The map $M_{\mathrm{red}}$ is $S$-equivariant and locally biholomorphic in a neighbourhood of a fixed point of $S$. The lemma follows.
\end{proof}
Using our explicit knowledge of the $G$-action on $Z$ we will show in the following that the linearization of the action of $S < G$ at a fixed point in the Del Pezzo surface $Z$ is not of the type described by the lemma above. We distinguish two cases when studying $Z$.

Let $Z \cong \mathbb{P}_2$ and $[x_0:x_1:x_2]$ denote homogeneous coordinates on $\mathbb{P}_2$ such that the action of $S <G$ on $\mathbb{P}_2$ is given by $[x_0:x_1:x_2] \mapsto [\zeta x_0, \zeta^2 x_1, \zeta ^4 x_2 ]$ where $\zeta$ is a $7^\text{th}$ root of unity. Using affine coordinates $z= \frac{x_1}{x_0}, w= \frac{x_2}{x_0}$ we check that the action of $S$ at $p_1=[1:0:0]$ is locally given by $(z,w) \mapsto (\zeta z, \zeta^3 w)$. This contradicts Lemma \ref{shapeofSaction}.

Let $Z \overset{q}{\to} \mathbb{P}_2$ be the double cover of $\mathbb{P}_2$ branched along Klein's quartic curve and 
let $[x_0:x_1:x_2]$ denote homogeneous coordinates on $\mathbb{P}_2$. As above, using affine coordinates $u= \frac{x_1}{x_0}, v= \frac{x_2}{x_0}$ we check that the action of $S$ in a neighbourhood of $[1:0:0]$ is locally given by $(u,v) \mapsto (\zeta u, \zeta^3 v)$. The branch curve $C_\mathrm{Klein} \subset \mathbb{P}_2$ is defined by the equation $u^3+uv^3+v$. In new coordinates $(\tilde u (u,v), \tilde v(u,v))= (u, u^3+uv^3+v)$ the branch curve is defined by $\tilde v = 0$ and the action of $S$ is given by $(\tilde u,\tilde v) \mapsto (\zeta \tilde u, \zeta^3 \tilde v)$. Consider the fixed point $[1:0:0] \in \mathbb P_2$ and its preimage  $p \in Z$. At $p$, coordinates $(z,w)$ can be chosen such that the covering map is locally given by $(z,w) \mapsto (z,w^2) = (\tilde u, \tilde v)$. It follows that the action of $S$ at $p \in Z$ is locally given by $(z,w) \mapsto (\zeta z, \zeta^5 w)$. This is again contrary to Lemma \ref{shapeofSaction}.

In summary, if $|\mathrm{Fix}_{\hat{Y}}(S)|=3$, the action of $S < G$ on the Del Pezzo surface $Z$ cannot be induced by a symplectic $C_7$-action on the K3-surface $X$. This proves the following lemma.

\begin{lemma}
$|\mathrm{Fix}_{\hat{Y}}(S)|\geq 4$.

\end{lemma}
\subsection{Lifting Klein's quartic}
The discussion of the previous section shows that there must be a step in the Mori reduction where the blow-down of a (-1)-curve identifies two $S$-fixed points. Let $z \in Z$ be a fixed point of $S$. Then, by equivariance, all points in the $G$-orbit of $z$ are obtained by blowing down (-1)-curves in the process of Mori reduction. 
If $Z \cong \mathbb{P}_2$, we denote by $C_\mathrm{Klein} \subset Z$ Klein's quartic curve. 
If $Z$ is the double cover of $\mathbb{P}_2$ branched along Klein's curve, we abuse notation and denote by $C_\mathrm{Klein}$ the ramification curve in $Z$. In the later case $C_\mathrm{Klein}$ is a $G$-invariant curve of genus 3 and self-intersection 8 by Lemma \ref{selfintbranch}.

Let $z \in \mathrm{Fix}_Z (S) \subset C_\mathrm{Klein}$ and consider the $G$-orbit $G\cdot z$. By invariance, $G\cdot z \subset C_\mathrm{Klein}$. The isotropy group $G_z$ must be cyclic and $G_z =S$ implies $|G\cdot z|=24$.  Let $B$ denote the strict transform of $C_\mathrm{Klein}$ in $\hat Y$. The curve $B$ is a smooth $G$-invariant curve of genus 3 and meets at least 24 Mori fibers. Applying Lemma \ref{selfintblowdown} to $M_\mathrm{red}(B)=C_\mathrm{Klein}$ we obtain
\[
B^2 \leq C_\mathrm{Klein}^2 -24 \leq -8.
\]
\begin{lemma}
The curve $B$ does not coincide with any of the curves of type $E$ or $F$. Its preimage $\hat B := \hat \pi^{-1}(B) \subset \hat X$ is a cyclic degree three cover of $B$ branched at $B \cap(\bigcup E_i \cup \bigcup F_j)$.
\end{lemma}
\begin{proof}
The curves $E_i \subset \hat Y$ are (-3)-curves whereas $B$ has self-intersection less than or equal to $-8$. Assume $B = F_j$ for some $j$. Then $\hat B$ is a curve of self-intersection less than or equal to $-4$ by Lemma \ref{selfintbranch} which is mapped biholomorphically to the K3-surface $X$. We obtain a contradiction since K3-surfaces do not admit curves of self-intersection less than $-2$. 
\end{proof}
Since  $\mathrm{Fix}_Z (S) \subset C_\mathrm{Klein}$ there are three fixed points of $S$ on $\hat B$. From $\mathrm{Fix}_{\hat{X}}(S) \subset \mathrm{Fix}_{\hat{X}} (D)$ it follows that $\hat \pi|_{\hat B}: \hat B \to B$ is branched at three or more points.  In particular, the curve $\hat{B}$ is connected. 
In the following, we will distinguish two cases: the curve $\hat B$ being reducible or irreducible. 
\subsection*{Case 1: The curve $\hat{B}$ is reducible}
\addcontentsline{toc}{subsection}{\hspace{1.1cm} Case 1: The curve $\hat B$ is reducible}
The three irreducible components $\hat{B}_i$, $i= 1,2,3$ of $\hat B$ are smooth curves which are mapped biholomorphically onto $B$. Since $B$ is exceptional, the configuration of curves $\hat{B}$ is also exceptional. 
It follows that the intersection matrix $(\hat{B}_i \cdot \hat{B}_j)_{ij}$ is negative definite. In the following we study the intersection matrix of $\hat B$ and will obtain a contradiction.

The restricted map $b_X: \hat B_i \to b_X(\hat B_i)$ is the normalization of $b_X(\hat B_i)$ and consequently the arithmetic genus of $b_X(\hat{B}_i)$ is given by the formula (cf. II.11 in \cite{BPV})
\[
g(b_X(\hat{B}_i)) = g(\hat {B_i}) + \delta(b_X(\hat{B}_i)),
\]
where the number $\delta$ is computed as 
$\delta(b_X(\hat{B}_i)) = \sum_{p \in b_X(\hat{B}_i)} \mathrm{dim}_\mathbb{C}({b_X}_*\mathcal{O}_{\hat{B}_i}/\mathcal{O}_{b_X(\hat{B}_i)})_p$.
Note that the sum can also be taken over the singular points $p \in b_X(\hat{B}_i)$ only, since smooth points do not contribute to the sum.
Since $X$ is a K3-surface, the adjunction formula for $b_X(\hat{B}_i)$ reads
\[
(b_X(\hat{B}_i))^2 = 2 g(b_X(\hat{B}_i)) -2 = 2g(\hat B_i) +2 \delta(b_X(\hat{B}_i)) -2.
\]
By Lemma \ref{selfintblowdown}, the self-intersection number $(b_X(\hat{B}_i))^2$ can be expressed in terms of the self-intersection $\hat{B}_i^2$ and intersection multiplicities $E_j\cdot \hat{B}_i$:
\[
(b_X(\hat{B}_i))^2 = \hat{B}_i^2 + \sum_j (\hat{E}_j \cdot \hat{B}_i)^2.
\]
It follows that the self-intersection number of $\hat{B}_i$ can be expressed as
\begin{equation}\label{selfintB}
\hat{B}_i^2 = 2g(\hat B_i) +2 \delta(b_X(\hat{B}_i)) -2 -\sum_j (\hat{E}_j\cdot \hat{B}_i)^2.
\end{equation}
For simplicity, we first consider the case where $\hat{B}_i$ has nontrivial intersection with only one curve of type $\hat{E}$. We refer to this curve as $\hat{E}$. The general case then follows by addition over all curves $\hat E_j$, the number $\delta$ for the full contraction $b_X$ is the sum of all numbers $\delta$ obtained when blowing down disjoint curves $\hat{E}_j$ stepwise.
\paragraph{Estimating the number $\delta$}
\begin{example}
Let $C= C_1 \cup C_2$ be a connected curve consisting of two irreducible components. Then the arithmetic genus of $C$ is calculated as
$g(C)= g(C_1) + g(C_2) + C_1 \cdot C_2 -1$.
The normalization $ \tilde C$ of $C$ is given by the disjoint union of the normalizations $\tilde C_i$ of $C_1$ and $C_2$. In particular,
$g(\tilde C) = g( \tilde C_1) + g(\tilde C_2) -1$,
so that $\delta(C) = \delta(C_1) + \delta(C_2) + C_1 \cdot C_2$
(cf. II.11 in \cite{BPV}).
\end{example}
 Since the number $\delta$ is a sum of contributions $ \delta_p$ at singular points $p$, we can calculate the number $\delta_p$ locally at each singularity where we decompose the germ of the curve as the union of irreducible components and use a formula generalizing the example above. We refer to an irreducible component of a curve germ realized in a open neighbourhood of the surface as a \emph{curve segment}.

In order to study the singularities of $b_X(\hat B_i)$ one needs to consider the points of intersection $\hat E \cap \hat B_i$. These points of intersection can be of different quality:
\begin{itemize}
\item
\textbf{Type $ m= 1$:} The intersection at $b \in \hat B_i$ is transversal and the local intersection multiplicity at $b$ is equal to 1. A neighbourhood of $b$ in $\hat B_i$ is mapped to a smooth curve segment in $b_X(\hat B_i)$. 
\item
\textbf{Type $m>1$:} The intersection at $b \in \hat B_i$ is of higher multiplicity $m(b)$, i.e., $\hat E$ is tangent to $\hat B_i$ and in local coordinates $(z,w)$ we may write $\hat E = \{ z=0\}$ and $\hat B_i = \{  z- w^m\}$. Blowing down $\hat E$ transforms a neighbourhood of $b$ into a a curve segment isomorphic to $\{ x^{m+1} -y^m =0\}$. For the singularity $(0,0)$ of this curve we calculate
\[
\delta_{(0,0)}= \frac{1}{2}m(m-1).
\]
\end{itemize}
Let $b_m$ denote the number of points  in $\hat E \cap \hat B_i$ with local intersection multiplicity $m$. For each point of intersection of $\hat E$ and $\hat B_i$ we obtain an irreducible component of the germ of $b_X(\hat B_i)$ at $p = b_X(\hat E)$. We compute $\delta_p$ by decomposing this germ and need to determine local intersection multiplicities of all combinations of irreducible components.
\begin{lemma}\label{normarg}
Two irreducible components of the germ of $b_X(\hat B_i)$ at $p$ corresponding to points in $\hat E \cap \hat B_i$ of type $m$ and $n$ meet with local intersection multiplicity greater than or equal to $mn$.
\end{lemma}
\begin{proof}
In order to determine the intersection multiplicity of two irreducible components corresponding to points of type $m$ and $n$, we write one curve as $\{ x^{m+1} -y^m =0\}$. The second curve can be expressed as $\{h_1(x,y)^{n+1} - h_2(x,y)^n=0\}$ where $(x,y) \mapsto (h_1(x,y),h_2(x,y))$ is a holomorphic change of coordinates. Now normalizing the first curve by $\xi \mapsto (\xi^m, \xi^{m+1})$ and pulling back the equation of the second curve to the normalization $\mathbb{C}$, we obtain the equation
$h_1(\xi^m, \xi^{m+1})^{n+1} - h_2(\xi^m, \xi^{m+1})^n=0$
which has degree at least $mn$ in $\xi$. It follows that the local intersection multiplicity is greater than or equal to $mn$. 
\end{proof}
Counting different types of intersections of irreducible components we obtain the following estimate for $\delta_p$
\begin{align*}
\delta_p &= \sum \delta_p(C_i) + \sum_{i\neq j}(C_i \cdot C_j)_p\\
&\geq \sum_{m\in\mathbb N}\frac{b_m}{2}m(m-1) +  \frac{1}{2}\sum_{m \in \mathbb N}b_m(b_m-1)m^2 + \sum_{m> n} b_mb_n mn
\end{align*}

where $\sum_{i\neq j}(C_i \cdot C_j)_p$ decomposes into intersections $(C_i \cdot C_j)_p$ of type $mm$ and intersections of type $mn$ for $m \neq n$. The formula above applies to each curve $\hat E_j$ having nontrivial intersection with $\hat B_i$.
Let $p_j$ be the point on $X$ obtained by blowing down $\hat E_j$ and let $b_m^j$ denote the number of points of type $m$ in $\hat B_i \cap \hat E_j$. Then
\begin{align*}
\delta ( b_X(\hat B_i)) &= \sum_j \delta_{p_j}(b_X(\hat B_i)) \\
&\geq \sum_j(\sum_{m\in\mathbb N}\frac{b_m^j}{2}m(m-1) + \frac{1}{2}\sum_{m\in \mathbb N}b_m^j(b_m^j-1)m^2 + \sum_{m> n} b_m^jb_n^j mn).
\end{align*}
Returning to the formula (\ref{selfintB}) for $\hat B_i^2$ we obtain
\begin{align*}
\hat{B}_i^2 &= 2g(\hat B_i) +2 \delta(b_X(\hat{B}_i)) -2 -\sum_j (\hat{E}_j \cdot \hat{B}_i)^2\\
& \geq \sum_j(\sum_{m\in\mathbb{N}}b_m^jm(m-1) + \sum_{m\in \mathbb{N}}b_m^j(b_m^j-1)m^2 + 2\sum_{m> n} b_m^jb_n^j mn)\\
&-2 - \sum_j(\sum_m b_m^jm)^2\\
& \geq -2-\sum_j \sum_m b_m^j m.
\end{align*}
As a next step, we will find a bound for $(\hat B_i \cdot \hat B_k)$ in the
case $i \neq k$. If a curve $\hat B_i$ intersects a ramification curve of type $\hat
E$ or $\hat F$ in a point $x$, then $(\hat B_i \cdot \hat B_k)_x \geq
1$. If $(\hat B_i \cdot \hat E_j)_x =m$, then for $k \neq i$ 
\[
(\hat B_k \cdot\hat E_j)_x = (\varphi_D( \hat B_i)  \cdot \hat E_j)_x= ( \varphi_D(\hat B_i)  \cdot \varphi_D(\hat E_j))_x = (\hat B_i \cdot \hat E_j)_x =m
\] 
where $\varphi_D \in D$ is a biholomorphic transformation and $E_j$
is in the fixed locus of $D$.
\begin{lemma}
Assume $\hat B_i$ meets a curve of type $\hat E$ or $\hat F$ in $x$ with local intersection multiplicity  $m$. Then $(\hat B_i \cdot \hat B_k)_x \geq  m$.
\end{lemma}
\begin{proof}
Let $\hat E$, $\hat F$ respectively, be locally given by
$\{z=0\}$. Then $\hat B_i$ is locally given by $\{z-w^m=0\}$ and $\hat
B_k$ by $\{h_1(z,w) - h_2(z,w)^m =0\}$ where $(z,w) \mapsto
(h_1(z,w), h_2(z,w))$ is, as in the proof of Lemma \ref{normarg}, a holomorphic
change of coordinates. Note that it stabilizes $\{z=0\}$,
i.e., $h_1(0,w)=0$ for all $w$ and we can write $h_1(z,w) =z
\tilde{h}_1(z,w)$. The intersection of $\hat B_i$ and $\hat B_j$ corresponds to the equation
$ w^m \tilde{h}_1(w^m, w) - h_2(w ^m, w)$
which is of degree greater than or equal to $m$. The lemma follows. 
\end{proof}

Summing over all points of intersection of $\hat B_i$ and $ \hat B_k$ one finds
$\hat B_i \cdot \hat B_k \geq \sum_j \sum_m b_m^j m$. 
Recall that by Lemma \ref{fixSinfixD}
$\mathrm{Fix}_{\hat X}(S)$ is contained in $\mathrm{Fix}_{\hat X}(D)$ and that the curve $B$ contains
three $S$-fixed points. Therefore, it intersects the
ramification locus of $\hat \pi$ in at least three points. At these points the three irreducible components of $\hat B$ must meet. In particular, $(\hat B_i, \hat B_k) \geq 3 $. This yields
\begin{align*}
&(1,1,1) \begin{pmatrix}
\hat B_1^2 & \hat B_1 \cdot \hat B_2 &\hat B_1 \cdot \hat B_3\\
\hat B_2 \cdot \hat B_1 & \hat B_2^2 & \hat B_2 \cdot \hat B_3\\
\hat B_3 \cdot \hat B_1 & \hat B_3 \cdot \hat B_2 & \hat B_3^2
        \end{pmatrix}
                \begin{pmatrix}
                1\\1\\1
                \end{pmatrix}\\
&= \hat B_1^2 + \hat B_2^2 + \hat B_3^2 + 2(\hat B_1 \cdot \hat B_2 + \hat B_2 \cdot \hat B_3 + \hat B_1 \cdot \hat B_3)\\
& \geq -6 - 3 \sum_j \sum_m b_m^j+3\sum_j \sum_m b_m^j m +(\hat B_1 \cdot \hat B_2 + \hat B_2 \cdot \hat B_3 + \hat B_1 \cdot \hat B_3) \\
&= -6 +(\hat B_1 \cdot \hat B_2 + \hat B_2 \cdot \hat B_3 + \hat B_1 \cdot \hat B_3)\\
&\geq 3.
\end{align*}
Hence, the intersection matrix $(\hat B_i \cdot \hat B_j)_{ij}$ is not negative-definite contradicting the fact that $\hat B$ is exceptional. It follows that the curve $\hat B$ must be irreducible.
\subsection*{Case 2: The curve $\hat B$ is irreducible}
\addcontentsline{toc}{subsection}{\hspace{1.1cm} Case 2: The curve $\hat B$ is irreducible}
Let $n: N \to \hat B$ be the normalization of $\hat B$. Since $b_X$ is a blow-up, $b_X \circ n: N \to b_X(\hat B)$ is the normalization of the curve $b_X(\hat B) \subset X$. It follows that
$g(b_X(\hat{B})) = g(N) + \delta(b_X(\hat{B}))$.
By adjunction, the self-inter\-sec\-tion of $b_X(\hat B)$ is given by
\[
(b_X(\hat B))^2 = 2g(b_X(\hat{B}))- 2 = 2g(N) + 2\delta(b_X(\hat{B})) -2.
\]
As above, by Lemma \ref{selfintblowdown},
$(b_X(\hat{B}))^2 = \hat{B}^2 + \sum_j (\hat{E}_j \cdot \hat{B})^2$.
Thus, the self-in\-ter\-sec\-tion of $\hat B$ can be expressed as 
\[
\hat{B}^2 = 2g(N) + 2\delta(b_X(\hat{B})) -2 - \sum_j (\hat{E}_j \cdot \hat{B})^2.
\]
Since the curve $\hat B$ is exceptional, this self-intersection number must be negative. By finding a lower bound for $\hat B^2$ we will obtain a contradiction.

Let us first examine the points of intersection $\hat B \cap \hat E$ for one curve $\hat E$ among the exceptional curves of the blow-down $b_X$. 
We consider the corresponding points of intersection of $B$ and $E$ in
$\hat Y$ and we choose coordinates $(\xi, \eta)$ such that $E$ is locally defined by $\{\xi =0\}$, the map $\hat \pi $ is locally given by $(z,w) \mapsto (z^3,w)=(\xi, \eta)$ and $B = \{f(\xi, \eta)=0\}$. It follows that $\hat B$ is locally defined by $\{ h= f\circ \hat\pi =0\}$. 

If $E$ and $B$ meet transversally, we know that the function $f(\xi, \eta)$ fulfills $\frac{\partial f}{\partial \eta}|_{(0,0)} \neq 0$. It follows that  $\frac{\partial h}{\partial w}|_{(0,0)} \neq 0$ and after a suitable change of coordinates $h(z,w)= z^m-w$. 

If $E$ and $B$ meet tangentially, we know that the function $f(\xi, \eta)$ fulfills $\frac{\partial f}{\partial \eta}|_{(0,0)} = 0$. Since $B$ is smooth, we know $\frac{\partial f}{\partial \xi}|_{(0,0)} \neq 0$. After a suitable change of coordinates $h(z,w)= z^3-w^n$ with $n>0$. Note that in both cases the coordinate change on $\hat X$ is such that $\hat E$ is still defined by $\{z=0\}$. This will be important when describing the blow-down $b_X$ of $\hat E$. 

Consider a curve segment $\{h=0\}$ in $\hat X$ and its image under the map $b_X$. If $h(z,w)= z^m-w$ then the corresponding smooth segment of $b_X(\hat B)$ is defined by $ x^{m+1} -y =0$. If $h(z,w)= z^3-w^n$ then the corresponding piece of $b_X(\hat B)$ is defined by $ x^{n+3} -y^n =0$ and has a singular point if $n>1$. 

Let $p =b_X(\hat E)$. We will determine $\delta_p$ by decomposing the germ of $b_X(\hat B)$ at $p$ into its irreducible components. There are three different types of such components:
\begin{enumerate}
\item 
smooth components locally defined by $ x^{m+1} -y =0$,
\item
singular components locally defined by $ x^{n+3} -y^n =0$ for $n>1$ not divisible by $3$,
\item
 triplets of smooth components locally defined by $ x^6 -y^3 =0$,
\item
 triplets of singular components locally defined by $ x^{n+3} -y^n =0$ for $n=3k$ and $k >1$. 
\end{enumerate}
The singularity in case 2) gives $\delta = \frac{n^2+n-2}{2}$. In case 4), each component is defined by an equation of type $x^{k+1}-y^k=0$ and the singularity of each component gives $\delta = \frac{k^2-k}{2}$. \\
In order to determine $\delta_p$ we need to specify intersection multiplicities for all combinations of irreducible components.
\begin{lemma}
The local intersection multiplicities of 
 pairs of irreducible components of the germ of $b_X(\hat B)$ at $p$ in general position are given by the following table.
\renewcommand{\baselinestretch}{1.5}
\begin{center}
\begin{table}[H]
\begin{center}
\begin{tabular}{c|c|c|c|c}
local equation & $ x^{m_1+1} -y  $ & $ x^{n_1+3} -y^{n_1}$ & $x^2-y$ & $x^{k_1+1}-y^{k_1}$\\ \hline
$ x^{m_2+1} -y $ & 1& $n_1$ & 1& $k_1$\\ \hline
$ x^{n_2+3} -y^{n_2}  $& $n_2$ & $n_1n_2$ & $n_2$ & $n_2k_1$ \\ \hline
$x^2-y $ & 1& $n_1$ & 1 \text{\ or\ }(2) & $k_1$ \\ \hline 
$x^{k_2+1}-y^{k_2}$ & $k_2$ & $k_2n_1$ & $k_2$ & $k_1k_2$ \text{\ or\ }($k^2+k$)
\end{tabular}
\end{center}
\end{table}
\end{center}
\renewcommand{\baselinestretch}{1.1}
\end{lemma}
Note that the local equations in the first row and column, although
all written as functions of $(x,y)$, describe the curve segments in
different choices of local coordinates. 
\begin{proof}[Sketch of proof]
As above, we rewrite one equation as
$f(h_1(x,y),h_2(x,y))$ where $(h_1,h_2)$ is a holomorphic change of local 
coordinates. The intersection multiplicities can then be calculated by
the method introduced in the proof of Lemma \ref{normarg}. 
Two irreducible components in a triplet of type 3) meet with
intersection multiplicity 2. Two irreducible components in a triplet
of type 4) meet with intersection multiplicity $k^2+k$. These
quantities are indicated in brackets as they differ from the intersection multiplicities of two irreducible components from different triplets.
\end{proof}
\begin{remark}
If two irreducible components of the germ of $b_X(\hat B)$ at $p$ are in special position, their local intersection multiplicity is greater than the value specified in the above table. In particular, the table gives lower bounds for the respective intersection numbers.
\end{remark}
Let $a$ denote the number of irreducible components of type 1), let $b_n$ the number of irreducible components of type 2) where $ n \not\in 3\mathbb{N}$, let $c \in 3\mathbb{N}$ denote the number of irreducible components of type 3) and let $d_k \in 3 \mathbb{N}$ denote the number of irreducible components of type 4). We summarize $e=a+c$.

A lower bound for $\delta_p$ is given by
\begin{align*}
\delta_p &\geq  \sum_n b_n \frac{n^2+n-2}{2} + \sum_k d_k \frac{k^2-k}{2}\\
       & + \frac{1}{2}e(e-1) + c +\sum_n e b_n n   +  \sum_k e d_k k\\
       & + \frac{1}{2}\sum_n b_n(b_n-1)n^2  +\sum_{n_1 > n_2} b_{n_1}b_{n_2} n_1n_2 +  \sum_{n,k} b_nd_k nk\\
        & + \frac{1}{2}\sum_k d_k(d_k-1)k^2 + \sum_k d_k k+ \sum_{k_1 > k_2} d_{k_1}d_{k_2} k_1k_2.
\end{align*}
For simplicity, we first consider only one curve $\hat E$ intersecting $\hat B$. The formula for $\hat B^2$ becomes 
\begin{align}
\hat{B}^2 &= 2g(N) + 2\delta(b_X(\hat{B})) -2 - (\hat{E} \cdot \hat{B})^2 \notag\\
&=2g(N) + 2\delta(b_X(\hat{B})) -2 - \underset{(\hat{E} \cdot \hat{B})^2}{\underbrace{(e+ \sum_n b_n n + \sum_k d_k k)^2}} \notag\\
&=2g(N)-2- e +2c+ \sum_k d_k k + \sum_n b_n(n-2) \notag \\
&\geq 2g(N)-2- a. \label{ineq}
\end{align}
The same formula also holds if we consider the general case of curves $\bigcup_i \hat E_i$ intersecting $\hat B$ since both the calculation of $\delta$ and the intersection number $\sum_i(\hat B, \hat E_i)^2$ can be obtained from the above by addition. The number $a$ now represents the number of points of type 1) in the union of curves $\hat E_i$. 

The map $n \circ \hat \pi: N \to \hat B \to B$ is a degree three cover of the smooth curve $B$ branched at $V \subset B$. The genus of $B$ is three, the topological Euler characteristic is $e(B)= -4$. Let $\tilde V := B \cap (\bigcup E_i \cup \bigcup F_j)$ denote the branch locus of $\hat \pi: \hat B \to B$. Then $V \subset \tilde V$ and $V$ must contain those points in $\tilde V$ which correspond to smooth points on $\hat B$. In partcular, $|V| \geq a$. 

The Euler characteristic of $N$ is given by
$e(N) = 3e(B) - 2|V| = -12 -2|V| = 2- 2g(N)$.
and inequality \eqref{ineq} becomes
\[
\hat B ^2 \geq 12 + 2|V|-a \geq 12 + |V| \geq 0 
\]
contradicting the fact that $\hat B$ is exceptional.
\subsection*{Conclusion} 
\addcontentsline{toc}{subsection}{Conclusion}
The above contradiction shows the non-existence of a K3-surface with an action of $G \times C_3$. This completes the prove of Theorem \ref{nonexist}.
%
%
%
%
\chapter{The alternating group of degree six}\label{chapterA6}
In the previous chapters we have considered symplectic automorphisms groups of K3-surfaces centralized by an antisymplectic involution, i.e., the groups under consideration were of the form $\tilde G = G \times \langle \sigma \rangle$ where $\tilde G_\mathrm{symp} = G$. In this chapter we wish to discuss more general automorphims groups $\tilde G$ of mixed type: if $\tilde G $ contains an antisymplectic involution $\sigma$ with fixed points we consider the quotient by $\sigma$. In general, if $\sigma$ does not centralize the group $\tilde G_\mathrm{symp}$ inside $\tilde G$, the action of $\tilde G_\mathrm{symp}$ does \emph{not} descend to the quotient surface. We therefore restrict our consideration to the centralizer $Z_{\tilde G}(\sigma)$ of $\sigma$ inside $\tilde G$ (or $\tilde G_\mathrm{symp}$) and study its action on the quotient surface. 

If we are able to describe the family of K3-surfaces with $Z_{\tilde G} (\sigma)$-symmetry, it remains to detect the surfaces with $\tilde G$-symmetry inside this family. 
This chapter is devoted to a situation where the group $\tilde G$ contains the alternating group of degree six.
Although, a precise classification cannot be obtained at present, we achieve an improved understanding of the equivariant geometry of K3-surfaces with $\tilde G$-symmetry and classify families of K3-surfaces with $Z_{\tilde G}(\sigma)$-symmetry (cf. Theorem \ref{classiA6}). In this sense, this closing chapter serves as an outlook on how the method  of equivariant Mori reduction allows generalization to more advanced classification problems.
\section{The group $\tilde A_6$}
We let $\tilde G$ be any finite group which fits into the exact sequence
\[
\{\mathrm{id}\} \to A_6 \to \tilde G \overset {\alpha}{\to} C_n \to \{\mathrm{id}\}.
\]
and in the following consider a K3-surface $X$ with an effective action of $\tilde G$. 
The group of
symplectic automorphisms $(\tilde{G})_{\text{symp}}$ in $\tilde G$
coincides with $A_6$. 

This particular situation is considered by Keum, Oguiso, and Zhang in \cite{KOZLeech} and \cite{KOZExten}. They lay special emphasis on the maximal possible choice of $\tilde G$ and therefore consider a group $\tilde G = \tilde A_6$ characterized by the exact sequence
\begin{equation}\label{tilde A6}
\{\mathrm{id}\} \to A_6 \to \tilde A_6 \overset {\alpha}{\to} C_4 \to \{\mathrm{id}\}.
\end{equation}

Let $N := \mathrm{Inn}(\tilde{A_6}) \subset
\mathrm{Aut}(A_6)$ denote the group of inner automorphisms of $\tilde A_6$ and let $\mathrm{int} : \tilde A_6 \to N$ be the homomorphisms mapping an element $g \in \tilde A_6$ to conjugation with $g$.
 It can be shown that the group $\tilde{A_6}$ is a semidirect product $A_6
\rtimes C_4$ embedded in $N\times C_4$ by the map $(\mathrm{int}, \alpha)$
(Theorem 2.3 in \cite{KOZExten}). By Theorem 4.1 in \cite{KOZExten} the group $N$ is isomorphic to $M_{10}$
  and the isomorphism class of $\tilde A_6$ is uniquely determined by \eqref{tilde A6} and the condition that it acts on a K3-surface.

In \cite{KOZLeech} a lattice-theoretic proof of the following classification result (Theorem 5.1, Theorem 3.1, Proposition 3.5) is given. 
\begin{theorem}
A K3 surface $X$ with an effective action of $\tilde A_6$ is isomorphic to the minimal desingularization of the surface in $\mathbb P_1 \times \mathbb P_2$ given by 
\[
S^2(X^3+Y^3 + Z^3) -3 (S^2 + T^2) XYZ =0.
\]
\end{theorem}
Although this realization is very concrete, the action of $\tilde A_6$ on this surface is hidden. The existence of an isomorphism from a K3-surface with $\tilde A_6$-symmetry to the surface defined by the equation above follows abstractly since both surfaces are shown to have the same transcendental lattice. It is therefore desirable to achieve a more geometric understanding of K3-surfaces with $\tilde A_6$-symmetry in general and in particular to obtain an explicit realization of $X$ where the action of $\tilde A_6$ is visible.

We let the generator of the
factor $C_4$ in the semidirect product $\tilde{A_6} = A_6
\rtimes C_4$ be denoted by $\tau$. The order four
automorphism $\tau$ is nonsymplectic and has fixed points. It follows that the antisymplectic
involution $\sigma := \tau^2$ fulfils
\[
\mathrm{Fix }_X(\sigma) \neq \emptyset.
\] 
Since $\sigma$ is mapped to the trivial automorphism in
$\mathrm{Out}(A_6) = \mathrm{Aut}(A_6)/\mathrm{int}(A_6) \cong C_2 \times C_2$
there exists $h \in A_6$ such that $
\mathrm{int}(h) = \mathrm{int}(\sigma) \in \mathrm{Aut}(A_6)$.
The antisymplectic involution $h \sigma$ centralizes $A_6$ in $\tilde{A_6}$. 
\begin{remark}
If $\mathrm{Fix}_X(h \sigma) \neq \emptyset$, we are in the
situation dealt with in Section \ref{A6Valentiner}, i.e., the K3-surface $X$ is an $A_6$-equivariant double
cover of $\mathbb{P}_2$ where $A_6$ acts as Valentiner's group and the branch locus is given by $F_{A_6}(x_1,x_2,x_3) = 10 x_1^3x_2^3+ 9 x_1^5x_3 + 9 x_2^3x_3^3-45 x_1^2 x_2^2 x_3^2-135 x_1 x_2 x_3^4 + 27 x_3^6$. By construction, there is an evident action of $A_6 \times C_2$ on the Valentiner surface\index{Valentiner surface!}, it is however not clear whether this surface admits the larger symmetry group $\tilde A_6$.
\end{remark}
In the following we assume that $h \sigma$ acts without fixed points on $X$ as otherwise the remark above yields an $A_6$-equivariant classification of $X$.
\subsection{The centralizer $G$ of $\sigma$ in $\tilde{A_6}$}\label{centralizer}
We study the quotient $ \pi : X \to X/\sigma = Y$. As mentioned above, the action of the centralizer of $\sigma$ descends to an action on $Y$. We therefore start by identifying the centralizer $G :=Z_{\tilde{A_6}}(\sigma)$ of $\sigma$ in
$\tilde{A_6}$. 
\begin{lemma}
The group $G$ equals $Z_{A_6}(\sigma) \rtimes C_4$ and
$Z_{A_6}(\sigma) = Z_{A_6}(h)$
\end{lemma}
\begin{proof}
The lemma follows from direct computations:
 we write an element of $\tilde A_6$ as $a\tau^k$ with $a \in A_6$. Then $a \tau^k$ is in $Z_{\tilde{A_6}}(\sigma)$ if and only if $a \tau^k \tau^2 = \tau^2 a \tau^k$. This is the case if and only if $a \tau^2 = \tau^2 a$, i.e., if $a \in Z_{A_6}(\sigma)$. Now $\langle \tau \rangle < Z_{\tilde{A_6}}(\sigma)$ implies the first part of the lemma. The second part follows from the equality $\mathrm{int}(\sigma) = \mathrm{int}(h)$. 
\end{proof}
\begin{lemma}
$Z_{A_6}(h) = D_8$
\end{lemma}
\begin{proof}
Since $\mathrm{int}(\sigma) = \mathrm{int}(h)$ and $\sigma^2 = \mathrm{id}$, it follows that $h^2$ commutes with any element in $A_6$. As $Z(A_6)= \{ \mathrm{id} \}$, it follows that $h$ is of order two. There is only one conjugacy class of elements of order two in $A_6$. We calculate $Z_{A_6}(h) = D_8$ for one particular choice of $h$.
Let 
\[
h=
\begin{pmatrix}
 1 & 2 & 3 & 4 & 5 & 6 \\
 3 & 4 & 1 & 2 & 5 & 6
\end{pmatrix}.
\]
Any element in the centralizer of $h$ must be of the form
\[
\begin{pmatrix}
 1 & 2 & 3 & 4 & 5 & 6 \\
 * & * & * & * & 5 & 6
\end{pmatrix}
\quad \text{or} \quad 
\begin{pmatrix}
 1 & 2 & 3 & 4 & 5 & 6 \\
 * & * & * & * & 6 & 5
\end{pmatrix}
\]
It is therefore sufficient to perform all calculations in $S_4$. If an element of $S_4$ is a composition of an even (odd) number of transpositions, the corresponding element of $Z_{A_6}(h)$ is given by completing it with the identity map (transposition map) on the fifth and sixth letter. 

Let 
\[
 g_1= 
\begin{pmatrix}
 1 & 2 & 3 & 4 \\
 2 & 1 & 4 & 3 
\end{pmatrix},
g_2= 
\begin{pmatrix}
 1 & 2 & 3 & 4 \\
 3 & 2 & 1 & 4 
\end{pmatrix},
g_3= 
\begin{pmatrix}
 1 & 2 & 3 & 4 \\
 1 & 4 & 3 & 2 
\end{pmatrix}.
\]
and check that $g_1,g_2,g_3 \in Z_{A_6}(h)$. Define $g_1 g_2 =:c$ and check 
\[
c= 
\begin{pmatrix}
 1 & 2 & 3 & 4 \\
 2 & 3 & 4 & 1
\end{pmatrix}, \quad
c^2 
= h.
\]
Now $g_3 c g_3 = c^3 $ and the subgroup of $S_4$ ($A_6$, respectively) generated by $c$ and $g_3$ is seen to be a dihedral group of order eight; $\langle g_3 \rangle \ltimes \langle c \rangle = D_8 < Z_{A_6}(h)$. In order to show equality, assume that $Z_{A_6}(h)$ is bigger. It then follows that the centralizer of $h$ in $S_4$ is a subgroup of order 12, in particular, it has a subgroup of order three. Going through the list of elements of order three in $S_4$ one checks that none commutes with $h$ and obtains a contradiction.
\end{proof}
Let $D_8=C_2\ltimes C_4$ where $C_2$ is generated by
$g= g_3$ and $C_4$ by $c$ and note that $c^2 =h$.
We study the action of $\tau $
on $D_8$ by conjugation. Since $C_4$ is the only cyclic subgroup
of order four in $D_8$, it is $\tau $-invariant. 
If $c$ is $\tau $-fixed, i.e. $\tau c = c \tau$, then 
\[
(\tau c)^2= c \tau \tau c = c \sigma c \overset{c \in Z(\sigma)}{=} \sigma c^2 =
\sigma h.
\]
In this case
$\tau c$ generates a cyclic group of order four acting
freely on $X$, a contradiction. 
So $\tau $ acts on $\langle c\rangle$
by $c \mapsto c^3$ and $c^2 \mapsto c^2$.
Now $\tau g\tau^{-1}=c^kg$ for some $k \in \{0,1,2,3\}$.
If $k=2$, then 
\[
(\tau g)^2= \tau g \tau g = \tau g \tau^{-1} \tau^2 g = c^2 g \sigma g
\overset{g \in Z(\sigma)}{=} c^2 \sigma = \sigma h
\]
and we obtain the same 
contradiction as above. So $k\in \{1,3\}$ and by choosing the 
appropriate generator of $\langle c\rangle$ we may assume that $k=3$.
The action of $\tau $
on $Z_{A_6}(h)=D_8$ given by $g\mapsto c^3 g$ and $c\mapsto c^3$. 
\begin{lemma}
$G'=\langle c\rangle$.
\end{lemma}
\begin{proof}
The commutator subgroup $G'$ is the smallest normal subgroup $N$ of $G$ such that $G/N$ is Abelian.
We use the above considerations about the action of $\tau$ on $D_8$ by conjugation.
The subgroup $\langle c \rangle$ is normal in $G = D_8 \rtimes \langle \tau \rangle$ and $G/ \langle c \rangle$ is seen to be Abelian. Since $G / \langle c^2 \rangle$ is not Abelian, $ G' \neq \langle c^2 \rangle$ and the lemma follows.
\end{proof}
\subsection{The group $H = G / \langle \sigma \rangle$}
We consider the quotient $Y=X/\sigma $ equipped with
the action of 
$G/\sigma =:H=Z_{\tilde A_6}(\sigma)/\langle \sigma \rangle=D_8\rtimes C_2$.  The
group $C_2$ is generated by $[\tau]_\sigma$.
For simplicity, we transfer the above notation from $G$ to $H$ by writing e.g. $ \tau$ for
$[\tau]_\sigma$. etc.. Since $\tau g\tau ^{-1}=c^3 g= g c$, it follows as above that 
$H'=\langle c\rangle$. 

Let $K < G$ be the cyclic group of order eight generated by $g \tau $.
\[
K = \{ \mathrm{id}, g \tau , c \sigma, g \tau^3 c, c^2, g \tau c^2,  \sigma
c^3, g c \tau^3  \}.
\]
We denote the image
of $K$ in $G/\sigma $ by the same symbol.
Since $[\sigma c ]_\sigma = [c]_\sigma \in K$ it contains $H'=\langle c\rangle$
and we can write
$H=\langle \tau \rangle \ltimes K=D_{16}$.  
\begin {lemma}\label{normal groups}
There is no nontrivial normal subgroup $N$ in $H$ with
$N\cap H'=\{\mathrm{id}\}$.
\end {lemma}
\begin {proof}
If such a group exists, first consider the case $N \cap K =
\{\mathrm{id}\}$. Then $N \cong C_2$ and 
$H= K \times N$ would be Abelian, a contradiction. If $N \cap K \neq 
\{\mathrm{id}\}$ then  $N \cap K = \langle (g \tau) ^k\rangle $ for some $k \in
\{1,2,4\}$. This implies $(g \tau )^4 =c^2 \in N$ and contradicts $N \cap H' = N
\cap \langle c \rangle = \emptyset$.
\end{proof}
The following observations strongly rely the assumption that $\sigma h$ acts freely on $X$.
\begin {lemma}\label{free on B}
The subgroup $H'$ acts freely on the branch set $B = \pi(\mathrm{Fix}_X(\sigma))$
in $Y$.
\end {lemma}
\begin {proof} 
If for some $b \in B$ the isotropy group $H'_b$ is nontrivial, then $c^2(b) = h(b)=b$ and
$\sigma h$ fixes the corresponding point $\tilde b\in X$.
\end {proof}
\begin {corollary}
The subgroup $H'$ acts freely on the set $\mathcal R$ 
of rational branch curves. In particular, the number of rational branch curves
$n$ is a 
multiple of four.
\end {corollary}
\begin {corollary}\label{tau-fixed}
The subgroup $H'$ acts freely on the set of $\tau $-fixed
points in $Y$.
\end {corollary}
\begin {proof}
We show $\mathrm{Fix}_Y(\tau) \subset B$.
Since $\sigma = \tau ^2$ on $X$, a $\langle \tau \rangle$-orbit $\{x, \tau x,
\sigma x, \tau^3 x \}$ in $X$ gives rise to a $\tau$-fixed point $y$ in the
quotient $Y = X /  \sigma$ if and only if $\sigma x = \tau x $. Therefore,
$\tau$-fixed points in $Y$ correspond to $\tau$-fixed points in $X$. By
definition $\mathrm{Fix}_X(\tau) \subset \mathrm{Fix}_X(\sigma)$ and the claim
follows. 
\end {proof}
\section {$H$-minimal models of $Y$}\label {reduction to del Pezzo}
Since $\mathrm{Fix}_X(\sigma) \neq \emptyset$, the quotient surface $Y$ is a smooth rational $H$-surface to which we apply the equivariant minimal model program. We denote by $Y_\mathrm{min}$ an $H$-minimal model of $Y$. It is known that $Y_\mathrm{min}$ is either a Del Pezzo surface or an $H$-equivariant conic bundle over $\mathbb P_1$.
\begin {theorem}\label {no equivariant fibration}
An $H$-minimal model $Y_{\mathrm {min}}$ does not admit an $H$-equivariant $\mathbb P_1$-fibration. In particular, $Y_{\mathrm {min}}$ is a Del Pezzo surface.
\end {theorem}
In order to prove this statement we begin with the following general fact (cf. Proof of Lemma \ref{fixed points on mori fibers}).
\begin {lemma}\label{no increase}
If $Y\to Y_{\mathrm {min}}$ is an $H$-equivariant Mori reduction and $A$ a cyclic subgroup of $H$, then 
\[
\vert \mathrm {Fix}_Y(A)\vert \geq 
\vert \mathrm {Fix}_{Y_{\mathrm {min}}}(A)\vert \,.
\]
\end {lemma} 
\begin {proof}
Each step of a Mori reduction is known to contract a disjoint union of (-1)-curves. 
It is sufficient to prove the statement for one step
in a Mori reduction. If such a step 
changes the
number of fixed points, then some Mori fiber $E$ of the reduction
is contracted to an $A$-fixed point.
The rational curve
$E$ is $A$-invariant and therefore contains two $A$-fixed points. The number of fixed points drops.
\end {proof}
Suppose that some $Y_\mathrm{min}$ is an $H$-equivariant conic bundle, i.e., there is an $H$-equivariant fibration 
$p :Y_{\mathrm {min}}\to \mathbb P_1$ with generic fiber $\mathbb P _1$. We
let
$p _*:H\to \mathrm {Aut}(\mathbb P_1)$ 
be the associated homomorphism.
\begin {lemma}\label{ker p*}
$\mathrm {Ker}(p_*)\cap H'=\{\mathrm {id}\}\,.$
\end {lemma}
\begin {proof}
The elements of $\mathrm {Ker}(p_*)$ fix two points
in every generic $p $-fiber. If $h = c^2 \in H' = \langle c \rangle$ fixes
points in every generic $p$-fiber, then $h$ acts trivially on a one-dimensional subset $C \subset Y$. Since $h=c^2$ acts symplectically on $X$ it has only
isolated fixed points in $X$. Therefore, on the preimage $\tilde C = \pi^{-1}(C) \subset X$, the action of $h$ coincides with the action of $\sigma$. But then $\sigma h | _{\tilde C} = \mathrm{id}| _{\tilde C}$ contradicts the assumption that $\sigma h$ acts freely on $X$.
\end {proof}
\begin{proof}[Proof of Theorem \ref{no equivariant fibration}]
Since there are no nontrivial normal
subgroups in $H$ which have trivial intersection with 
$H'$ (Lemma \ref{normal groups}), it follows from Lemma \ref{ker p*} that $\mathrm {Ker}(p_*)=
\{ \mathrm{id} \}$, i.e., the group $H$
acts effectively on the base. 

We regard $H$ as the semidirect product 
$H=\langle \tau \rangle \ltimes K$, where $K=C_8$
is described above. The group $H$ acts on the base as a dihedral group and therefore $\tau$ exchanges the $K$-fixed points. We will obtain a contraction by
analyzing the $K$-actions on the fibers over its two fixed
points.  Since $\tau $ exchanges these fibers, it is
enough to study the $K$-action on one of them which
we denote by $F$.

By Lemma \ref{singular fibers of conic bundle} there are two situations which we must consider:
\begin{enumerate}
 \item 
 $F$ is a regular fiber of $Y_{\mathrm {min}}\to \mathbb P_1$.
 \item
 $F=F_1\cup F_2$ is the union of two (-1)-curves intersecting transversally in
one point. 
\end{enumerate}

We study the fixed points
of $c$, $h=c^2$ and $g \tau $ in $Y_{\mathrm {min}}$. Recall that
in $X$ the symplectic transformation $c$ has precisely four fixed
points and $h$ has precisely eight fixed points. This set of eight
points is stabilized by the full centralizer of $h$, in particular by
$K = \langle  g \tau  \rangle \cong C_8$.

Since $h \sigma$ acts by assumption freely on $X$, it follows that
$\sigma$ acts freely on the set of $h$-fixed points in $X$.  If $hy=y$ for some $y \in Y$, then the preimage of $y$ in $X$ consists of two elements $x_1,\sigma x_1=x_2$. If these form an $\langle h \rangle $-orbit, then both are $ \sigma h $-fixed, a contradiction. It follows that $\{x_1,x_2 \} \subset \mathrm{Fix}_X(h)$  and the number of $h$-fixed points in $Y$ is precisely four. In particular, $h$ acts effectively on any curve in $Y$. 

Let us first consider 
Case 2 where $F= F_1 \cup F_2$ is reducible. 
Since $\langle c\rangle $ is the only subgroup
of index two in $K$, it follows that $\langle c \rangle $ stabilizes $F_i$ and both $c$ and $h$ have
three fixed points in $F$ (two on each irreducible component, one is the point of intersection $F_1 \cap F_2$), i.e., six fixed points on $F \cup \tau F \subset Y_\mathrm{min}$. This is contrary to Lemma \ref{no increase}
 because $h$ has at most four fixed
points in $Y_{\mathrm {min}}$.

If $F$ is regular (Case 1),
then the cyclic group $K$ has two fixed points on the rational curve
$F$. Since $h \in K$, the four $K$-fixed points on $F \cup \tau F$ are
contained in the set of $h$-fixed points on $Y_\mathrm{min}$. As
$|\mathrm{Fix}_{Y_\mathrm{min}}(h)| \leq 4$, the $K$-fixed points
coincide with the four $h$-fixed points in $Y_\mathrm{min}$;
\[
\mathrm{Fix}_{Y_\mathrm{min}}(h)= \mathrm{Fix}_{Y_\mathrm{min}}(K). 
\]
In particular, the Mori reduction does not affect the four $h$-fixed
points $\{y_1, \dots y_4\}$ in $Y$. By equivariance of the reduction,
the group $K$ acts trivially on this set of four points. Passing to
the double cover $X$, we conclude that the action of $g \tau \in K$ on
a preimage $\{x_i , \sigma x_i\}$ of $y_i$ is either trivial or
coincides with the action of $\sigma$. In both cases it follows that
$(g \tau )^2 = c\sigma$ acts trivially on the set of $h$-fixed points in $X$. As $\mathrm{Fix}_X(c) \subset \mathrm{Fix}_X(h)$, this is contrary to the fact that $\sigma$ acts freely on $\mathrm{Fix}_X(h)$.
\end{proof}
In the following we wish to identify the Del Pezzo surface $Y_\mathrm{min}$. For thus, we use the Euler characteristic
formulas,
\[
24= e(X)= 2 e(Y)-2n + \underset{\text{if $D_g$ is present}}{\underbrace{2g-2}},  
\] 
where $D_g \subset B$ is of general type, $g = g(D_g) \geq 2$, and 
\[
e(Y)= e(Y_{\mathrm{min}}) + m, 
\]
where $m = | \mathcal E|$ denotes the total number of Mori fibers.
For convenience we introduce the difference
$\delta =m -n$. 
If a branch curve $D_g$ of general type is present, then
$
13-g-\delta =e(Y_{\mathrm {min}})
$
and if it is not present
$
12-\delta =e(Y_{\mathrm {min}})
$.
\begin {proposition}
For every Mori fiber $E$ the orbit $H.E$ consists
of at least four Mori fibers.
\end {proposition}
\begin {proof}
We need to distinguish three cases:
\[
\text{1.)}\,\, E\cap B \neq \emptyset\text{\ and\ } E \not\subset B; \quad\quad
\text{2.)}\,\,  E \subset B; \quad\quad
\text{3.)}\,\,  E \cap B = \emptyset
\]
{\bf Case 1 }
Since $H'$ acts freely on the branch curves
and $E$ meets $B$ in at most two points,
we know $\vert H'.E\vert \ge 2$.
If $\vert H.E\vert =2$, then the isotropy group
$H_E$ is a normal subgroup of index two which necessarily
contains the commutator group $H'$, a contradiction.

{\bf Case 2 }
We show that the $H'$-orbit of $E$ consists of four Mori fibers.
If it consisted of less than four Mori fibers, the stabilizer $H'_E \neq \{\mathrm{id}\}$ of $E$ in $H'$ would fix two points in $E \subset B$. This contradicts Lemma \ref{free on B}.

{\bf Case 3 }
 All Mori fibers disjoint from $B$ have self-intersection (-2) and meet exactly one Mori fiber of the previous steps of the reduction in exactly one point. 
If $E \cap B = \emptyset$ there is a chain of Mori fibers $E_1, \dots, E_k =E$ connecting $E$ and $B$ with the following properties:
The Mori fiber $E_1$ is the only one to have nonempty intersection with $B$ and is the first curve of this configuration to be blown down in the reduction process. The curves fulfil $(E_i, E_{i+1})=1$ for all $i \in \{1, \dots, k-1\}$ and  $(E_i,E_j)=0$ for all $j\neq i+1$. The curves are blown down subsequently and meet no Mori fibers outside this chain.

The $H$-orbit of this union of Mori fibers consists of at least four copies of this chain. This is due to that fact that the $H$-orbit of $E_1$ consists of at least four Mori fibers by Case 1. In particular, the $H$-orbit of $E$ consists of at least four copies of $E$.
\end{proof}
\begin {corollary}
The difference $\delta $ is a non-negative multiple $4k$
of four. If $\delta =0$, then $X$ is a double cover of $Y = Y_\mathrm{min} = \mathbb P_1 \times \mathbb P_1$ branched along a curve of genus nine. 
\end {corollary}
\begin {proof}
Above we have shown that $m$ and 
and $n$ are multiples of four. Therefore $\delta =4k$.

If $\delta$ was negative, i.e., 
$m < n$, there is no configuration of Mori fibers meeting the rational branch curves such that the corresponding contractions transform the (-4)-curves in $Y$ to curves on a Del Pezzo surface $Y_\mathrm{min}$. It follows that $\delta$ is non-negative.

If $\delta= 0$, then $n= m=0$ and $Y$ is $H$-minimal. The commutator subgroup $H' \cong C_4$ acts freely on the branch locus $B$ implying $e(B)\in \{0,-8,-16, \dots \}$. Since the Euler characteristic of the Del Pezzo surface $Y$ is at least 3 and at most 11,
\[
6 \leq 2e(Y)= 24 + e(B) \leq 22,
\]
we only need to consider the case $e(Y)\in \{4,8\}$ and $B=D_{g}$ for $g \in \{9,5\}$. 

The automorphism group of a Del Pezzo surface of degree 4 is  $C_2^4 \rtimes \Gamma$ for $\Gamma \in \{C_2,C_4, S_3, D_{10} \}$.  If $D_{16} < C_2^4 \rtimes \Gamma$ then  $ A := D_{16} \cap C_2^4 \lhd D_{16}$ and $A$ is either trivial or isomorphic to $C_2$. In both case $D_{16} / A$ is not a subgroup of $\Gamma$ in any of the cases listed above. Therefore, $e(Y) \neq 8$.

A Del Pezzo surface of degree 8 is either the blow-up of $ \mathbb P_2$ in one point or $\mathbb P_1 \times \mathbb P_1$. Since the first is never equivariantly minimal, it follows that $Y \cong \mathbb P_1 \times \mathbb P_1$ and $g(B)=9$.
\end {proof}
\begin {theorem}
Any $H$-minimal model $Y_\mathrm{min}$ of $Y$ is 
$\mathbb P_1\times \mathbb P_1$ .
\end {theorem}

\begin {proof}
Suppose $\delta \neq 0$.
Since $\delta \ge 4$, it follows that $e(Y_\mathrm{min})=13-g-\delta \le 7$ if a branch curve $D_g$ of general type is present, and $e(Y_\mathrm{min})=12-\delta \le 8$
if not. We go through the list of
of Del Pezzo surfaces with $e(Y_\mathrm{min}) \leq 8$.
\begin {itemize}
\item
If $e(Y_{\mathrm {min}})=8$, i.e., $\mathrm{deg}(Y_\mathrm{min}) = 4$, then the possible automorphism groups are very limited and we have alredy noted above that
$D_{16}$ does not occur. 
\item
If $e(Y_{\mathrm {min}})=7$, then 
$\mathrm {Aut}(Y_{\mathrm {min}})=S_5$.  Since $120$ is not
divisible by $16$, we see that a Del Pezzo surface of degree five does not admit an effective action of the group $H$.
\item
If $e(Y_{\mathrm {min}})=6$, then 
$A:=\mathrm {Aut}(Y_{\mathrm {min}})=
(\mathbb C^*)^2\rtimes (S_3\times C_2)$.
We denote by $A^\circ \cong (\mathbb C^*)^2$ the connected component of $A$.
If $q :A\to A/A^\circ$ is the canonical quotient homomorphism
then $q (H') <  q(A)'\cong C_3$. Consequently 
$H'=C_4 < A^\circ$. We may realize $Y_{\mathrm {min}}$ as $\mathbb P_2$ blown up at the three
corner points and $A^\circ$ as the space of diagonal
matrices in $\mathrm {SL}_3(\mathbb C)$.
Every possible representation of $C_4$ in this group
has ineffectivity along one of the lines joining corner points.
But, as we have seen before, the elements of $H'$, in particular $c^2 = h$, have only isolated fixed points
in $Y_{\mathrm {min}}$.  
\item
A Del Pezzo surface obtained by blowing up one or two points in $\mathbb P_2$ is never $H$-minimal and therefore does not occur
\item 
Finally, $Y_{\mathrm {min}}\not=\mathbb P_2$:
If $e(Y_\mathrm{min}) =3$ then either $\delta =9$ (if $D_g$ is not present), a contradiction to $\delta = 4k$, or $g + \delta =10$. In the later case, $\delta =4,8$ forces $g= 6,2$. In both cases, the Euler characteristic $2-2g$ of $D_g$ is not divisible by 4. This contradicts the fact that $H'$ acts freely on $D_g$. 
\end {itemize}
We have hereby excluded all possible Del Pezzo surfaces except $\mathbb P_1\times \mathbb P_1$ and the proposition follows. 
\end {proof}
\section{Branch curves and Mori fibers}
We let $M : Y \to Y_\mathrm{min} = \mathbb P_1 \times \mathbb P_1$ denote an $H$-equivariant Mori reduction of $Y$. 
\begin{lemma}
The length of an orbit of Mori fibers is at least eight.
\end{lemma}
\begin{proof}
Consider the action of $H$ on $\mathbb P_1 \times \mathbb P_1$. 
Both canonical projections are equivariant with respect to the commutator subgroup $H'= \langle c \rangle  \cong C_4$. Since $c^2 \in H'$ does not act trivially on any curve in $Y$ or $Y_\mathrm{min}$, it follows that $H'$ has precisely four fixed points in $Y_\mathrm{min} =\mathbb P_1 \times \mathbb P_1$. 
Since $h = c^2$ has precisely four fixed points in $Y$ and $\mathrm{Fix}_Y (H') = \mathrm{Fix}_Y (c) \subset  \mathrm{Fix}_Y (c^2)$, we conclude that $H'$ has precisely four fixed points in $Y$ and it follows that the Mori fibers do not pass through $H'$-fixed points. Note that the $H'$-fixed points in $Y$ coincide with the $h$-fixed points. 

Suppose there is an $H$-orbit $H.E$ of Mori fibers of length strictly less then eight and let $p = M(E)$. We obtain an $H$-orbit $H.p$ in $\mathbb P_1 \times \mathbb P_1$ with $|H.p| \leq 4$. Now $| K.p| \leq 4$ implies that $K_p \neq \{\mathrm{id}\}$, in particular, $h= c^2 \in K_p$. It follows that $p$ is a $h$-fixed point. This contradicts the fact that the Mori fibers do not pass through fixed points of $h$.
\end{proof}
\begin{corollary}
The total number $m$ of Mori fibers equals 0, 8, or 16..
\end{corollary}
\begin{proof}
A total number of 24 or more Mori fibers would require 16 rational curves in $B$. This contradicts the bound for the number of connected components of the fixed point set of an antisymplectic involution on a K3-surface (cf. Corollary \ref{atmostten}) 
\end{proof}
Recalling that the number of rational branch curves is a multiple of four, i.e., $n \in \{0,4,8\}$ and  
using the fact $m \in \{0,8,16\}$ along with $m \leq n+9$, we conclude that the surface $Y$ is of one of the following types.
\begin{enumerate}
 \item
$m=0$\\ 
The quotient surface $Y$ is $H$-minimal. The map $X \to Y \cong \mathbb P_1 \times \mathbb P_1$ is branched along a single curve $B$. This curve $B$ is a smooth $H$-invariant curve of bidegree $(4,4)$.
\item
$m=8$ and $e(Y) = 12$\\
The surface $Y$ is the blow-up of $\mathbb P_1 \times \mathbb P_1$ in an $H$-orbit consisting of eight points. 
\begin{enumerate}
 \item 
If the branch locus $B$ of $X \to Y$ contains no rational curves, then $e(B)=0$ and $B$ is either an elliptic curve or the union of two elliptic curves defining an elliptic fibration on $X$. 
\item
If the branch locus $B$ of $X \to Y$ contains rational curves, their number is exactly four (Observe that eight or more rational branch curves of self-intersection (-4) cannot be modified sufficiently and mapped to curves on a Del Pezzo surface by contracting eight Mori fibers). It follows that the branch locus is the disjoint union of an invariant curve of higher genus and four rational curves. 
\end{enumerate}
\item
$m=16$ and $e(Y) =20$\\
The map $X \to Y$ is branched along eight disjoint rational curves. 
\end{enumerate}
We can simplify the above situation by studying rational curves in $B$, their intersection with Mori fibers and their images in $\mathbb P_1 \times \mathbb P_1$.
\begin{proposition}
If $e(Y)=12$, then $n=0$.
\end{proposition}
\begin{proof}
Suppose $n \neq 0$ and let $C_i \subset Y$ be a rational branch curve. Since $C_i^2 =-4$ and $M(C_i) \subset \mathbb P_1 \times \mathbb P_1$ has self-intersection $\geq 0$ it must meet the union of Mori fibers $\bigcup E_j$. 
All possible configurations of Mori fibers yield image curves $M(C_i)$ of self-inter\-sec\-tion $\leq 4$. If $M(C_i)$ is a curve a bidegree $(a,b)$, then, by adjunction.
\[
 2g(M(C_i)) -2 = (M(C_i))^2 + (M(C_i) \cdot K_{\mathbb P_1 \times \mathbb P_1} )= 2ab -2a-2b,
\]
and $(M(C_i))^2 = 2ab \leq 4$ implies 
 that $g(M(C_i))=0$. In particular, $M(C_i)$ must be nonsingular. Hence each Mori fiber meets $C_i$ in at most one point. It follows that $C_i$ meets four Mori fibers, each in one point, and $(M(C_i))^2 =0$. 
Curves of self-intersection zero in $\mathbb P_1 \times \mathbb P_1$ are fibers of the canonical projections $\mathbb P_1 \times \mathbb P_1 \to \mathbb P_1$. 
The curve $C_1$ meets four Mori fibers $E_1, \dots E_4$ and each of these Mori fibers meets some $C_i$ for $i \neq 1$. After renumbering, we may assume that $E_1$ and $E_2$ meet $C_2$ and therefore $M(C_1)$ and $M(C_2)$ meet in more than one point, a contradiction. It follows that $e(Y) = 12$ implies $n=0$
\end{proof}
\begin{proposition}
If $e(Y)=20$, then $Y$ is the blow-up of $\mathbb P_1 \times \mathbb P_1$ in sixteen points
\[
 \{p_1, \dots p_{16}\} = (F_1 \cup F_2 \cup F_3 \cup F_4) \cap  (F_5 \cup F_6 \cup F_7 \cup F_8),
\]
 where $F_1, \dots F_4$ are fibers of the canonical projection $\pi_1$ and $F_5, \dots F_8$ are fibers of $\pi_2$. The branch locus is given by the proper transform of $\bigcup F_i$ in $Y$.
\end{proposition}
\begin{proof}
We denote the eight rational branch curves by $C_1, \dots C_8$. The Mori reduction can have two steps. A slightly more involved study of possible configurations of Mori fibers shows that $0 \leq (M(C_i))^2 \leq 4$. 
As above $M(C_i)$ is seen to be nonsingular and each Mori fiber can meet $C_i$ in at most one point. Any configuration of curves with this property yields $(M(C_i))^2=0$ and $F_i = M(C_i)$ is a fiber of a canonical projection $\mathbb P_1 \times \mathbb P_1 \to \mathbb P_1$. 

If there are Mori fibers disjoint from $B$ these are blown down in the second step of the Mori reduction. Let $E_1, \dots, E_8$ denote the Mori fibers of the first step and $\tilde E_1, \dots, \tilde E_8$ those of the second step. We label them such that $\tilde E_i$ meets $E_i$. Each curve $E_i$ meets two rational branch curves $C_i$ and $C_{i+4}$ and their images $F_i = M(C_i)$ and $F_{i+4}=M(C_{i+4})$ meet with multiplicity $\geq 2$. This is contrary to the fact that they are fibers of the canonical projections. It follows that there are no Mori fibers disjoint from $B$ and all 16 Mori fibers are contrancted simultaniously. There is precisely one possible configuration of Mori fibers on $Y$ such that all rational brach curves are mapped to fibers of the canonical projections of $\mathbb P_1 \times \mathbb P_1$: The curves $C_1, \dots C_4$ are mapped to fibers of $\pi_1$ and $C_5, \dots, C_8$ are mapped to fibers of $\pi_2$. The Mori reduction contracts 16 curves to the 16 points of intersection $\{p_1, \dots p_{16} \} = (\bigcup_{i=1}^4 F_i) \cap(\bigcup_{i=5}^8 F_i)  \subset \mathbb P_1 \times \mathbb P_1$. 
\end{proof}
Let us now restrict our attention to the case where the branch locus $B$ is the union of two linearly equivalent elliptic curves and exclude this case.
\subsection{Two elliptic branch curves}
In this section we prove:
\begin{theorem}\label{two elliptic branch curves}
$\mathrm{Fix}_X(\sigma)$ is not the union of two elliptic curves.
\end{theorem}
We assume the contrary, let  $\mathrm{Fix}_X(\sigma) = D_1 \cup D_2$ with $D_i$ elliptic and let 
$f :X\to \mathbb P_1$ denote the elliptic fibration
defined by the curves $D_1$ and $D_2$.
Recall that $\sigma $ acts effectively on the base $\mathbb P_1$ as otherwise $\sigma$ would act trivially in a neighbourhood of $D_i$ by a linearization argument (cf. Theorem \ref{FixSigma}).
It follows that the group of order
four generated by $\tau $ acts effectively on $\mathbb P_1$.

Let $I$ be the ineffectivity of the induced $G$-action on the base
$\mathbb P_1$. We regard
$G=C_4\ltimes D_8$ where $C_4=\langle \tau \rangle$ and $D_8$ is the centralizer of $\sigma$ in $A_6$ (cf. Section \ref{centralizer}) and define $J:=I\cap D_8$. First, note that $I$ is nontrivial:
\begin {lemma}
The group $G$ does not act effectively on $\mathbb P_1$, i.e., $I \neq \{ \mathrm{id}\}$. 
\end {lemma}
\begin {proof}
If $G$ acts effectively on $\mathbb P_1$, then $G$ is among the groups specified in Remark \ref{autP_1}. In our special case $|G| = 32$ and $G$ would have to be cyclic or dihedral. 
Since the group $G$ does not contain a cyclic group of order 16, this is a contradiction.
\end {proof}
\begin {lemma}
The intersection $J=I\cap D_8$ is nontrivial.
\end {lemma}
\begin{proof}
Assume the contrary and let $J = I \cap D_8 = \{e\}$. We consider the quotient $G \to G/D_8 \cong C_4$  and see that either $I \cong C_2$ or $I \cong C_4$.
\begin{itemize}
\item 
If $I\cap D_8=\{e\}$ and $I\cong C_2$, we write $I=\langle \sigma \xi\rangle $ with $\xi \in D_8$
an element of order two. Now $I$ is normal if and only if $\xi =h$, i.e., $I=\langle \sigma h\rangle$.
In this case, since $\sigma h \notin K$, the image of $K$ in $G/I$ is a normal subgroup of
index two and one checks that $G/I\cong D_{16}$ . The group $K$ is mapped
injectively into $G/I$.  The equivalence
relation defining this quotient identifies $\sigma $ and $h$ and
both are in the image of $K$. So $h$-fixed points in $X$ must lie in the fibers over the $\sigma$-fixed
points in $\mathbb P_1$, i.e., the $\sigma $-fixed points sets $D_1, D_2$.
Since $\sigma h$ acts freely on $X$, this is a contradiction.
\item 
If $I\cap D_8=\{e\}$ and $I\cong C_4$ we write
$I =\langle \tau \xi\rangle$ and 
show that for no choice of $\xi$ the group $I =\langle \tau \xi\rangle$ is normal in $G$:
If $\xi =c^kg$, then $\langle \tau \xi\rangle =K$ is of order eight.
If $\xi =c^k$, then $\langle \tau \xi\rangle$ is of order
four and has trivial intersection with $D_8$. It is however not
normalized by $g$.
\end{itemize}
As we obtain contradictions in all cases, we see that the intersection $J=I\cap D_8$ is nontrivial.
\end{proof}
In the following, we consider the different possibilities for the order of $J$ and show that in fact none of these occur.

If $|J|=8$ then $D_8 \subset I$. Recall that any automorphism group of an elliptic curve splits into an Abelian part acting freely and a cyclic part fixing a point. Since $D_8$ is not Abelian, any $D_8$-action on the fibers of $f$ must have points with nontrivial isotropy and gives rise to a positive-dimensional fixed point set of some subgroup of $D_8$ on $X$ contradicting the fact that $D_8$ acts symplectically on $X$. 
It follows that the maximal possible order of $J$ is four.
\begin{lemma}\label{I does not contain c}
 The ineffectivity $I$ does not contain $\langle c\rangle$.
\end{lemma}
\begin{proof}
Assume the contrary and consider the fixed points of $c^2$. If a $c^2$-fixed point lies at a smooth point of a fiber of $f$, then the linearization of the $c^2$-action at this fixed point gives rise to a positive-dimensional fixed
point set in $X$ and yields a contradiction.
It follows that
the fixed points of $c^2$ are contained
in the singular $f $-fibers.
Since $\langle \tau \rangle $ normalizes $\langle c\rangle$ and the $\langle \tau \rangle $-orbit of a singular
fiber consists of four such fibers, we must only consider two cases:
\begin {enumerate}
\item
The eight $c^2$-fixed points are contained in four singular
fibers (one $\langle \tau \rangle $-orbit of fibers), each of these fibers contains two $c^2$-fixed points.
\item
The eight $c^2$-fixed points are contained in eight singular fibers
(two $\langle \tau \rangle$-orbits).
\end {enumerate}  
Note that $\langle c^2 \rangle$ is normal in $I$
 and therefore $I$ acts on the set of 
$\langle c^2\rangle $-fixed points. 
In the second case, all eight $c^2$-fixed points are
also $c$-fixed. This is contrary to $c$ having only four fixed
points and therefore the second case does not occur.

The first case does not occur for similar reasons: If
$c^2$ has exactly two fixed points $x_1$ and $x_2$ in some 
fiber $F$, then $\langle c\rangle $ either acts transitively
on $\{x_1,x_2\}$ or fixes both points. Since $\mathrm{Fix}_X(c) \subset \mathrm{Fix}_X(c^2)$ and $\langle c\rangle $
must have exactly one fixed point on $F$, this is impossible.
\end {proof}
\begin {corollary}
$|J| \neq 4$.
\end {corollary}
\begin {proof}
Assume $|J| = 4$. Using $\tau$ we check that no subgroup of $D_8$ isomorphic to $C_2 \times C_2$ is  normal in $G$. It follows that the group $\langle c\rangle$ is the only order four subgroup of $D_8$ which is normal in $G$ and therefore $J = \langle c\rangle$. By the lemma above this is however impossible.
\end {proof}
It remains to consider the case where $|J|=2$. 
The only normal subgroup
of order two in $D_8$ is $J=\langle h\rangle$.
\begin {lemma}
If $|J|=2$, then $I=\langle \sigma c\rangle$.
\end {lemma}
\begin {proof}
We first show that $|J|=2$ implies $|I|=4$:
If  $|I|=2$, then $ I =  \langle h \rangle$ and $G / I = C_4 \ltimes (C_2\times C_2)$. Since this group does not act effectively on $\mathbb P_1$, this is a contradiction. 
If  $|I| \geq 8$, then $G/I$ is Abelian and therefore $I$ contains the commutator subgroup $G' = \langle c \rangle$. This contradicts Lemma \ref{I does not contain c}.
It follows that $|I|=4$ and either $I \cong C_4$ or $I \cong C_2 \times C_2$. In the later case, the only possible choice is $I = \langle \sigma \rangle \times \langle h \rangle$ which contradicts the fact that $ \sigma$ acts effectively on the base. 
It follows that
$I=\langle \sigma \xi \rangle$, where $\xi ^2=h$
and therefore $\xi =c$.
\end {proof}
Let us now consider the action of $G$ on $X$ with
$I=\langle \sigma c\rangle $.
Recall that 
the cyclic group $\langle \tau \rangle$ acts effectively on the base 
and has two fixed points there. Since $\sigma =\tau ^2$, these
are precisely the two $\sigma $-fixed points. In particular,
$\langle \tau \rangle$ stabilizes both $\sigma $-fixed point 
curves $D_1$ and $D_2$ in $X$.  
Furthermore, the transformations $\sigma c$ and $c$ stabilize $D_i$ for $i =1,2$. 
Since the only fixed points of $c$ in $\mathbb P_1$ are the
images of $D_1$ and $D_2$,
$$
\mathrm {Fix}_X(c)\subset D_1\cup D_2=\mathrm {Fix}_X(\sigma).
$$
 On the other hand, we know
that $\mathrm {Fix}_X(c)\cap \mathrm {Fix}_X(\sigma )=\emptyset$. 
Thus $I=\langle \sigma c\rangle $ is not possible
and the case $|J|=2$ does not occur.

We have hereby eleminated all possibilities for $|J|$ and completed the proof of Theorem \ref{two elliptic branch curves}.

\section{Rough classification of $X$}
We summerize the observations of the previous section in the following classification result.
\begin{theorem}\label{roughclassiA6}
Let $X$ be a K3-surface with an effective action of the group $G$ such that $\mathrm{Fix}_X(h\sigma) = \emptyset$. Then $X$ is one of the following types:
\begin{enumerate}
 \item 
 a double cover of $\mathbb P_1 \times \mathbb P_1$ branched along a smooth $H$-invariant curve of bidegree (4,4). 
\item
 a double cover of a blow-up of $\mathbb P_1 \times \mathbb P_1$ in eight points and branched along a smooth elliptic curve $B$. The image of $B$ in $\mathbb P_1 \times \mathbb P_1$ has bidegree (4,4) and eight singular points.
\item 
 a double cover of a blow-up $Y$ of $\mathbb P_1 \times \mathbb P_1$ in sixteen points $\{p_1, \dots p_{16}\} = (\bigcup_{i=1}^4 F_i) \cap(\bigcup_{i=5}^8 F_i)$, where $F_1, \dots F_4$ are fibers of the canonical projection $\pi_1$ and $F_5, \dots F_8$ are fibers of $\pi_2$. The branch locus ist given by the proper transform of $\bigcup F_i$ in $Y$. The set $\bigcup F_i$ is an invariant reducible subvariety of bidegree (4,4).
\end{enumerate}
\end{theorem}
\begin{proof}
It remains to consider case 2. and show that the image of $B$ in $\mathbb P_1 \times \mathbb P_1$ has bidegree (4,4) and eight singular points.
We prove that each Mori fiber $E$ meets the branch locus $B$ either in two points or once with multiplicity two, i.e., we need to check that $E$ may not meet $B$ transversally in exactly one point.
If this was the case, the image $M(B)$ of the branch curve is a smooth $H$-invariant curve of bidegree $(2,2)$. 
The double cover $X'$ of $\mathbb P_1 \times \mathbb P_1$ branched along the smooth curve $M(B)= C_{(2,2)}$ is a smooth surface. Since $X$ is K3 and therefore minimal the induced birational map $X \to X'$ is an isomorphism. This is a contradiction since $X'$ is not a K3-surface.

As each Mori fiber meets $B$ with multiplicity two, the self-intersection number of $M(B)$ is 32 and $M(B)$ is a curve of bidegree (4,4) with eight singular points. These singularities are either nodes or cusps depending on the kind of intersection of $E$ and $B$. We obtain a diagram 
\[
\begin{xymatrix}{
X_\mathrm {sing}\ar[d]^{2:1} & X \ar[d]^{2:1} \ar[l]^{\text{desing.}}\\
C_{(4,4)} \subset \mathbb P_1 \times \mathbb P_1 & \ar[l]_>>>>>>{M} Y \supset B}
\end{xymatrix}
\]
\end{proof}
In order to obtain a description of possible branch curves, we study the action of $H$ on $\mathbb P_1 \times \mathbb P_1$ and its invariants.
\subsection{The action of $H$ on $\mathbb P_1 \times \mathbb P_1$}
Recall that we consider the dihedral group $H \cong D_{16}$ generated by $\tau g$ of order eight and $\tau$. For convenience, we recall the group structure of $H$:
\begin {align*}
c= ( g \tau )^2, & \quad \tau g \tau = gc,\\
g^2 = \mathrm{id}, & \quad \tau c \tau  = c^3,\\
c^4 = \mathrm{id}, & \quad \tau ^2 = \mathrm{id}. 
\end {align*}
In this section, we prove:
\begin{proposition}
In appropriately chosen coordinates the action of $H$ on $\mathbb P_1\times \mathbb P_1$ given by
\begin {itemize}
\item
$
c([z_0:z_1],[w_0:w_1])=
([iz_0:z_1],[-iw_0:w_1])
$
\item
$
\tau ([z_0:z_1],[w_0:w_1])=
([z_1:z_0],[iw_1:w_0])
$
\item
$
g([z_0:z_1],[w_0:w_1])=
([w_0:w_1],[z_0:z_1])\,.
$
\end {itemize}
\end{proposition}
\begin{proof}[Sketch of proof]
 First note that the index two subgroup $H_1$ of $H$ preserving the canonical projections is generated by $\tau$ and $c$, i.e,
$H_1 = \langle \tau \rangle \ltimes \langle c \rangle \cong D_8$. We begin by choosing coordinates such that
$$
c([z_0:z_1],[w_0:w_1])=
([\chi _1(c)z_0:z_1],[\chi _2(c)w_0:w_1])
$$
where $\chi _i : H' \to S^1$ are faithful characters. Since $\tau$ acts transitively on the set of $H'$-fixed points, we conclude that after an appropriate change of coordinates not affecting the $H'$-action
$$
\tau ([z_0:z_1],[w_0:w_1])=
([z_1:z_0],[w_1:w_0]).
$$
The automorphism $g$ permutes the factors of $\mathbb P_1 \times \mathbb P_1$, stabilizes the fixed point set of $H'$ and fulfills $gcg^{-1}= c^3$ and $g \tau g^{-1} = c\tau$. Therefore, one finds
\begin {itemize}
\item
$
c([z_0:z_1],[w_0:w_1])=
([iz_0:z_1],[-iw_0:w_1])
$
\item
$
\tau ([z_0:z_1],[w_0:w_1])=
([z_1:z_0],[w_1:w_0])
$
\item
$
g([z_0:z_1],[w_0:w_1])=
([\lambda w_0:w_1],[\lambda ^{-1}z_0:z_1])\,,
$
where $\lambda ^2=i$.
\end {itemize} 
We introduce a change of coordinates such that $g$ is of the simple form
\[
g([z_0:z_1],[w_0:w_1]) = ([w_0:w_1],[z_0:z_1]).
\]
This does affect the shape of the $\tau$-action and yields the action of $H$ described in the propostion. 
\end{proof}
\subsection{Invariant curves of bidegree $(4,4)$}
Given the action of $H$ on $\mathbb P_1 \times \mathbb P_1$ discussed above, we wish to study the invariants and semi-invariants of bidegree $(4,4)$.
The space of $(a,b)$- bihomogeneous polynomials in $[z_0 : z_1][w_0 : w_1]$ is denoted by $\mathbb C_{(a,b)} ([z_0 : z_1][w_0 : w_1])$. 

 An invariant curve $C$ is given by a $D_{16}$-eigenvector $f \in \mathbb C_{(4,4)} ([z_0 : z_1][w_0 : w_1])$. The kernel of the $D_{16}$-representation on the line $\mathbb C f$ spanned $f$ contains the commutator subgroup $H' = \langle c \rangle $. It follows that $f$ is a linear combination of $c$-invariant monomials of bidegree $(4,4)$. These are
\begin{align*}
z_0^4w_0^4,\,
z_0^4w_1^4, \,
z_1^4w_0^4, \, 
z_1^4w_1^4, \,
z_0^2z_1^2w_0^2w_1^2, \,
z_0^3z_1w_0^3w_1, \,
z_0z_1^3w_0w_1^3.
\end{align*}
The polynomials
\begin{align*}
f_1 =z_0^4w_0^4 + z_1^4w_1^4, \quad 
f_2 =z_0^4w_1^4 + z_1^4w_0^4, \quad
f_3 =z_0^3z_1w_0^3w_1 -i z_0z_1^3w_0w_1^3
\end{align*}
span the space of $D_{16}$-invariants. Semi-invariants are appropiate linear combinations of  
\begin{align*}
g_1 =z_0^4w_0^4 - z_1^4w_1^4, \quad
g_2 =z_0^4w_1^4 - z_1^4w_0^4,\quad
g_3 =z_0^3z_1w_0^3w_1 +i  z_0z_1^3w_0w_1^3,\quad
g_4 =z_0^2z_1^2w_0^2w_1^2.
\end{align*}
Note
{\begin{displaymath}
\begin{array}{llll}
\tau (g_1) = -g_1, & \tau (g_2) = -g_2, & \tau (g_3) = -g_3, & \tau (g_4) = -g_4, \\
g(g_1) = g_1, & g(g_2) = -g_2, & g(g_3) = g_3, & g(g_4) = g_4.
              \end{array}
              \end{displaymath}
}
It follows that a $D_{16}$-invariant curve of bidegree $(4,4)$ in $\mathbb P_1 \times \mathbb P_1$ is of the following three types
\begin{align*}
C_a &= \{a_1 f_1 + a_2 f_2 + a_3 f_3 = 0\}, \\
C_b &= \{b_1 g_1 + b_3 g_3 + b_4 g_4 =0\}, \\
C_0 &= \{g_2 =0\}.
\end{align*}
\subsection{Refining the classification of $X$}
Using the above description of invariant curves of bidegree (4,4) we may refine Theorem \ref{roughclassiA6}.
\subsubsection{Reducible curves of bidegree $(4,4)$}
\begin{theorem}
Let $X$ be a K3-surface with an effective action of the group $G$ such that $\mathrm{Fix}_X(h\sigma) = \emptyset$. If $e(X/\sigma) = 20$, then $X/\sigma$ is equivariantly isomorphic to the blow up of $\mathbb P_1 \times \mathbb P_1$ in the singular points of the curve $C = \{f_1-f_2=0\}$ and $X \to Y$ is branched along the proper transform of $C$ in $Y$.
\end{theorem}
\begin{proof}
It follows from Theorem \ref{roughclassiA6} that $X$ is the double cover of $\mathbb P_1 \times \mathbb P_1$ blown up in sixteen points. These sixteen points are the points of intersection of eight fibers of $\mathbb P_1 \times \mathbb P_1$, four for each of fibration. 

By invariance these fibers lie over the base points $[1:1], [1:-1], [1: i], [1:-1]$ and the configurations of eight fibres is defined by the invariant polynomial $f_1-f_2$. 

The double cover $X \to Y$ is branched along the proper transform of this configuration of eight rational curves. This proper transform is a disjoint union of eight rational curves in $Y$, each with self-intersection (-4). 
\end{proof}
\subsubsection{Smooth curves of bidegree $(4,4)$}
\begin{theorem}
Let $X$ be a K3-surface with an effective action of the group $G$ such that $\mathrm{Fix}_X(h\sigma) = \emptyset$. If $X/\sigma \cong \mathbb P_1 \times \mathbb P_1$, then after a change of coordinates the branch locus is $C_a$ for some $a_1,a_2,a_3 \in \mathbb C$.
\end{theorem}
\begin{proof}
 The surface $X$ is a double cover of $\mathbb P_1 \times \mathbb P_1$ branched along a smooth $H$-invariant curve of bidegree (4,4). The invariant (4,4)-curves $C_b$ and $C_0$ discussed above are seen to be singular at $([1:0],[1:0])$ or $([1:0],[0:1])$. 
\end{proof}
Note that the general curve $C_a$ is smooth. We obtain a 2-dimensional family $\{C_a\}$ of smooth branch curves and a corresponding family of K3-surfaces $\{X_{C_a}\}$.
\subsubsection{Curves of bidegree $(4,4)$ with eight singular points}
It remains to consider the case 2. of the classification.
Our aim is to find an example of a K3-surface $X$ such that $X/\sigma = Y $ has a nontrivial Mori reduction $M: Y  \to \mathbb P_1 \times \mathbb P_1= Z$ contracting a single $H$-orbit of Mori fibers consisting of precisely 8 curves. In this case the branch locus $B \subset Y$ is mapped to a singular $(4,4)$-curve $C= M(B)$ in $Z$. The curve $C$ is irreducible and has precisely 8 singular points along a single $H$-orbit in $Z$. 

As we have noted above, many of the curves $C_a,C_b,C_0$ are seen to be singular at  $([1:0],[1:0])$ or $([1:0],[0:1])$. Since both points lie in $H$-orbits of length two, these curves are not candidates for our construction. This argument excludes the curves $C_b, C_0$ and $C_a$ if $a_1 = 0$ or $a_2 = 0$.

For $C_a$ with $a_3=0$ one checks that $C_a$ has singular points if and only if $a_1 = -a_2$, i.e., if $C_a$ is reducible. It therefore remains to consider curves $C_a$ where all coefficients $a_i \neq 0$. We choose $a_3=1$.
\begin{lemma}
If $a_i\neq 0$ for $i=1,2,3$, then $C_a$ is irreducible.
\end{lemma}
\begin{proof}[Sketch of proof]
First note that $C_a$ does not pass through $([1:0],[1:0])$ or $([1:0],[0:1])$. Therefore, possible singularities or points of intersection of irreducible components come in orbits of length eight.
Assume that $C_a$ is reducible, consider the decomposition into irreducible components and the $H$-action on it. A curve of type $(n,0)$ is always reducible for $n>1$ and therefore does not occur in the decomposition.

If $C_a$ contains a (2,2)-curve $C_a^{(2,2)}$, then the $H$-orbit of $C_a^{(2,2)}$ has length $\leq2 $ and $C_a^{(2,2)}$ is stable with respect to the subgroup $H' = \langle c \rangle $ of $H$. All $c$-semi-invariants of bidegree (2,2) are, however, reducible. Similary, all $c$-semi-invariants of bidegree (1,2) or (2,1) are reducible an therefore $C$ does not have a curve of this type as an irreducible component. 

The curve $C_a$ is not the union of a (1,3)- and a (3,1)-curve, since 
their intersection number is 10 and contradicts invariance. Similarly one excludes the union of a (1,1) and a (3,3)-curve. 

If $C_a$ is a union of (1,1) or (1,0) and (0,1)-curves, one checks by direct computation that the requirement that $C_a$ is $H$-invariant gives strong restrictions and finds that in all cases at least one coefficient $a_i$ has to be zero. 
\end{proof}
One possible choice of an orbit of length eight is given by the orbit through a $\tau$-fixed point $p_\tau = ([1:1],[\pm \sqrt{i}:1])$. One checks that $p_\tau \in C_a$ for any choice of $a_i$. However, if we want $C_a$ to be singular in $p_\tau$, then $a_2=0$. It then follows that $C_a$ is singular at points outside $H .p_\tau$. It has more than eight singular points  and is therefore reducible.

All other orbits of length eight are given by orbits through $g$-fixed points $p_x= ([1:x],[1:x])$ for $x \neq 0$. One can choose coefficients $a_i(x)$ such that $C_{a(x)}$ is singular at $p_x$ if and only if $x^8 \neq 1$. If the curve $C_{a(x)}$ is irreducible, then it has precisely eight singular points $H.p_x$ of multiplicity 2 (cusps or nodes) and the double cover of $\mathbb P_1 \times \mathbb P_1$ branched along $C_{a(x)}$ is a singular K3-surface with precisely eight singular points. We obtain a diagram
\[
\begin{xymatrix}{
X_\mathrm {sing}\ar[d]^{2:1} & X \ar[d]^{2:1} \ar[l]^{\text{desing.}}\\
C_{(4,4)} \subset \mathbb P_1 \times \mathbb P_1 & \ar[l]^<<<<<{M} Y \supset B}
\end{xymatrix}
\]
If $p_x$ is a node in $C_{a(x)}$, then the corresponding singularity of $X_\mathrm{sing}$ is resolved by a single blow-up. The (-2)-curve in $X$ obtained from this desingularization is a double cover of a (-1)-curve in $Y$ meeting $B$ in two points.

If $p_x$ is a cusp in $C_{a(x)}$, then the corresponding singularity of $X_\mathrm{sing}$ is resolved by two blow-ups. The union of the two intersecting (-2)-curves in $X$ obtained from this desingularization  is a double cover of a (-1)-curve in $Y$ tangent to $B$ is one point.

The information determining whether $p_x$ is a cusp or a node is encoded in the rank of the Hessian of the equation of $C_{a(x)}$ at $p_x$. The condition that this rank is one gives a nontrivial polynomial condition. For a general irreducible member of the family $\{C_{a(x)} \, | \, x\neq 0, \,  x^8 \neq 1 \}$ the singularities of $C_{a(x)}$ are nodes.

We let $q$ be the polynomial in $x$ that vanishes if and only if the rank of the Hessian of $C_{a(x)}$ at $p_x$ is one. It has degree 24, but 16 of its solutions give rise to reducible curves $C_{a(x)}$. The remaining eight solution give rise to four different irreducible curves. These are identified by the action of the normalizer of $H$ in $\mathrm{Aut}(\mathbb P_1 \times \mathbb P_1)$ and therefore define equivalent K3-surfaces.

We summarize the discussion in the following main classification theorem.
\begin{samepage}
\begin{theorem}\label{classiA6}
Let $X$ be a K3-surface with an effective action of the group $G$ such that $\mathrm{Fix}_X(h \sigma) = \emptyset$. Then $X$ is an element of one the following families of K3-surfaces:
\begin{enumerate}
\item 
the two-dimensional family  $\{X_{C_a}\}$ for $C_a$ smooth,
\item
the one-dimensional family of minimal desingularization of double covers of $\mathbb P_1 \times \mathbb P_1$ branched along curves in $\{C_{a(x)} \, | \, x\neq 0, \,  x^8 \neq 1 \}$. The general curve $C_{a(x)}$ has precisely eight nodes along an $H$-orbit. Up to natural equivalence there is a unique curve $C_{a(x)}$ with eight cusps along an $H$-orbit.
\item
the trivial family consisting only of the minimal desingularization of the double cover of $\mathbb P_1 \times \mathbb P_1$ branched along the curve $C_a = \{f_1-f_2=0\}$ where $a_1 =1, a_2 =-1, a_3=0$.
\end{enumerate}
\end{theorem}
\end{samepage}
\begin{corollary}
 Let $X$ be a K3-surface with an effective action of the group $\tilde A_6$. If $\mathrm{Fix}_X(h \sigma) = \emptyset$, then $X$ is an element of one the families 1. -3. above. If $\mathrm{Fix}_X(h \sigma) \neq \emptyset$, then $X$ is $A_6$-equivariantly isomorphic to the Valentiner surface.
\end{corollary}
\section{Summary and outlook}
Recall that our starting point was the description of K3-surfaces with $\tilde A_6$-symmetry. Using the group structure of $\tilde A_6$ we have divided the problem into two possible cases corresponding to the question whether $\mathrm{Fix}_X(h\sigma)$ is empty or not. If it is nonempty, the K3-surface with $\tilde A_6$-symmetry is the Valentiner surface discussed in Section \ref{A6Valentiner}. 
If is is empty, our discussion in the previous sections has reduced the problem to finding the $\tilde A_6$-surface in the families of surfaces $X_{C_a}$ with $D_{16}$-symmetry.

It is known that a K3-surface with $\tilde A_6$-symmetry has maximal Picard rank 20. This follows from a criterion due to Mukai (cf. \cite{mukai}) and is explicitely shown in \cite{KOZLeech}. 

All surfaces $X_{C_a}$ for $C_a \subset \mathbb P_1 \times \mathbb P_1$ a (4,4)-curve are elliptic since the natural fibration of $\mathbb P_1 \times \mathbb P_1$ induces an elliptic fibration on the double cover (or is desingularization). 

A possible approach for finding the $\tilde A_6$-example inside our families is to find those surfaces with maximal Picard number by studying the elliptic fibration.
It would be desirable to apply the following formula for the Picard rank of an elliptic surface $f: X \to \mathbb P_1$ with a section (cf. \cite{shiodainose}):
\[
 \rho(X) = 2 + \mathrm{rank}(MW_f) + \sum_i (m_i-1) 
\]
where the sum is taken over all singular fibers, $m_i$ denotes the number of irreducible components of the singular fiber and $\mathrm{rank}(MW_f)$ is the rank of the Mordell-Weil group of sections of $f$. The number two in the formula is the dimension of the hyperbolic lattice spanned by a general fiber and the section. 

First, one has to ensure that the fibration under consideration has a section. One approach to find sections is to consider the quotient $q:\mathbb P_1 \times \mathbb P_1 \to \mathbb P_2$ and the image of the curve $C_a$ inside $\mathbb P_2$. If we find an appropiate bitangent to $q(C_a)$ such that its preimage in $\mathbb P_1 \times \mathbb P_1$ is everywhere tangent to $C_a$, then its preimage in the double cover of $\mathbb P_1 \times \mathbb P_1$ is reducible and both its components define sections of the elliptic fibration. For $C_a$ the curve with eight nodes the existence of a section (two sections) follows from an application of the Pl\"ucker formula to the curve $q(C_a)$ with 3 cusps and its dual curve. 

As a next step, one wishes to understand the singular fibers of the elliptic fibrations. Singular fibers occur whenever the branch curve $C_a$ intersects a fiber $F$ of the $\mathbb P_1 \times \mathbb P_1$ in less than four points. Depending on the nature of intersection $F \cap C_a$ one can describe the corresponding singular fiber of the elliptic fibration. For $C_a$ the curve with eight cusps one finds precisely eight singular fibres of type $I_3$, i.e., three rational curves forming a closed cycle. In particular, the contribution of all singular fibres $\sum_i (m_i-1)$ in the formula above is 16. In the case where $C_a$ is smooth or has eight nodes, this contribution is less. 

In order to determine the number $ \rho(X_{C_a})$ it is neccesary to either understand the Mordell-Weil group and its $\mathrm{rank}(MW_f)$ or to find curves which give additional contribution to $\mathrm{Pic}(X_{C_a}$ not included in $2 + \sum_i (m_i-1)$. 

In conclusion, the method of equivariant Mori reduction applied to quotients $X/\sigma$ yields an explicit description of a families of K3-surfaces with $D_{16} \times \langle \sigma \rangle$-symmetry and by construction, the K3-surface with $\tilde A_6$-symmetry is contained in one of these families.  It remains to find criteria to characterize this particular surface inside this family. The possible approach by understanding the function 
\[
 a \mapsto \rho( X_{C_a}) 
\]
using the elliptic structure of $X_{C_a}$ requires a detailed analysis of the Mordell-Weil group. 
%
%
%
%
\appendix
\chapter{Actions of certain Mukai groups on projective space}
In this appendix, we derive the unique action of the group $N_{72}$ on $\mathbb P_3$ and the unique action of $M_9$ on $\mathbb P_2$ in the context of Sections \ref{N72} and \ref{M9}. We consider the homomorphism $\mathrm{SL}_n(\mathbb C) \to \mathrm{PSL}_n(\mathbb C)$ and determine preimages $\tilde g \in \mathrm{SL}_n(\mathbb C)$ of the generators $g \in G  \subset \mathrm{PSL}_n(\mathbb C)$. 
Our considerations benefit from fact that both actions are induced by symplectic actions of the corresponding group on a K3-surface $X$. 
\section{The action of $N_{72}$ on $\mathbb P_3$}\label{N72appendix}
One can calculate explicitly the realization of the $N_{72}$-action on $\mathbb P_3$ by using the decomposition $C_3^2 \rtimes D_8$ where $D_8 = C_2\ltimes (C_2 \times C_2) = \mathrm{Aut}(C_3^2)$. For each generator of $N_{72}$ we will specify the corresponding element in $\mathrm{SL}_4(\mathbb C)$. We denote the center of $\mathrm{SL}_4(\mathbb C)$ by $Z$.
Recall that the action of $D_8 = C_2\ltimes (C_2 \times C_2) =
\langle \alpha \rangle  \ltimes ( \langle \beta \rangle \times \langle \gamma \rangle)$ on $C_3 \times C_3$ is given by
\[
\alpha(a,b) = (b,a), \quad \beta(a,b)=(a^2,b), \quad \gamma(a,b) = (a,b^2).
\]

In suitably chosen coordinates the generator $a$ of $C_3^2$ can be represented as 
\begin{align*}
\tilde a=
 \begin{pmatrix}
\xi & 0 & 0 & 0  \\
0 & \xi ^2& 0 & 0\\
0&0& 1 & 0\\
0&0&0&1
 \end{pmatrix}
\end{align*}
where $\xi$ is a third root of unity.
Next we wish to specify $\gamma$ in $\mathrm{SL}_4(\mathbb C)$. We know that $a \gamma = \gamma a$, i.e.,  $\tilde a \tilde \gamma \tilde a^{-1} \tilde \gamma^{-1} \in Z$, and $\gamma$ is seen to be of the form 
\begin{align*}
\tilde \gamma =
 \begin{pmatrix}
* & 0 & 0 & 0  \\
0 & * & 0 & 0\\
0&0& * & *\\
0&0& * & *
 \end{pmatrix}
\end{align*}
where $*$ denotes a nonzero matrix entry.
Since $a$ and $b$ commute in $N_{72}$, we know that $\tilde a \tilde b \tilde a^{-1} \tilde b^{-1} \in Z$  and 
\begin{align*}
\tilde b =
 \begin{pmatrix}
* & 0 & 0 & 0  \\
0 & * & 0 & 0\\
0&0& * & *\\
0&0& * & *
 \end{pmatrix}
\end{align*}
Since $\gamma$ acts on $ b$ by $\gamma b \gamma = b^{-1} = b^2$, it follows that 
\begin{align*}
\tilde b =
 \begin{pmatrix}
1 & 0 & 0 & 0  \\
0 & 1 & 0 & 0\\
0&0& * & *\\
0&0& * & *
 \end{pmatrix}.
\end{align*}
We apply a change of coordinates affecting only the lower $(2 \times 2)$-block of $b$ and therefore not affecting the shape of $a$ auch that
\begin{align*}
\tilde b= 
 \begin{pmatrix}
1 & 0 & 0 & 0\\
0 & 1 & 0 & 0\\
0 & 0 & \xi & 0 \\
0 & 0 & 0 & \xi^2
 \end{pmatrix}.
\end{align*}
It follows that $\alpha$ interchanges the two $(2 \times 2)$-blocks of the matrices $a$ and $b$ and 
\begin{align*}
\tilde \alpha= 
 \begin{pmatrix}
0 & 0 & 1 & 0\\
0 & 0 & 0 & 1\\
1 & 0 & 0 & 0 \\
0 & 1 & 0 & 0
 \end{pmatrix}.
\end{align*}
Finally, $\gamma$ and $\beta$ can be put into the form
\begin{align*}
\tilde \gamma =
 \begin{pmatrix}
1 & 0 & 0 & 0  \\
0 & 1 & 0 & 0\\
0&0& 0 & 1\\
0&0& 1 & 0
 \end{pmatrix}, \quad
\tilde \beta =
 \begin{pmatrix}
0 & 1 & 0 & 0  \\
1 & 0 & 0 & 0\\
0&0& 1 & 0\\
0&0& 0 & 1
 \end{pmatrix}.
\end{align*}
\subsection{Invariant quadrics and cubics}
Let $f \in \mathbb C_2[x_1:x_2:x_3:x_4]$ be a semi-invariant homogeneous polynomial of degree two, 
\[
 f = \sum_i a_i x_1^2 + \sum _{i\neq j} b_{ij} x_ix_j.
\]
If $a_i \neq 0$ for some $1 \leq i \leq 4$, then semi-invariance with respect to the transformations $\alpha, \beta, \gamma$ yields
$a_1 = a_2 = a_3 = a_4$. It follows that $f$ is not semi-invariant with respect to $a$. 

If $b_{13} \neq 0$, then semi-invariance with respect to the transformations $\alpha, \beta, \gamma$ yields $b_{13} = b_{23} = b_{24} = b_{14}$. As above, the polynomial $f$ is not semi-invariant with respect to $a$.

Therefore, if $f$ is semi-invariant, then $a_i = b_{13} = b_{23} = b_{24} = b_{14} =0$ and $b_{12} = b_{34}$. In particular, all degree two semi-invariants are in fact invariant.
There is a unique $N_{72}$-invariant quadric hypersurface in $\mathbb P_3$ given by the equation $x_1 x_2 + x_3 x_4$ .

Analogous considerations show that a semi-invariant polynomial of degree three is a multiple of $f_\mathrm{Fermat} = x_1^3 +x_2^3 + x_3^3 +x_4^3$ and the Fermat cubic $\{f_\mathrm{Fermat} =0\}$ is seen to be the unique $N_{72}$-invariant cubic hypersurface in $\mathbb P_3$.
\section{The action of $M_9$ on $\mathbb P_2$}\label{M9 on P2}
We consider the decompostion of $M_9 = (C_3 \times C_3) \rtimes Q_8$. The generators of $ C_3 \times C_3$ are denoted $a$ and $b$ and the generators of $Q_8$ are denoted by $I, J, K$. Recall $ I^2 = J^2 = K^2 = IJK = -1$. We choose the factorization of $C_3 \times C_3$ such that $-1$  acts as
\[
 (-1)a(-1) = a^2, \quad  (-1)b(-1) = b^2.
\]
Furthermore, $I a (-I) = b$ and  $Ja(-J) = b^2 a$.

We repeatedly use the fact that the action of $M_9$ is induced by a symplectic action of $M_9$ on a K3-surface $X$ which is a double cover of $\mathbb P_2$. 

We begin by fixing a representation of $a$. Since $a$ may not have a positive dimensional set of fixed points in $\mathbb P_2$, it follows that in appropriately chosen coordinates 
\[
\tilde a=
 \begin{pmatrix}
1 & 0 & 0 \\
0 & \xi & 0\\
0&0& \xi^2
 \end{pmatrix}, 
\]
where $\xi$ is third root of unity. 

As a next step, we want to specify a representation of $b$ inside $\mathrm{SL}_3(\mathbb C)$. Since $a$ and $b$ commute in $\mathrm{PSL}_3(\mathbb C)$, we know that
\[
\tilde a \tilde b \tilde a^{-1} \tilde b^{-1} = \xi^k \mathrm{id}_{\mathbb C^3}
\]
for $k \in \{0, 1,2\}$.
Note that $\tilde b$ is not diagonal in the coordinates chosen above since this would 
 give rise to $C_3^2 $-fixed points in $\mathbb P_2$. As these correspond to $C_3^2$-fixed points on the double cover $X \to Y$ and a symplectic action of $C_3^2 \nless \mathrm{SL}_2(\mathbb C)$ on a K3-surface does not admit fixed points, this is a contradiction. An explicit calculation yields
\begin{align*}
\tilde b = \tilde b_1=
 \begin{pmatrix}
0 & 0 & * \\
* & 0 & 0\\
0 & * & 0 
 \end{pmatrix} \quad \text{or} \quad
\tilde b = \tilde b_2= 
 \begin{pmatrix}
0 & * & 0 \\
0 & 0 & *\\
* & 0 & 0 
 \end{pmatrix}.
\end{align*}
We can introduce a change of coordinates commuting with $\tilde a$ such that
\begin{align*}
\tilde b = \tilde b_1 =
 \begin{pmatrix}
0 & 0 & 1 \\
1 & 0 & 0\\
0 & 1 & 0 
 \end{pmatrix} \quad \text{or} \quad
\tilde b = \tilde b_2= 
 \begin{pmatrix}
0 & 1 & 0 \\
0 & 0 & 1\\
1&0&0 
 \end{pmatrix}.
\end{align*}
Since $\tilde b_1 = \tilde b_2^2$, the two choices above correspond to choices of generators $b$ and $b^2$ of $\langle b \rangle$ and are therefore equivalent. In the following we fix the second choice of $b$.
A direct computation yields that the element $-1$ must be represented in the form
\begin{align*}
 \begin{pmatrix}
* & 0 & 0 \\
0 & 0 & *\\
0 & * & 0 
 \end{pmatrix} \quad \text{or} \quad
 \begin{pmatrix}
0 & * & 0 \\
* & 0 & 0\\
0 & 0 & * 
 \end{pmatrix} \quad \text{or} \quad
 \begin{pmatrix}
0 & 0 & * \\
0 & * & 0\\
* & 0 & 0 
 \end{pmatrix}.
\end{align*}
After reordering the coordinates, we can assume that 
\begin{align*}
\widetilde{-1} = 
 \begin{pmatrix}
* & 0 & 0 \\
0 & 0 & *\\
0 & * & 0 
 \end{pmatrix}.
\end{align*}
The relation $(-1)b(-1) = b^2$ yields
\begin{align*}
\widetilde{-1} = 
 \begin{pmatrix}
-1 & 0 & 0 \\
0 & 0 & \eta\\
0 & \eta^2 & 0 
 \end{pmatrix}.
\end{align*}
for some third root of unity $\eta$. The element $I$ fulfills $I a (-I) = b$ and, using the representation of $a$ and $b$ given above, we conclude 
\begin{align*}
\widetilde{I} = \frac{1}{\xi -\xi ^2}
 \begin{pmatrix}
1 & 1 & 1\\
\zeta^2 & \zeta^2 \xi &  \zeta^2 \xi^2 \\
\zeta  & \zeta \xi^2 & \zeta \xi
 \end{pmatrix}
\end{align*}
for some third root of unity $\zeta$. Now $I^2 = -1$ implies $\zeta = 1$ and $\eta =1$. Analogous considerations yield the following shape of $J$:
\begin{align*}
\widetilde{J} = \frac{1}{\xi -\xi ^2}
 \begin{pmatrix}
1 & \xi & \xi\\
\xi^2 &  \xi &   \xi^2 \\
\xi^2 &  \xi^2 & \xi
 \end{pmatrix}.
\end{align*}
In appropiately chosen coordinates the action on $M_9$ is precisely of the type claimed in Section \ref{M9}. 
\nocite{atlas, jamesliebeck}
\newpage
\addcontentsline{toc}{chapter}{Bibliography}
\newcommand{\etalchar}[1]{$^{#1}$}
\providecommand{\bysame}{\leavevmode\hbox to3em{\hrulefill}\thinspace}
\end{document}